\documentclass[a4paper]{amsart}
\usepackage{amsmath}
\usepackage{amsthm}
\usepackage{amssymb}
\usepackage{amscd}
\usepackage{latexsym}
\usepackage{amsfonts}
\usepackage{amsmath}

\newcommand\veclam{{\vec \lambda}}

\newcommand\Vlam{{{\Cal V}_{\vec \lambda}}}

\newcommand\Vdaglam{{\Cal V}_{\vec \lambda}^{\dagger}}

\newcommand\Hdaglam{{\Cal H}_{\vec \lambda}^{\dagger}}

\newcommand\Hlam{{{\Cal H}_{\vec \lambda}}}

\newcommand\Ob{{{\Cal O}_{\Cal B}}}

\newcommand\Oc{{{\Cal O}_{\Cal C}}}

\newcommand\sumjn{\sum_{j=1}^N}

\newcommand\bojn{\bigoplus_{j=1}^N}

\newcommand\Res{\operatornamewithlimits{Res}}
\newcommand\Dim{\operatornamewithlimits{dim}}

\newcommand\bolc{{\mathbf C}}

\newcommand\goX{{\goth X}}

\newcommand\gog{{\goth g}}

\newcommand\goF{{\goth F}}

\newcommand\Hom{{\hbox{\rm Hom}}}

\newcommand\Bosonnormalord{\,\lower.8ex \hbox{$\circ$} \llap{\raise.8ex\hbox{$\circ$}} \,}

\newtheorem{proposition}{Proposition}[section]
\newtheorem{lemma}{Lemma}[section]
\newtheorem{theorem}{Theorem}[section]
\newtheorem{corollary}{Corollary}[section]

\theoremstyle{definition}
\newtheorem{definition}{Definition}[section]
\newtheorem{remark}{Remark}[section]

\newcommand{\eproof}{\begin{flushright} $\square$ \end{flushright}}

\newcommand{\good}{good\mbox{ }}
\newcommand{\goodp}{good.}
\newcommand{\stable}{stable\mbox{ }}

\newcommand{\fm}{saturated\mbox{ }}
\newcommand{\fmp}{saturated.\mbox{ }}
\newcommand{\fmc}{saturated,\mbox{ }}
\newcommand{\abelian}{abelian\mbox{ }}

\newcommand{\e}[1]{\mathbf #1}
\newcommand{\goth}[1]{\mathfrak #1}

\renewcommand{\varinjlim}{\mathop{\fam0 \stackrel{\textstyle Lim}{\textstyle \longleftarrow}}}

\newcommand{\Vdagtlam}{{\mathcal V}_{\e {{\tilde \lambda}}}^{\dagger}}
\newcommand{\VdagSlam}{{\mathcal V}_{S(\vec \lambda)}^{\dagger}}
\newcommand{\HSlam}{{\mathcal H}_{S(\vec \lambda)}}

\newcommand{\Spofvlam}{\Vdaglam}
\newcommand{\Spofvlamp}{{\mathcal V}_{\e \lambda'}^{\dagger}}
\newcommand{\Vdagmulam}{{\mathcal V}_{\mu, \mu^{\dagger}, \vec{\lambda}}^{\dagger}}
\newcommand{\ab}{\mathop{\fam0 ab}\nolimits}
\newcommand{\Vdagab}{{\mathcal V}_{\ab}^{\dagger}}

\newcommand{\Hdagtlam}{{\mathcal H}_{\e {{\tilde \lambda}}}^{\dagger}}
\newcommand{\Hlamone}{{{\mathcal H}_{{\vec \lambda}_1}}}
\newcommand{\Hlamtwo}{{{\mathcal H}_{{\vec \lambda}_2}}}
\newcommand{\mC}{{\mathcal C}}

\newcommand{\tgoF}{{\tilde {\mathfrak F}}}

\newcommand{\im}{\mathop{\fam0 Im}\nolimits}
\newcommand{\spane}{\mathop{\fam0 Span}\nolimits}
\newcommand{\re}{\mathop{\fam0 Re}\nolimits}

\newcommand{\ad}{\mathop{\fam0 ad}\nolimits}

\newcommand{\id}{\mathop{\fam0 Id}\nolimits}
\newcommand{\Lie}[1]{\mbox{\sf #1}}

\newcommand{\Aut}{\mathop{\fam0 Aut}\nolimits}

\newcommand{\R}[1]{\mathop{\fam0 {\mathbb R}^{#1}}\nolimits}
\newcommand{\RPP}{\mathop{\fam0 {\mathbb R}_+^{P}}\nolimits}
\newcommand{\C}[1]{\mathop{\fam0 {\mathbb C}^{#1}}\nolimits}

\newcommand{\Z}{{\mathbb Z}}
\newcommand{\tD}{{\tilde D}}

\newcommand{\ra}{\mathop{\fam0 \rightarrow}\nolimits}

\newcommand{\tnabla}{\mathop{\fam0 {\tilde \nabla}}\nolimits}
\newcommand{\tc}{\mathop{\fam0 {\tilde c}}\nolimits}
\newcommand{\cT}{\mathop{\fam0 {\mathcal T}}\nolimits}
\newcommand{\ctT}{\mathop{\fam0 {\mathcal T}^{(r)}}\nolimits}
\newcommand{\tPsi}{\mathop{\fam0 {\tilde \Psi}}\nolimits}
\def\l{\lambda}\def\m{\mu}
\def\L{\Lambda}

\newcommand{\fg}{\mathop{\fam0 {\mathfrak g}}\nolimits}
\newcommand{\cv}{\mathop{\fam0 {\upsilon}}\nolimits}
\newcommand{\Si}{\Sigma}
\newcommand{\Sib}{{\mathbf \Sigma}}

\def\vQ{\vec Q}
\def\vtQ{\vec {\widetilde Q}}
\def\tx{\tilde x}
\def\vs{\vec s}
\def\vh{\vec h}
\def\vth{\vec {\widetilde h}}
\def\vxi{\vec \xi}
\def\vtxi{\vec {\widetilde \xi}}
\def\veta{\vec \eta}
\def\vteta{\vec {\widetilde \eta}}
\def\veclam{{\vec \lambda}}
\def\Vlam{{{\mathcal V}_{\vec \lambda}}}
\def\Vdag{{\mathcal V}^{\dagger}}
\def\Vdaglam{{\mathcal V}^{\dagger}}
\def\Vdaglam{{\mathcal V}_{\vec \lambda}^{\dagger}}

\def\Vdaglamone{{\mathcal V}_{\vec \lambda_1}^{\dagger}}

\def\Vdaglamtwo{{\mathcal V}_{\vec \lambda_2}^{\dagger}}

\def\Spofv{{\mathcal V}^{\dagger}}
\def\Spofvlam{{\mathcal V}_{\l}^{\dagger}}
\def\Spofvlambar{{\mathcal V}_{{\bar \l}}^{\dagger}}
\def\Spofvlamone{{\mathcal V}_{\l_1}^{\dagger}}
\def\Spofvlami{{\mathcal V}_{\l_i}^{\dagger}}
\def\Spofvlamp{{\mathcal V}_{\l'}^{\dagger}}
\def\Spofvlampone{{\mathcal V}_{\l'_1}^{\dagger}}
\def\Spofvlamtwo{{\mathcal V}_{\l_2}^{\dagger}}
\def\Spofvlamptwo{{\mathcal V}_{\l'_2}^{\dagger}}
\def\Hdaglam{{\mathcal H}_{\vec \lambda}^{\dagger}}
\def\Hdaglamone{{\mathcal H}_{\vec \lambda_1}^{\dagger}}

\def\Hdaglamtwo{{\mathcal H}_{\vec \lambda_2}^{\dagger}}
\def\Hlam{{{\mathcal H}_{\vec \lambda}}}

\def\Ob{{{\mathcal O}_{\mathcal B}}}
\def\Oc{{{\mathcal O}_{\mathcal C}}}
\def\Ocal{{\mathcal O}}
\def\sumjn{\sum_{j=1}^N}
\def\bojn{\bigoplus_{j=1}^N}
\def\Res{\operatornamewithlimits{Res}}
\def\ad{\operatorname{ad}}

\def\bolc{{\mathbb C}}

\def\goX{{\mathfrak X}}
\def\gog{{\mathfrak g}}
\def\goF{{\mathfrak F}}
\def\tgoF{ {\tilde {\mathfrak F}}}
\def\:{:}

\def\Hom{\operatorname{Hom}}

\def\Bosonnormalordconstruction#1{\vcenter{\hbox{\ooalign{%
\raise.8ex\hbox{$#1\circ$}\crcr\lower.8ex\hbox{$#1\circ$}}}}}

\def\hZ{{{\mathbf Z}_h}}
\def\cWd{{\mathcal W}^{\dagger}}
\def\bC{{\mathbf C}}
\def\bZ{{\mathbf Z}}
\def\cW{{\mathcal W}}
\def\cF{{\mathcal F}}
\def\cFd{{\mathcal F}^\dagger}
\def\cWd{{\mathcal W}^{\dagger}}
\def\ovpsi{{\overline{\psi}}}
\def\normalord{\,\lower.8ex \hbox{$\cdot$} \llap{\raise.8ex\hbox{$\cdot$}} \,}

\def\gF{\mathfrak{F}}
\def\gX{\mathfrak{X}}
\def\oe{{\overline{e}}}
\def\cO{{{\mathcal O}}}
\def\cV{{{\mathcal V}}}
\def\cD{{\mathcal D}}
\def\cC{{\mathcal C}}
\def\cB{{\mathcal B}}
\def\cL{{\mathcal L}}

\begin{document}

\title[Geometric construction of Modular functors]{Geometric construction of modular functors
from Conformal Field Theory}

\author{J{\o}rgen Ellegaard Andersen}
\address{Department of Mathematics\\
        University of Aarhus\\
        DK-8000, Denmark}
\email{andersen{\@@}imf.au.dk}

\author{Kenji Ueno}
\address{Department of Mathematics\\
        Faculty of Science, Kyoto University\\
        Kyoto, 606-01 Japan}
\email{ueno{\@@}math.kyoto-u.ac.jp} %\dedicatory{In preparation.}
%\date{October 23, 1995}

\thanks{This research was conducted
 partly by the first author for the Clay Mathematics Institute at University of California,
 Berkeley and for MaPhySto --
A Network for Mathematical Physics and Stochastics, funded by The
Danish National Research Foundation. The second author is partially supported by Grant in Aid
for Scientific Research {\sc no}. 14102001 of JSPS}

\begin{abstract} We give a geometric construct of a modular
functor for any simple Lie-algebra and any level by twisting the
constructions in \cite{TUY} and \cite{Ue2} by a certain fractional
power of the abelian theory first considered in \cite{KNTY} and
further studied in \cite{AU1}.
\end{abstract}
%\newpage
\maketitle

\tableofcontents

%\addtolength{\textheight}{-10cm}
%\newpage
%\addtolength{\textheight}{-10cm}

\section{Introduction}

This is the second paper in a series of three papers (\cite{AU1} and \cite{AU3})
in which we provide
a geometric construction of modular functors and topological
quantum field theories from conformal field theory building on the
constructions in \cite{TUY} and \cite{Ue2} and \cite{KNTY}. In this
paper we provide the geometric construction of a modular
functor $V^{\fg}_{\ell}$ for each simple Lie algebra $\fg$ and a
positive integer $\ell$ (the {\em level}).  In our third paper \cite{AU3}
in this series we give
an explicit isomorphism of the modular functor underlying the
Reshetikhin-Turaev TQFT for $U_q(\Lie{sl}(n))$ and the one
constructed in this paper for the Lie algebra $\Lie{sl}(n)$. This
uses the Skein theory approach to the Reshetikhin-Turaev TQFT of
Blanchet, Habegger, Masbaum and Vogel \cite{BHMV1}, \cite{BHMV2}
and \cite{Bl1}. In particular we use Blanchet's \cite{Bl1} constructions of the
Hecke-category and its associated modular tensor categories. This
construction is really a generation of the BHMV-construction of the
$U_q(sl_2(\C{} ))$-Reshetikhin-Turaev TQFT \cite{BHMV2} to the $U_q(sl_n(\C{} ))$-case.
As a consequence of this, we construct a duality and a unitary structure on our modular
functor in the case of $\fg= \Lie{sl}(n)$.

By a very general
construction any modular functor with duality induces a
topological quantum field theory in dimension $2+1$ by the work of
Kontsevich \cite{Kont1} and Walker \cite{Walker} and also Grove
\cite{Grove}. By applying this to the modular functor
$V^{\fg}_{\ell}$, for $\fg= \Lie{sl}(n)$, we get a TQFT for each
$\ell$. We also prove in \cite{AU3} that this TQFT is isomorphic
to the Reshetikhin-Turaev TQFT for $U_q(\Lie{sl}(n))$
at level $\ell$.

Let us now describe our construction.
Fix a simple Lie algebra $\fg$ and normalize the invariant inner
product on it by requiring the highest root to have length squared
equal to $2$. Let $\ell$ be a positive integer and consider the
finite (label) set $P_{\ell}$ of integrable highest weight
representations at level $\ell$ of the affine Lie algebra of
$\fg$. By the usual highest weight vector representations, this
finite set $P_{\ell}$ is naturally identified with a subset of the
dominant integrable weights of $\fg$ (see formula
(\ref{labelset})).

The main idea is to construct a modular functor $V^{\fg}_{\ell}$
by associating to each labeled marked surface the {\em space of
vacua} using the given labels for the Lie algebra $\fg$ at level
$\ell$ for some complex structure on the marked surface, as
defined in \cite{TUY} \cite{Ue} and \cite{Ue2}. In order to make this
construction independent of the complex structure it must be
understood in terms of bundles with connections over the hole of
Teichm\"{u}ller space of the surface, relying on parallel transport to
provide the required identifications between the different spaces
of vacua. The consistencies of these identifications translates to
flatness requirements on these connections. However, the {\em
sheaf of vacua} construction in \cite{TUY} and \cite{Ue2} gives a
bundles with a connections, which is only projectively flat, over
Teichm\"{u}ller space of the surface. By tensoring this bundle with a
line bundle with a connection with the opposite curvature, we get
a flat bundle over Teichm\"{u}ller space and the vector space we
associate to the labeled marked surface is the vector space of
covariant constant sections of this resulting bundle. This line
bundle is constructed as a fractional power of a certain rank $1$
\abelian sheaf of vacua, which we considered in the first paper in
this series \cite{AU1} from the same point of view as \cite{TUY}
and \cite{Ue2}. The extraction of this fractional power brings in
central extensions of mapping classes as the natural morphisms on
which the resulting functor is defined.

The construction and properties of this flat bundle primarily rely
on the complex algebraic constructions and results of \cite{TUY}
and \cite{Ue2} on the the sheaf of vacua construction yielding a
conformal field theory for each simple Lie algebra $\fg$ and level
$\ell$. For the $1$-dimensional correction theory, we draw on the
work \cite{AU1}, which in turn relies on \cite{KNTY}.

The definition of the functor $V^{\fg}_{\ell}$ requires only considerations
of smooth families of Riemann surfaces (with some extra structure
which is specified in section \ref{New1}) over smooth complex manifolds. This is
the setting for section \ref{New1} through to section \ref{construction}.
However, in order to define the glueing isomorphism, which a modular functor
is required to have, we need to discuss certain very simple families of stable curves,
which contains so call nodal curves. These are described and considered
in section \ref{shofvacandglue} and in the Appendix to this paper.

The paper is organized as follows.
In section \ref{AxiomsMF} we give the axioms for a modular
functor. We introduce the notion of a {\em marked} surface, which
is a closed smooth oriented surface, with a finite subset of
points with projective tangent vectors and a Lagrangian subspace of
the first integer homology of the surface. These form a category
on which there is the operation of disjoint union and the
operation of orientation reversal. There is also the process of
glueing on this category. If we have a finite set, we can label
the finite set of points on a marked surface by elements from this
finite {\em label} set and get the category of labeled marked
surfaces. A modular functor based on some finite label set, is a
functor from the category of labeled marked surfaces to the
category of finite dimensional complex vector spaces, which takes
the disjoint union operation to the tensor product operation and
which takes the glueing process to a certain direct sum
construction, such that some compatibility holds, as described
in details in Definition \ref{DefMF}. A modular functor is said to
be with duality if further the operation of orientation reversal
is taken to the operation of taking the dual vector space.

In sections \ref{New1} to \ref{New5} we describe in detail how any
simple Lie algebra and a level $\ell$, via  the sheaf of vacua
constructions in \cite{Ue2} yields a holomorphic vector bundles
with a projectively flat connection over Teichm\"{u}ller spaces of
pointed surfaces equipped with symplectic basis of the first
homology. In sections \ref{New3ab} to \ref{prefsec} we describe in
detail how the abelian sheaf of vacua constructions in \cite{AU1}
yields a holomorphic line bundles with a projectively flat
connection and a preferred non-vanishing section over
Teichm\"{u}ller spaces of pointed surfaces equipped with a
symplectic basis of the first homology.

In section \ref{construction} we describe our global geometric
construction of a modular functor for any simple Lie algebra and a
level $\ell$.  Theorem \ref{mainconstT}
and \ref{mainconstTab} summarizes the constructions from sections
\ref{New1} to
\ref{prefsec}. The preferred section of the abelian theory allows
us to construct a certain fractional power of this line bundle
as stated in  Theorem \ref{fracpowerab}, which we tensor onto this
holomorphic vector bundle, so as to obtain a holomorphic vector
bundle with a {\em flat} connection over Teichm\"{u}ller space. The
modular functor is then defined (Definition \ref{def.main}) by
taking covariant constant sections of this flat bundle. The
section ends with the construction of the disjoint union
isomorphism. The glueing isomorphism is
constructed in section \ref{shofvacandglue}, where we also prove
the needed properties of glueing.

In section \ref{verification} we establish all the axioms
of a modular functor is satisfied based on the main results
of the preceding sections.

We have included an Appendix, which recalls the nessessary
definitions regarding nodal curves, families of stable curves and the
glueing construction.

\section{The axioms for a modular functor}\label{AxiomsMF}

We shall in this section give the axioms for a modular functor.
These are due to G. Segal and appeared first in \cite{Se}. We
present them here in a topological form, which is due to K. Walker
\cite{Walker}. See also \cite{Grove}. We note that similar, but
different, axioms for a modular functor are given in \cite{Tu} and
in \cite{BB}. It is however not clear if these definitions of a
modular functor is equivalent to ours.

Let us start by fixing a bit of notation. By a closed surface we mean a
smooth real two dimensional manifold.
For a closed oriented surface $\Si$ of genus $g$ we have the
non-degenerate skew-symmetric intersection pairing
\[(\cdot,\cdot) : H_1(\Si,\Z) \times H_1(\Si,\Z) \ra \Z.\]
Suppose $\Si$ is connected. In this case a Lagrangian subspace $L\subset H_1(\Si,\Z)$
is by definition a subspace, which is maximally isotropic with respect to the
intersection pairing. - A
$\Z$-basis $(\vec \alpha, \vec \beta) = (\alpha_1,\ldots, \alpha_g,\beta_1, \ldots
\beta_g)$ for $H_1(\Si,\Z)$ is called a symplectic basis if
\[(\alpha_i,\beta_j) = \delta_{ij}, \quad (\alpha_i,\alpha_j) = (\beta_i,\beta_j) = 0,\]
for all $i,j = 1, \ldots, g$.

If $\Si$ is not connected, then $H_1(\Si,\Z) =
\oplus_i H_1(\Si_i,\Z)$, where $\Si_i$ are the connected
components of $\Si$. By definition a Lagrangian subspace is in this
paper a subspace of the form $L = \oplus_i L_i$, where $L_i\subset
H_1(\Si_i,\Z)$ is Lagrangian. Likewise a symplectic basis for
$H_1(\Si,\Z)$ is a $\Z$-basis of the form $((\vec \alpha^i, \vec
\beta^i))$, where $(\vec \alpha^i, \vec
\beta^i)$ is a symplectic basis for $H_1(\Si_i,\Z)$.

For any real vector space $V$, we define $PV = (V-\{0\})/\R{}_+.$

\begin{definition}\label{DefPointS}
A {\em pointed surface} $(\Si,P)$ is an oriented closed
surface $\Si$ with a finite set $P\subset \Si$ of points. A
pointed surface is called {\em \stable}if the Euler characteristic
of each component of the complement of the points $P$ is negative.
A pointed surface is called {\em \fm}if each component of $\Si$
contains at least one point from $P$.
\end{definition}

\begin{definition}\label{DefMorPointS}
A {\em morphism of pointed surfaces} $f :(\Si_1,P_1) \ra
(\Si_2,P_2)$ is an isotopy class of orientation preserving
diffeomorphisms which maps $P_1$ to $P_2$. Here the isotopy is
required not to change the induced map of
 the first order Jet at $P_1$ to the first order Jet at $P_2$.
\end{definition}

\begin{definition} \label{msurface}

A {\em marked surface\/} $ {\Sib} = (\Si, P, V, L)$ is an oriented
closed smooth surface $\Si$ with a finite subset $P \subset \Si$
of points with projective tangent vectors $V\in \sqcup_{p \in
P}PT_{p}\Si$ and a Lagrangian subspace $L \subset H_1(\Si,\Z)$.
\end{definition}

\begin{remark} {\em The notions of \stable and \fm marked surfaces are
defined just like for pointed surfaces. }\end{remark}

\begin{definition} \label{mmorphism}

A {\em morphism\/} $\e f : {\Sib}_1 \to {\Sib}_2$ of marked
surfaces ${\Sib}_i = (\Si_i,P_i,V_i,L_i)$ is an isotopy class of
orientation preserving diffeomorphisms $f : \Si_1 \to \Si_2$ that
maps $(P_1,V_1)$ to $(P_2,V_2)$ together with an integer $s$.
Hence we write $\e f = (f,s)$.
\end{definition}

\begin{remark} {\em Any marked surface has an underlying pointed
surface, but a morphism of marked surfaces does not quit induce a
morphism of pointed surfaces, since we only require that the
isotopies preserve the induced maps on the projective tangent
spaces. }\end{remark}

Let $\sigma$ be Wall's signature cocycle for triples of Lagrangian
subspaces of $H_1(\Si,\R{})$ (See \cite{Wall}).

\begin{definition} \label{composition}
Let $\e f_1 = (f_1,s_1) : {\Sib}_1 \to {\Sib}_2$ and $\e f_2 =
(f_2,s_2) : {\Sib}_2 \to {\Sib}_3$ be morphisms of marked surfaces
${\Sib}_i = (\Si_i,P_i,V_i,L_i)$ then the {\it composition\/} of
$\e f_1$ and $\e f_2$ is $$ \e f_2 \e f_1 = (f_2 f_1, s_2 + s_1 -
\sigma((f_2f_1)_*L_1, f_{2*}L_2,L_3)). $$
\end{definition}

With the objects being marked surfaces and the morphism and their
composition being defined as in the above definition, we have
constructed
 the category of marked surfaces.

The mapping class group $\Gamma({\Sib})$ of a marked surface
${\Sib} = (\Si,L)$ is the group of automorphisms of ${\Sib}$. One
can prove that $\Gamma({\Sib})$ is a central extension of the
mapping class group $\Gamma(\Si)$ of the surface $\Si$ defined by
the 2-cocycle $c : \Gamma({\Sib}) \to \mathbb Z$, $c(f_1,f_2) =
\sigma((f_1f_2)_*L,f_{1*}L,L)$. One can also prove that this
cocycle is equivalent to the cocycle obtained by considering
two-framings on mapping cylinders (see \cite{At1} and \cite{A}).

Notice also that for any morphism $(f,s) : \Sib_1 \to \Sib_2$, one
can factor
\begin{eqnarray*}
(f,s) &=& \left((\id,s') : \Sib_2 \to \Sib_2\right) \circ
(f,s-s')\\ &=& (f,s-s') \circ \left((\id,s') : \Sib_1 \to
\Sib_1\right).
\end{eqnarray*}
In particular $(\id,s) : {\Sib} \to {\Sib}$ is $(\id,1)^s$.

\begin{definition} \label{disjunion}
The operation of {\em disjoint union of marked surfaces} is $$
(\Si_1,P_1,V_1,L_1)
 \sqcup (\Si_2,P_2,V_2,L_2) = (\Si_1 \sqcup \Si_2,P_1 \sqcup P_2,V_1\sqcup V_2,L_1 \oplus
L_2). $$

Morphisms on disjoint unions are accordingly $(f_1,s_1) \sqcup
(f_2,s_2) = (f_1 \sqcup f_2,s_1 + s_2)$.
\end{definition}

We see that disjoint union is an operation on the category of
marked surfaces.

\begin{definition}\label{or}
Let ${\Sib}$ be a marked surface. We denote by $- {\Sib}$ the
marked surface obtained from ${\Sib}$ by the {\em operation of
reversal of the orientation}. For a morphism $\e f = (f,s) :
{\Sib}_1 \to {\Sib}_2$ we let the orientation reversed morphism be
given by
 $- \e f = (f,-s) : -{\Sib}_1 \to -{\Sib}_2$.
\end{definition}

We also see that orientation reversal is an operation on the
category of marked surfaces. Let us now consider glueing of marked
surfaces.

Let $(\Si, \{p_-,p_+\}\sqcup P,\{v_-,v_+\}\sqcup V,L)$ be a marked
surface, where we have selected an ordered pair of marked points
with projective tangent vectors $((p_-,v_-),(p_+,v_+))$, at which
we will perform the glueing.

Let $c : P(T_{p_-}\Si) \ra P(T_{p_+}\Si)$ be an orientation
reversing projective linear isomorphism such that $c(v_-) = v_+$.
Such a $c$ is called a {\em glueing map} for $\Si$. Let
$\tilde{\Si}$ be the oriented surface with boundary obtained from
$\Si$ by blowing up $p_-$ and $p_+$, i.e.
\[\tilde{\Si} = (\Si -\{p_-,p_+\})\sqcup P(T_{p_-}\Si)\sqcup P(T_{p_+}\Si),\]
with the natural smooth structure induced from $\Si$. Let now
$\Si_c$ be the closed oriented surface obtained from $\tilde{\Si}$
by using $c$ to glue the boundary components of $\tilde{\Si}$. We
call $\Si_c$ the glueing of $\Si$ at the ordered pair
$((p_-,v_-),(p_+,v_+))$ with respect to $c$.

Let now $\Si'$ be the topological space obtained from $\Si$ by
identifying $p_-$ and $p_+$. We then have natural continuous maps
$q : \Si_c \ra \Si'$ and $n : \Si \ra \Si'$. On the first homology
group $n$ induces an injection and $q$ a surjection, so we can
define a Lagrangian subspace $L_c \subset H_1(\Si_c,\Z)$ by $L_c =
q_*^{-1}(n_*(L))$. We note that the image of $P(T_{p_-}\Si)$ (with
the orientation induced from $\tilde{\Si}$) induces naturally an
element in $H_1(\Si_c,\Z)$ and as such it is contained in $L_c$.

\begin{remark}{\em \label{remarkglue2}
If we have two glueing maps $c_i : P(T_{p_-}\Si) \ra
P(T_{p_+}\Si),$ $i=1,2,$ we note that there is a diffeomorphism
$f$ of $\Si$ inducing the identity on
$(p_-,v_-)\sqcup(p_+,v_+)\sqcup(P,V)$ which is isotopic to the
identity among such maps, such that $(df_{p_+})^{-1} c_2 df_{p_-}
= c_1$. In particular $f$ induces a diffeomorphism $f : \Si_{c_1}
\ra \Si_{c_2}$ compatible with $f : \Si \ra \Si$, which maps
$L_{c_1}$ to $L_{c_2}$. Any two such diffeomorphisms of $\Si$
induces isotopic diffeomorphisms from $\Si_1$ to
$\Si_2$.}\end{remark}

\begin{definition} \label{glueing}
Let ${\Sib} = (\Si, \{p_-,p_+\}\sqcup P,\{v_-,v_+\}\sqcup V,L)$ be
a marked surface. Let $$c : P(T_{p_-}\Si) \ra P(T_{p_+}\Si)$$ be a
glueing map and $\Si_c$ the glueing of $\Si$ at the ordered pair
$((p_-,v_-),(p_+,v_+))$ with respect to $c$. Let $L_c \subset
H_1(\Si_c,\Z)$ be the Lagrangian subspace constructed above from
$L$. Then the marked surface ${\Sib}_c = (\Si_c,P,V,L_c)$ is
defined to be the {\em glueing} of ${\Sib}$ at the ordered pair
$((p_-,v_-),(p_+,v_+))$ with respect to $c$.
\end{definition}

We observe that glueing also extends to morphisms of marked
surfaces which preserves the ordered pair $((p_-,v_-),(p_+,v_+))$,
by using glueing maps which are compatible with the morphism in
question.

We can now give the axioms for a 2 dimensional modular functor.

\begin{definition} \label{DefLS}

A {\em label set\/} $\L$ is a finite set furnished with an
involution $\l \mapsto \hat \l$ and a trivial element $1$ such
that $\hat 1 = 1$.
\end{definition}

\begin{definition} \label{lmsurface}

Let $\L$ be a label set. The category of {\em $\L$-labeled marked
surfaces\/} consists of marked surfaces with an element of $\L$
assigned to each of the marked point and morphisms of labeled
marked surfaces are required to preserve the labelings. An
assignment of elements of $\L$ to the marked points of ${\Sib}$ is
called a labeling of ${\Sib}$ and we denote the labeled marked
surface by $({\Sib},\l)$, where $\l$ is the labeling.
\end{definition}

We define a labeled pointed surface similarly.

\begin{remark}{\em
The operation of disjoint union clearly extends to labeled marked
surfaces. When we extend the operation of orientation reversal to
labeled marked surfaces, we also apply the involution $\hat \cdot$
to all the labels. }\end{remark}

\begin{definition} \label{DefMF}
A {\em modular functor\/} based on the label set $\L$ is a functor
$V$ from the category of labeled marked surfaces to the category
of finite dimensional complex vector spaces satisfying the axioms
MF1 to MF5 below.
\end{definition}

\subsubsection*{MF1} {\it Disjoint union axiom\/}: The operation of disjoint
union of labeled marked surfaces is taken to the operation of
tensor product, i.e. for any pair of labeled marked surfaces there
is an isomorphism $$ V(({\Sib}_1,\l_1) \sqcup ({\Sib}_2,\l_2)) )
\cong V({\Sib}_1,\l_1) \otimes V({\Sib}_2,\l_2). $$ The
identification is associative.

\subsubsection*{MF2} {\it Glueing axiom\/}: Let ${\Sib} $ and ${\Sib}_c$ be
marked surfaces such that ${\Sib}_c$ is obtained from ${\Sib} $ by
glueing at an ordered pair of points and projective tangent
vectors with respect to a glueing map $c$. Then there is an
isomorphism $$ V({\Sib}_c,\lambda) \cong \bigoplus_{\m \in \L}
V({\Sib},\m,\hat \m,\l), $$ which is associative, compatible with
glueing of morphisms, disjoint unions and it is independent of the
choice of the glueing map in the obvious way (see remark
\ref{remarkglue2}).

\subsubsection*{MF3} {\it Empty surface axiom\/}: Let $\emptyset$ denote
the empty labeled marked surface. Then $$ \dim V(\emptyset) = 1.
$$

\subsubsection*{MF4} {\it Once punctured sphere axiom\/}: Let $\Sib = (S^2,
\{p\},\{v\},0)$ be a marked sphere with one marked point. Then $$
\dim V(\Sib,\l) = \left\{ \begin{array}{ll} 1,\qquad &\l = 1\\
0,\qquad & \l \ne 1.\end{array}\right. $$

\subsubsection*{MF5} {\it Twice punctured sphere axiom\/}: Let $\Sib = (S^2,
\{p_1,p_2\},\{v_1,v_2\},\{0\})$ be a marked sphere with two marked
points. Then $$ \dim V(\Sib,(\l,\mu)) = \left\{ \begin{array}{ll}
1, \qquad &\l = \hat \mu\\ 0,\qquad &\l \ne \hat
\mu.\end{array}\right. $$

In addition to the above axioms one may has extra properties,
namely

\subsubsection*{MF-D} {\it Orientation reversal axiom\/}:
The operation of orientation reversal of labeled marked surfaces
is taken to the operation of taking the dual vector space, i.e for
any labeled marked surface $({\Sib},\l)$ there is a pairing $$
\langle \cdot,\cdot\rangle : V({\Sib},\l) \otimes V(-{\Sib},\hat
\l) \ra \C{}, $$ compatible with disjoint unions, glueings and
orientation reversals (in the sense that the induced isomorphisms
$ V({\Sib},\l) \cong V(-{\Sib},\hat \l)^*$ and $V(-{\Sib},\hat \l)
\cong V({\Sib},\l)^*$ are adjoints).

\vskip.4cm

 and

\subsubsection*{MF-U} {\it Unitarity axiom\/} Every vector
space $V({\Sib},\l)$ is furnished with a hermitian inner product
$$ ( \cdot,\cdot ) : V({\Sib},\l) \otimes \overline{V({\Sib},\l)}
\to {\mathbb C} $$ so that morphisms induces unitary
transformation. The hermitian structure must be compatible with
disjoint union and glueing. If we have the orientation reversal
property, then compatibility with the unitary structure means that
we have a commutative diagrams $$\begin{CD} V({\Sib},\l) @>>\cong>
V(-{\Sib},\hat \l)^*\\ @VV\cong V @V\cong VV\\
\overline{V({\Sib},\l)^*} @>\cong>> \overline{V(-{\Sib},\hat \l)},
\end{CD}$$
where the vertical identifications come from the hermitian
structure and the horizontal from the duality.

The rest of the paper is concerned with the detailed geometric construction of modular functors
using conformal
field theory. However, we shall assume the reader is familiar with
\cite{Ue2} and \cite{AU1} and freely use the notations of these
two papers in this paper.

\section{Teichm\"{u}ller space and families of pointed Riemann Surfaces with formal neighbourhoods}\label{New1}

Let us first review some basic Teichm\"{u}ller theory. Let $\Si$
be a closed oriented smooth surface and let $P$ be finite set of
points on $\Si$.

\begin{definition}\label{mc}
A {\em marked Riemann surface} $\e C$ is a Riemann surface $C$ with
a finite set of marked points $Q$ and non-zero tangent vectors $W \in
T_{Q}C= \bigsqcup_{q\in Q} T_{q}C$.
\end{definition}

\begin{definition}\label{morphmc}
A {\em morphism} between marked Riemann surface is a biholomorphism of the
underlying Riemann surface which induces a bijection between the two sets of
marked points and tangent vectors at the marked points.
\end{definition}

The notions of \stable and \fm is defined just like for pointed
surfaces.

\begin{definition} \label{cs}
A {\em complex structure} on $(\Si,P)$ is a marked Riemann surface
$\e C = (C, Q, W)$ together with an orientation preserving
diffeomorphism $\phi : \Si \ra
C$ mapping the points $P$ onto the points $Q$.
Two such complex structures $\phi_{j} : (\Si,P) \ra
\e C_{j} = (C_{j},Q_{j},W_{j})$ are {\em equivalent} if there exists a
morphism of marked Riemann surfaces
$$\Phi : \e C_{1} \ra \e C_{2}$$
such that $\phi_{2}^{-1} \Phi \phi_{1} : (\Si,P) \ra
(\Si,P)$ is isotopic to the identity through maps inducing the identity on the first
order neighbourhood of $P$.
\end{definition}

We shall often in our notation suppress the diffeomorphism, when
we denote a complex structure on a surface.

\begin{definition}\label{Teichmsp}
The {\em Teichm\"{u}ller space} $\cT_{(\Si,P)}$ of the pointed surface $(\Si,P)$ is by
definition the set of equivalence classes of complex structures on
$(\Si,P)$.
\end{definition}

We note there is a natural projection map from $\cT_{(\Si,P)}$ to $T_P\Si = \sqcup_{p\in P}
T_{p}\Si,$ which we call $\pi_P$.

\begin{theorem}[Bers]
There is a natural structure of a finite dimensional complex analytic manifold on
Teichm\"{u}ller space $\cT_{(\Si,P)}$.
Associated to any morphism of pointed surfaces $f: (\Si_1,P_1) \ra (\Si_2,P_2)$
 there is a biholomorphism $f^* : \cT_{(\Si_1,P_1)}\ra
 \cT_{(\Si_2,P_2)}$ which is induced by mapping a complex structure $\e C = (C, Q,W)$,
 $\phi : \Si_1 \ra C$ to $\phi \circ f^{-1} : \Si_2 \ra C$. Moreover,
 compositions of morphisms go to compositions of induced
 biholomorphisms.
\end{theorem}

 There is an
action of $\RPP$ on $\cT_{(\Si,P)}$ given by scaling the tangent vectors. This
action is free and the quotient $\ctT_{(\Si,P)} =
\cT_{(\Si,P)}/\RPP$ is a smooth manifold, which we call the {\em reduced} Teichm\"{u}ller space
of the pointed surface $(\Si, P)$.
Moreover the projection map $\pi_P$ descend
to a smooth projection map from $\ctT_{(\Si,P)}$
to $\sqcup_{p\in P}
P(T_{p}\Si)$, which we denote $\pi^{(r)}_P$.
We denote the fiber of this map
over $V\in \sqcup_{p\in P}
P(T_{p}\Si)$ by $\cT_{(\Si,P,V)}$. Teichm\"{u}ller space of a marked surface $\Sib = (\Si, P, V,
L)$ is by definition $\cT_{\Sib}=\cT_{(\Si,P,V)}$, which we call the Teichm\"{u}ller space of the
marked surface. Morphisms of marked
surfaces induce diffeomorphism of the corresponding
Teichm\"{u}ller spaces of marked surfaces, which of course also behaves well under
composition. We observe that the self-morphism $(\id,
s)$ of a marked surface acts trivially on the associated Teichm\"{u}ller space for all integers $s$.
General Teichm\"{u}ller theory implies that

\begin{theorem}\label{contractT}
The Teichm\"{u}ller space $\cT_{\Sib}$ of any marked surface $\Sib$
is
contractible.
\end{theorem}

Now let us recall the definition of a formal neighbourhood of a
point on a Riemann surface.

\begin{definition}\label{curvwfn}
Let $C$ be a Riemann surface and $q$ a point on $C$. Let $\mathcal O_{C,q}$
be the stalk of $\mathcal O_C$ at $q$ and let
$\mathfrak{m}_q$ the maximal ideal in $\mathcal O_{C,q}$. We note that $\mathfrak{m}_q^n$, $n=0,1,2, \ldots$,
gives a filtration of $\mathcal O_{C,q}$. A formal $n$'th-order
neighbourhood at $q$ is a filtration preserving isomorphism
\[\mathcal O_{C,q}/\mathfrak{m}_q^{n+1} \cong \mathbb C [[ \xi ]] / (\xi^{n+1}).\]
Let $\hat{\mathcal O}_{C,q} = \lim_{n\ra \infty}\mathcal O_{C,q}/\mathfrak{m}_q^n$ be
the completion of $\mathcal O_{C,q}$
with respect to the filtration. A {\em formal neighbourhood\/} (or {\em formal coordinate\/})
at $q$ is a filtration preserving isomorphism
\[ \eta : \hat{\mathcal O}_{C,q}\cong \mathbb C[[\xi]].\]
\end{definition}

We note that we have a
canonical isomorphism

\begin{eqnarray*}
   {\Ocal}_{C,q}/{\mathfrak m}_q^{2}& \simeq &
       \C{}\oplus T^*_{q}C,\\
       f &\mapsto& (f(q), df_q).
\end{eqnarray*}

Hence a formal $1$'st order neighbourhood induces and is
determined by an isomorphism of $T^*_q C$ with $\C{}$. Hence a
formal $1$'st order neighbourhood determines and is determined by
a non-zero vector in $T^*_q C$, specified by the property that it
maps to $1\in \C{}$ or equivalently a vector in $T_q C$ pairing
to unity with this vector.

\begin{definition}
A pointed Riemann surface with formal neighbourhoods
$$\goX = (C; q_1,\ldots, q_N; \eta_1, \ldots, \eta_N)$$
is the following data: A Riemann surface $C$, an ordered $N$-tuple of
$N$ distinct points $(q_1,\ldots, q_N)$ on $C$ together with
formal neighbourhoods
\[\eta_j : \hat{\mathcal O}_{C,q_j}\cong \mathbb C[[\xi_j]]\]
for $j=1,\ldots, N$.
\end{definition}

We remark that a pointed Riemann surface with formal
neighbourhoods is an "$N$-pointed smooth curve with formal
neighbourhoods" in the sense of Definition 1.1.3. of \cite{Ue2}.

\begin{definition}\label{v1order}
For a pointed Riemann surface with formal neighbourhoods $\goX$,
we
denote by $c(\goX)$ the underlying marked Riemann surface. For a
labeled pointed Riemann surface with formal neighbourhoods $(\goX, \vec \l)$,
we denote by $c(\goX,\vec \l) = (c(\goX),\l)$ the underlying
labeled marked Riemann surface. Here $\l$ denotes the labeling of the marked
points of $c(\goX)$ induced by $\vec \l$.
\end{definition}

\begin{definition}
A family of  pointed Riemann Surfaces with formal
neighbourhoods
$$\goF = ( \pi \: \mathcal C \rightarrow \mathcal B; \vs; \veta
)$$
 is the following date:
\begin{itemize}
\item Connected complex manifolds $\mathcal C$ and $\mathcal B$, such
that $\Dim_{\mathbb C} C = \Dim_{\mathbb C} \mathcal B + 1$.
\item A holomorphic submersion $\pi : \mathcal C \ra \mathcal B$.
\item Holomorphic sections $s_j$, $j = 1,\ldots, N$ of $\pi$.
\item Filtered $\mathcal O_{\mathcal B}$-algebra isomorphisms
\[\eta_j  : \widehat{\mathcal{O} }_{/s_j} = \varprojlim_{n \to \infty}
\mathcal{O} _Y/I_{j}^{n} \simeq \mathcal{O} _{\mathcal B}[[\xi]],
\]
where $I_{j}$ is the defining ideal of $s_j({\mathcal B})$ in ${\mathcal C}$, $j = 1, \ldots, N$.
\end{itemize}

\end{definition}

Note that a family of pointed Riemann Surfaces with formal
neighbourhoods is a "Family of $N$-pointed smooth curves with formal
neighbourhoods" as in Definition 1.2.1 in \cite{Ue2}. See also the Appendix
at the end of this paper.

Let $\Si$ be a closed oriented smooth surface and let $P$ be finite set of $N$ marked points on
$\Si$, i.e. $(\Si,P)$ is a pointed surface.

For a connected smooth complex manifold $\mathcal B$ let
$Y = \Si \times \mathcal B$.

Let
$\goF = ( \pi \: \mathcal C \rightarrow \mathcal B; \vs; \veta )$
be a family of pointed Riemann surfaces with formal
neighbourhoods and assume we have a smooth fiber
preserving diffeomorphism $\Phi_\goF$ from
$Y$ to $\mathcal C$ taking the marked points to the sections
$\vs$ and inducing the identity on $\mathcal B$. This data induces
a unique holomorphic map $\Psi_\goF$ from $\mathcal B$ to the
Teichm\"{u}ller space $\cT_{(\Si,P)}$ of the
surface $(\Si,P)$ by the universal property of Teichm\"{u}ller space.

\begin{definition}\label{famonsurf}
The pair $(\goF,\Phi_\goF)$ is called a family of pointed Riemann
surfaces
with formal neighbourhoods {\em on} $(\Si,P)$. If $P'\subset P$ is a
strict subset, we say that $(\goF,\Phi_\goF)$ is called a family of pointed Riemann surfaces
with formal neighbourhoods {\em over} $(\Si,P')$.
\end{definition}

Often we will suppress $\Phi_\goF$ in our notation and just
write $\goF$ is a family of pointed Riemann
surfaces with formal neighbourhoods on $(\Si,P)$.

\begin{definition}\label{goodfamily}
If a family $\goF = ( \pi \: \mathcal C \rightarrow \mathcal B; \vs; \veta )$
of pointed Riemann
surfaces
with formal neighbourhoods on $(\Si,P)$,
as above, has the properties, that the base $\mathcal B$ is biholomorphic to
an open ball and that the induced map $\Psi_\goF$
is a biholomorphism onto an open subset of Teichm\"{u}ller space
$\cT_{(\Si,P)}$ then the family is said to be \goodp
\end{definition}

Note that if a family of pointed Riemann
surfaces
with formal neighbourhoods on $(\Si,P)$ is {\em versal} around
some point $b\in B$, in the sense of Definition 1.2.2 in
\cite{Ue2}, then there is a open ball around $b$ in $B$,
such that the restriction of the family to this neighbourhood is
good.

\begin{proposition}\label{coverT}
For a stable and saturated pointed surface $(\Sigma,P)$ the Teichm\"{u}ller space $\cT_{(\Si,P)}$
can be covered by images of such \good families.
\end{proposition}

This follows from Theorem 1.2.9 in \cite{Ue2}.

Suppose now that we have two \stable and \fm families $\goF_i$, $i=1,2$
with the property that they have the same image $\Psi_{\goF_1}({\mathcal B}_1) =
\Psi_{\goF_2}({\mathcal B}_2)$ in Teichm\"{u}ller space
$\cT_{(\Si,P)}$ and that $\goF_2$ is a \good family.

\begin{proposition}\label{famequivalence}
For such a pair of families there exists a unique fiber preserving
biholomorphism $\Phi : \mathcal C_1\ra \mathcal C_2$ covering
$\Psi^{-1}_{\goF_2}\Psi_{\goF_1}$ such that $\Phi^{-1}_{\goF_2}
\Phi \Phi_{\goF_1} : (Y, P) \ra (Y,P)$ is isotopic to
$\Psi^{-1}_{\goF_2}\Psi_{\goF_1}\times \id$ through such fiber
preserving maps inducing the identity on the first order neighbourhood of
$P$.
\end{proposition}

This follows from uniqueness of the
$\Phi$ in Definition \ref{cs}.

We note that there is some permutation $S$ of $\{1,\ldots, N\}$ such that $(S \Phi^* (\veta_2))^{(1)} =
(\veta_1)^{(1)}$, i.e. $S \Phi^* (\veta_2)$ induce
the same first order formal neighbourhoods as $\veta_1$ does.

Suppose now $f$ is an orientation preserving diffeomorphism from $(\Si_1, P_1)$ to $(\Si_2,
P_2)$. Let $\goF_1$ be a family of pointed Riemann surfaces with formal
neighbourhoods of $(\Si_1, P_1)$. By composing $\Phi_{\goF_1}$
with $f^{-1}\times \id$ we get a family of pointed Riemann surfaces with formal
neighbourhoods of $(\Si_2, P_2)$. We note that $\Psi_{\goF_2}= f^* \circ \Psi_{\goF_1}$, where
$f^*$ is the induced map between the Teichm\"{u}ller spaces. This operation on families
clearly behaves well under compositions of diffeomorphisms.

\section{The space of vacua associated to a labeled marked Riemann Surface}

\subsection{Affine Lie algebras and integrable highest weight modules}

  In this section we recall the basic facts
about integrable highest weight representations
of affine Lie algebras.
For the details of integrable highest weight representations of
affine
Lie algebras we refer the reader to Kac's book [Ka].

   Let $\goth g$ be a simple Lie algebra over the complex numbers
 $\mathbb C$, which we fix throughout the paper. Let
$\goth h$ be its Cartan subalgebra. By $\Delta$ we denote the root
system of
$(\goth g , \goth h)$. We have the root space decomposition
$$
    \goth g = \goth h \oplus \displaystyle{\sum_{\alpha \in
        \Delta}} \goth g_\alpha .
$$

    Let ${\goth h}_{\mathbf R}^{*}$ be the linear span of $\Delta$
over $\mathbf R$.
Fix a choice of positive roots $\Delta_+$.

 Let $(\phantom{X}, \phantom{X} )$ be a
constant multiple of the Cartan-Killing form of the simple Lie algebra $\goth g$.
For each element of $\lambda \in \goth h$, there exists a unique element
$H_\lambda \in \goth h^*$  such that
$$
    \lambda(H) = ( H_\lambda, H)
$$
for all $H \in \goth  h$.  For $\alpha \in \Delta$, $H_\alpha$ is called the {\it root
vector\/}   corresponding to the root $\alpha$.

  On $\goth h^*$ we introduce
an inner product by
\begin{equation}
   (\lambda, \mu) = ( H_\lambda, H_\mu ).   \label{2.1-1}
\end{equation}
 Let us normalize the inner product $(\phantom{X}, \phantom{X})$ by
 requiring that $\theta$, the
highest (or longest) root, has length squared
\begin{equation}
   (\theta ,\, \theta) = 2. \label{2.1-2}
\end{equation}

 Let $V_\lambda$ be the irreducible left $\goth g$-module
 of
highest weight $\lambda$.
It is well-known that a finite dimensional irreducible left
$\goth g$-module is a highest weight module and two irreducible
left $\goth g$-modules are isomorphic if and only if they have
the same highest weight. A weight $\lambda \in
\goth h_{\mathbf R}^*$ is called an integral weight, if
$$
     2(\lambda, \alpha)/(\alpha, \alpha) \in \mathbf Z
$$
for any $\alpha \in \Delta$. A weight $\lambda \in
\goth h_{\mathbf R}^*$ is called a dominant weight, if
$$
  w(\lambda) \leq \lambda
$$
for any element $w$ of the Weyl group $W$ of $\goth g$.
By $P_+$ we denote the set of dominant integral weights of
$\goth g$. A weight $\lambda$ is the highest weight of an
irreducible left $\goth g$-module if and only if $\lambda
\in P_+$.

Let $w$ be longest element of $W$. Then we define an involution
$\dagger$ on $P_+$ by
\begin{equation}\label{involution}
\lambda^\dagger = - w(\lambda).
\end{equation}
One has that the opposite of dual of the left-$\goth g$-module $V_\lambda$
is isomorphic to left-$\goth g$-module $V_{\lambda^\dagger}$, meaning there exists a
non-degenerate $\goth g$-invariant perfect pairing
\[(\phantom{X}, \phantom{X}) : V_\lambda \otimes V_{\lambda^\dagger} \ra {\mathbb C}.\]

As mentioned in the introduction, we will need to fix $|0\rangle\in V_0 \setminus \{0\}$ and
\[|0_{\lambda,\lambda^\dagger}\rangle\in (V_\lambda\otimes V_{\lambda^\dagger})^{\goth g} \setminus \{0\},\]
where we put $|0_{\lambda \lambda^\dagger }\rangle= |0\rangle \otimes |0\rangle$ for $\lambda=0$ .
Fixing such a vector is of course equivalent to fixing the above
mentioned pairing.

    By $\mathbf C[[ \xi ]]$ and $\mathbf C((\xi))$ we mean the ring
of formal power series
in $\xi$ and the field of formal Laurent power series in $\xi$,
respectively.

\begin{definition}\label{D2.1.2}
{\rm The affine Lie algebra $\widehat{\goth g}$ over
$\mathbf  C((\xi))$  associated with $\goth g$ is defined to be
$$
  \widehat{\goth g} = \goth g \otimes {\mathbf  C((\xi))} \oplus
\mathbf C c
$$
where $c$ is an element of the center of $\widehat{\goth g}$
and the Lie algebra structure is
given by
\begin{eqnarray*}
 [X\otimes f(\xi), \,  Y\otimes g(\xi)] =
     [X, Y]\otimes f(\xi)g(\xi) + c\cdot (X,Y) \Res_{\xi = 0}
         (g(\xi)df(\xi))
\end{eqnarray*}
for}
$$
 X, \, Y \in \goth g, \, f(\xi),\, g(\xi) \in \mathbf  C((\xi)).
$$
\end{definition}
Put
\begin{equation}
\widehat{\goth g}_+ = \goth g \otimes \mathbf C[[\xi]]\xi,
\quad \widehat{\goth g}_- =\goth g \otimes
      \mathbf C[\xi^{-1}]\xi^{-1}.  \label{2.1-5}
\end{equation}

We regard $\widehat{\goth g}_+$ and
$\widehat{\goth g}_-$ as Lie subalgebras of $\widehat{\goth g}$.
We have a decomposition
\begin{equation}
 \widehat{\goth g} = \widehat{\goth g}_+ \oplus
 \goth g \oplus \mathbf C c \oplus \widehat{\goth g}_-.
\label{2.1-6}
\end{equation}

   Let us fix a positive integer  $\ell$
(called the  {\sl level\/}) and put
\begin{equation}\label{labelset}
   P_\ell = \{\, \lambda \in P_+ \, | \,
   0\le(\theta, \lambda) \le\ell \, \}.
\end{equation}
For all levels $\ell$ we observe that $\dagger$ takes $P_\ell$ to
it self.

For each element $\lambda \in P_\ell$ we shall define the Verma
module $\mathcal M_\lambda$ as follows.
Put

$$
  \widehat{\goth p}_+ := \widehat{\goth g}_+ \oplus \goth g \oplus \mathbf C \cdot c.
$$
Then $\widehat{\goth p}_+$ is a Lie subalgebra of $\widehat{\goth g}$. Let
$V_\lambda$ is the irreducible left $\goth g$-module of
highest weight $\lambda$. The action of $\widehat{\goth p}_+$ on
$V_\lambda$ is defined as

\begin{eqnarray*}
  c v & =& \ell v \quad \text{for all $v \in V_\lambda$}\\
  a v &= &0 \quad \text{for all
    $a \in \widehat{\goth g}_+$ and $v \in V_\lambda$}
\end{eqnarray*}

Put

\begin{equation}
   \mathcal M_\lambda := U( \widehat{\goth g})
         \otimes_{\widehat{ \goth p}_+}  V_\lambda. \label{2.1-14}
\end{equation}

Then $\mathcal M_\lambda$ is a left $\widehat{\goth g}$-module and
is called a {\it Verma module\/}. The Verma module
$\mathcal M_\lambda$ is not irreducible and contains the maximal
proper submodule $\mathcal J_\lambda$. The quotient module
$\mathcal H_\lambda :=  \mathcal M_\lambda/\mathcal J_\lambda$ has the
following properties.
\begin{theorem} \label{T2.1.4} For each $\lambda \in P_\ell$,
the  left $\widehat{\goth g}$-module ${\mathcal H}_\lambda$ is the
unique left $\widehat{\goth g}$-module
{\rm (}called the
{\rm integrable highest weight $\widehat{\goth g}$-module)}
satisfying the following
properties.
\begin{itemize}
\item $V_\lambda = \{ \, |v\rangle \in {\mathcal H}_\lambda \, | \;
  \widehat{\goth g}_+ |v\rangle = 0 \, \}$
is the irreducible left $\goth g$-module with highest
weight $\lambda$.
\item The central element $c$ acts on $\mathcal H_\lambda$ as
 $\ell \cdot \hbox{\rm id}$.
\item $\mathcal H_\lambda$ is generated by $V_\lambda$ over
$\widehat{\goth g} _-$ with only
one relation
\begin{equation}
  (X_\theta \otimes \xi^{-1} ) ^{\ell - (\theta, \lambda) + 1} |
  \lambda\rangle = 0
\label{2.1-15}
\end{equation}
where $X_\theta \in \goth g$ is the element corresponding to
the maximal root
$\theta$ and $|\lambda\rangle \in V_\lambda$ is a highest
weight vector.
\end{itemize}
\end{theorem}
The theorem says that the maximal proper submodule
$\mathcal J_\lambda$ is given by
\begin{equation}
\mathcal J_\lambda = U(\widehat{\goth p}_- )|J_\lambda \rangle
\label{2.1-16}
\end{equation}
where we put
\begin{equation}
 |J_\lambda \rangle =
 (X_\theta \otimes \xi ^{-1} ) ^{\ell - (\theta, \lambda) + 1} |
  \lambda\rangle .   \label{2.1-17}
\end{equation}
For the details see (10.4.6) in \cite{Ka}.

Similarly we have the integrable lowest weight right
 $\widehat{\goth g}$-module
${\mathcal H}_{\lambda}^\dagger$ which will be discussed below.

\subsection{The Segal-Sugawara construction}
We use the following notation
\begin{eqnarray*}
   X(n)   & = & X \otimes \xi^n, \quad X \in \goth g \\
   X(z)  & = &\sum_{n \in {\mathbf Z}} X(n) z^{-n-1}
\end{eqnarray*}
where $z$ is a variable. The normal ordering
$\Bosonnormalord \phantom{X}\Bosonnormalord $ is defined
by
$$
  \Bosonnormalord X(n) Y(m)\Bosonnormalord  = \begin{cases} X(n) Y(m), &n<m, \\
                 \frac{1}{2}(X(n)Y(m) + Y(m)X(n)) & n=m, \\
                Y(m) X(n) & n > m. \end{cases}
$$
Note that, if $n>m$ and $X=Y$,  we have
\begin{equation}
    \Bosonnormalord X(n) X(m) \Bosonnormalord  =
                             X(n) X(m) - n \delta_{n+m,0}(X,X) \cdot c.
                             \label{2.2-1}
\end{equation}

\begin{definition}\label{EM}
The
{\em energy-momentum} tensor $T(z)$  of level $\ell$ is defined by
$$
  T(z) = \frac{1}{2(g^{*} + \ell)} \sum_{a=1}^{\dim \goth g}
 \Bosonnormalord J^a(z)J^a(z)\Bosonnormalord
$$
where $\{J^1,J^2,\ldots, J^{\dim \goth g}\}$ is an orthonormal basis of $\goth g$ with respect
to the Cartan-Killing form $(\phantom{X},\phantom{X})$ and $g^{*}$ is the dual
Coxeter number of $\goth g$.
\end{definition}

  Put
\begin{equation}
  L_n = \frac{1}{2(g^{*} + \ell)}\displaystyle{\sum_{m \in {\mathbf Z}}
      \sum_{a=1}^{\dim \goth g}} \Bosonnormalord J^a(m)
 J^a(n-m)\Bosonnormalord .  \label{2.2-2}
\end{equation}

Then we have the expansion
$$
  T(z) = \sum_{n \in {\mathbf Z}}L_n z^{-n-2}.
$$
  The operator $L_n$ is called the $n$'th Virasoro operator and it acts on
$\mathcal H_\lambda $.

   For $X \in \goth g$, $f= f(z) \in \bolc (( z))$  and
$\underline{\ell} = \ell(z)\displaystyle{\frac{d}{dz}}\in
\bolc (( z ))\displaystyle{\frac{d}{dz}}$
we use the following notation.
\begin{eqnarray*}
  X[f] &= &\Res_{z=0}(X(z)f(z)dz) \\
  T[\underline{\ell}] &=& \Res_{z=0}(T(z)\ell(z)dz).
\end{eqnarray*}

In particular, we have that
\begin{equation}
    L_0 = T[\xi \dfrac d{d \xi }]. \label{2.2-6}
\end{equation}

  To define a filtration $\{F_{\bullet}\}$ on $\mathcal H_\lambda $,
we first define
the subspace ${\mathcal H_\lambda }(d)$ of $\mathcal H_\lambda $
for a non-negative integer $d$ by
\begin{equation}
 \mathcal H_\lambda  (d) = \{ \, |v\rangle \in \mathcal H_\lambda  \,|\,
  \;\;  L_0 |v\rangle = ( d + \Delta_\lambda)|v\rangle \, \}
   \label{2.2-8}
\end{equation}
where
\begin{equation}
\label{2.2.8a}
 \Delta_\lambda = \frac{(\lambda, \lambda) + 2(\lambda, \rho)}
   {2(g^{*} + \ell)}, \quad \rho = \frac{1}{2}
 \sum_{\alpha \in \Delta_+}\alpha.
\end{equation}

The subspaces $\mathcal H_\lambda (d)$ are finite
dimensional vector space and one has that
$$
       \mathcal H_\lambda  = \bigoplus_{d=0}^\infty \mathcal H_\lambda (d).
$$

Now we define the filtration $\{ F_p\mathcal H_\lambda \}$ by
\begin{equation}
    F_p\mathcal H_\lambda  = \sum_{d=0}^p \mathcal H_\lambda  (d).
    \label{2.2-9}
\end{equation}

  Put
\begin{equation}
{\mathcal H}_\lambda^{\dag} (d) = \Hom_{\mathbf C}(\mathcal H_\lambda  (d), \bolc).
\label{2.2-10}
\end{equation}
Then the dual space ${\mathcal H}_\lambda^{\dag}$ of ${\mathcal H}_\lambda$
is defined to be
\begin{equation}
  {\mathcal H}_\lambda^{\dag}  = \Hom_\bolc(\mathcal H_\lambda , \bolc)
 =  \prod _{d=0}^\infty
        {\mathcal H}_\lambda^{\dag} (d) . \label{2.2-11}
\end{equation}
By our definition ${\mathcal H}_\lambda^{\dag}$ is a right $\widehat{\gog}$-module.
A decreasing filtration $\{F^p{\mathcal H}_\lambda^{\dag}\}$ is defined by
\begin{equation}
   F^p{\mathcal H}_\lambda^{\dag} = \prod_{d \ge p}
          {\mathcal H}_\lambda^{\dag} (d) . \label{2.2-12}
\end{equation}
There is a unique canonical perfect  bilinear pairing
\begin{equation}
\langle \phantom{X}|\phantom{X}\rangle : {\mathcal H}_\lambda^{\dag} \times
         \mathcal H_\lambda  \longrightarrow \bolc, \label{2.2-13}
\end{equation}
given on $V_\lambda^\dagger \otimes V_\lambda$ by evaluation and which satisfies the following equality for each $a \in \widehat{\gog}$.
$$
   \langle u|av\rangle = \langle ua|v\rangle, \quad
\text{for all } \langle u| \in {\mathcal H}_\lambda^{\dag} \; \text{and }
       |v\rangle \in \mathcal H_\lambda  \, .
$$
Put
$$
V_\lambda^{\dag} = \{\,\langle v| \in {\mathcal H}_\lambda^{\dag}\, |\,
\;\;    \langle v|\widehat{\gog}_- = 0 \; \}.
$$
It is easy to show that $V_\lambda^{\dag} = {\mathcal H}_\lambda^{\dag} (0)$ and
$V_\lambda^{\dag}$ is the irreducible right $\goth g$-module with lowest weight
$\lambda$. The integrable highest weight right $\widehat{\gog}$-module with lowest
weight $\lambda$ is generated by $V_\lambda^{\dag}$ over $\widehat{\gog}_+$ with only one
relation
$$
    \langle \lambda|(X_{-\theta} \otimes \xi)^{\ell -(\theta,\lambda) +1} = 0.
$$

   Now let us introduce the left $\goth g$-module structure on
$\mathcal H_\lambda^\dagger$ by
$$
   X(n)\langle \Phi| := - \langle \Phi|X(-n) .
$$
It is easy to check that this indeed defines the left
$\goth g$-module structure on $\mathcal H_\lambda^\dagger$.

 Now we give the relationship of the {\it left}  $\goth g$-module
$\mathcal H_\lambda^\dagger$ and $\mathcal H_{\lambda^\dagger}$.

\begin{lemma} \label{L2.2.12}  There exists a unique canonical bilinear pairing
$$
  (\phantom{X}|\phantom{X}) :
\mathcal H_\lambda \times {\mathcal H}_{\lambda^\dagger} \rightarrow \mathbf C
$$
such that we have
$$
 (X(n)u|v) + (u|X(-n)v) = 0
$$
for any $X \in \gog$, $n \in {\mathbf Z}$, $| u \rangle \in \mathcal H_\lambda$,
$| v \rangle \in \mathcal H_{\lambda^\dagger}$, the pairing is zero on
$\mathcal H_\lambda(d) \times \mathcal H_{\lambda^\dagger}(d')$, if $d \neq d'$ and it evaluates
to $1$ on $|0_{\lambda,\lambda^\dagger}\rangle$.
\end{lemma}
\begin{corollary}\label{C2.2.13}
This pairing induces a canonical left $\goth g$-module isomorphism
$$
\mathcal H_\lambda^\dagger \simeq \widehat{\mathcal H}_{\lambda^\dagger}
$$
where $\widehat{\mathcal H}_{\lambda^\dagger}$ is the completion of
${\mathcal H}_{\lambda^\dagger}$ with respect to the filtration $\{F_p\}$.
\end{corollary}

\subsection{The space of vacua}

\begin{definition} \label{D3.1.1} The Lie algebra $\widehat{\goth g}_N$
is defined as
$$
  \widehat{\goth g}_N = \bigoplus_{j=1}^N \goth g \otimes_\bolc
 \bolc((\xi_j)) \oplus \bolc c
$$
with the following commutation relations.
\begin{equation}
[(X_j \otimes f_j), (Y_j \otimes g_j)] = ([X_j, Y_j] \otimes f_j g_j)
 + c \sum_{j=1}^N (X_j,Y_j) \Res_{\xi_j = 0}(g_j df_j)
 \label{3.1-1}
\end{equation}
where $(a_j)$ means $(a_1, a_2, \ldots, a_N)$ and $c$ belongs to the
center of $\widehat{\goth g}_N$.
\end{definition}

Let $\mathfrak X = (C; q_1,q_2,
\ldots, q_N; \eta_1, \eta_2, \ldots, \eta_N)$ be a pointed saturated
Riemann Surface with formal neighbourhoods and define
\begin{equation}
  \widehat{\goth g}(\goth X) = \goth g \otimes_\bolc H^0(C,\mathcal O_C (*\sumjn q_j)).
\label{3.1-2}
\end{equation}
We have the natural embedding
$$
t = \oplus t_i : H^0(C,\mathcal O_C (*\sumjn q_j)) \hookrightarrow
\bigoplus_{j=1}^N\bolc((\xi_j))
$$
given by Laurent expansion using the formal neighbourhoods.
In the following we often regard $H^0(C,\mathcal O_C (*\sumjn q_j))$ as a subspace
of $\displaystyle{\bigoplus_{j=1}^N}\bolc((\xi_j))$.
One has by lemma 1.1.15 in \cite{Ue2}, that  $\widehat{\goth g}(\goth X)$ is a
Lie subalgebra of
$\widehat{\goth g}_N$.

    Let us fix a non-negative integer $\ell$. For each $\vec \lambda =
(\lambda_1, \ldots, \lambda_N) \in (P_\ell)^N$, the left
$\widehat{\goth g}_N$-module $\Hlam$ and a right $\widehat{\goth g}_N$-%
module $\Hdaglam$ are defined by
\begin{eqnarray*}
 \Hlam & = & \mathcal H_{\lambda_1} \otimes_\bolc \cdots \otimes_\bolc
         \mathcal H_{\lambda_N} \\
  \Hdaglam & = &{\mathcal H}_{\lambda_1}^\dagger \widehat{\otimes}_\bolc
  \cdots \widehat{\otimes}_\bolc {\mathcal H}_{\lambda_N}^\dagger.
 \end{eqnarray*}
The hats over the tensor product means that the algebraic
tensor product has been completed with respect to the induced
filtration.

For each element $X_j \in \goth g$, $f(\xi_j) \in \bolc((\xi_j))$, the action
$\rho_j$ of $X_j[f_j]$ on $\Hlam$ is given by
\begin{equation}
\rho_j(X_j[f_j])|v_1 \otimes \cdots \otimes v_N\rangle
   = |v_1 \otimes \cdots \otimes v_{j-1} \otimes
(X_j[f_j])v_j \otimes v_{j+1} \otimes \cdots v_N\rangle
\label{3.1-3}
\end{equation}
where $|v_1 \otimes \cdots \otimes v_N\rangle$ means
$|v_1\rangle \otimes \cdots \otimes |v_N\rangle$, $|v_j\rangle \in
\mathcal H_{\lambda_j}$.
The left $\widehat{\goth g}_N$-action is given by
\begin{equation}
(X_1\otimes f_1, \ldots, X_N \otimes f_N)|v_1\otimes \cdots v_N\rangle
 = \sumjn \rho_j(X_j[f_j])|v_1 \otimes \cdots v_N\rangle.
 \label{3.1-4}
\end{equation}
Similarly, the right $\widehat{\goth g}_N$-action on $\Hdaglam$ is defined
 by
\begin{equation}
\langle u_1 \otimes \cdots u_N| (X_1 \otimes f_1, \ldots, X_N \otimes f_N)
 = \sumjn \langle u_1 \otimes \cdots u_N|\rho_j(X_j[f_j]).
 \label{3.1-5}
\end{equation}
As a Lie subalgebra, $\widehat{\goth g}(\goth X)$ operates on
 $\Hlam$ and $\Hdaglam$ as
\begin{equation}
 (X \otimes f)| v_1 \otimes \cdots \otimes v_N\rangle  =
 \sumjn \rho_j(X \otimes t_j(f))|v_1 \otimes \cdots v_N\rangle
 \label{3.1-6}
\end{equation}
and as
\begin{equation}
\langle u_1 \otimes \cdots \otimes u_N| (X \otimes f)  =
 \sumjn \langle u_1 \otimes \cdots \otimes u_N| \rho_j(X \otimes t_j(f)).
\label{3.1-7}
\end{equation}
The pairing $\langle \phantom{X}|\phantom{X}\rangle$ introduced in
\ref{2.2-11}
induces a perfect bilinear pairing
\begin{equation}
\langle \phantom{X}|\phantom{X}\rangle :
\Hdaglam \times \Hlam  \rightarrow \bolc  \label{3.1-8}
\end{equation}
given by
\begin{equation*}
(\langle u_1 \otimes \ldots \otimes u_N|,|v_1\otimes \ldots \otimes u_N
\rangle)  \rightarrow \langle u_1|v_1\rangle\langle u_2|v_2\rangle \cdots
 \langle u_N|v_N\rangle
\end{equation*}
which is $\widehat{\goth g}_N$-invariant:
$$
\langle\Psi(X_j\otimes f_j)|\Phi\rangle = \langle\Psi |(X_j\otimes f_j)
\Phi\rangle.
$$

  Now we are ready to define the space of vacua attached to $\goth X$.
\begin{definition}\label{D3.1.3}
Assume that $\goth X$ is saturated. Put
\begin{equation}
\Vlam(\goth X) = \Hlam / \widehat{\goth g}(\goth X)\Hlam .
\label{3.1-9}
\end{equation}
\end{definition}
The vector space $\Vlam(\goth X)$ is called the {\it space of
 covacua}
attached to $\goth X$.
The {\it space of vacua} attached to $\goth X$ is defined as
\begin{equation}
  \Vdaglam (\goth X) = \Hom_\bolc(\Vlam(\goth X), \bolc).
  \label{3.1-10}
\end{equation}

One gets that
\begin{equation}
\Vdaglam(\goth X) = \{ \; \langle \Psi | \in \Hdaglam \;|
\langle \Psi|\;\widehat \gog(\goX) = 0 \;\}.
\label{3.1-11}
\end{equation}
Moreover, the pairing (\ref{3.1-8}) induces a perfect pairing
\begin{equation}
\langle\phantom{X}|\phantom{X}\rangle : \Vdaglam (\goX) \times \Vlam(\goX)
\rightarrow \bolc. \label{3.1-12}
\end{equation}

The following theorem is proved in \cite{Ue2}.

\begin{theorem} \label{T3.1.5} The vector spaces $\Vlam(\goX)$ and $\Vdaglam(\goX)$ are
finite-dimensional.
\end{theorem}

\subsection{Propagation of vacua}

For a pointed Riemann Surface with formal neighbourhoods
$\goth X = (C;q_1,\ldots, q_N;
\eta_1,\ldots, \eta_N)$ let $q_{N+1}$ be a point on $C \setminus \{q_1,\ldots, q_N\}$ and
$\eta_{N+1}$ a formal neighbourhood of $C$ at $q_{N+1}$. Put
$$
   {\widetilde {\goth X}} = (C; q_1, \ldots ,
q_N, q_{N+1};
  \eta_1, \ldots , \eta_N, \eta_{N+1}).
$$

Since there is a canonical  inclusion
\begin{eqnarray*}
     \Hlam &\longrightarrow &\Hlam \otimes {\mathcal H}_0 \\
     |v \rangle &\longrightarrow & |v \rangle \otimes |0 \rangle
\end{eqnarray*}

we have a canonical surjection
$$
 \widehat{\iota}^* \; : \; \Hdaglam \widehat{\otimes} {\mathcal H}_0^\dagger
    \longrightarrow   \Hdaglam \; .
$$
\begin{theorem}\label{T3.2.1} The canonical surjection
$\widehat{\iota}^*$ induces  a canonical {\em Propagation of vacua} isomorphism
$$
 P_{\widetilde{\goth X},\goth X} : \mathcal V_{\veclam,0}^\dagger(\widetilde{\goth X})
 \ra
\mathcal V_{\veclam}^\dagger(\goth X).
$$
\end{theorem}

\subsection{Change of formal neighbourhoods}

\medskip

We let ${\mathcal D}$ be the automorphism
group $\Aut \C{}((\xi))$ of the field $\C{}((\xi))$ of formal Laurent
series as a $\C{}$-algebra.
There is a natural isomorphism
\begin{eqnarray*}
{\mathcal D} & \simeq & \{ \; \sum_{n=0}^\infty a_n \xi^{n+1} \; | \; a_0 \ne 0\; \}\\
 h & \mapsto & h(\xi)
\end{eqnarray*}
where the composition $h \circ g$ of $h$, $g \in {\mathcal D}$ corresponds to
the formal power series $h(g(\xi))$.

Put
$$
{\mathcal D}^p = \{ \;h \in  {\mathcal D}\; | \; h(\xi) = \xi + a_p\xi^{p+1} + \cdots\;\}
$$
for a positive integer $p$. Then we have a filtration
$$
{\mathcal D}= {\mathcal D}^0 \supset {\mathcal D}^1 \supset {\mathcal D}^2 \supset \ldots
$$
Also let
\begin{eqnarray*}
\underline{d} &=& \C{}[[\xi]]\xi\frac{d}{d\xi} \\
\underline{d}^p &=& \C{}[[\xi]]\xi^{p+1}\frac{d}{d\xi} \quad p=0,1,2,\ldots
\end{eqnarray*}
Then, we have a filtration
$$
\underline{d} = \underline{d}^0 \supset \underline{d}^1 \supset \underline{d}^2
\supset \cdots
$$

For any $\underline{l} \in \underline{d}$ and $f(\xi) \in\C{}[[\xi]]$ define
$\exp (\underline{l})(f(\xi))$ by
$$
\exp (\underline{l})(f(\xi)) = \sum_{k=0}^\infty \frac{1}{k!}
(\underline{l}^k f(\xi)).
$$
Set
\begin{eqnarray*}
{\mathcal D}_+^0 &=& \{ \;h \in  {\mathcal D}\; | \; h(\xi) = a\xi +
a_1\xi^{2} + \cdots,  \quad a>0 \;\} \\
\underline{d}_+^0 &=& \{ \; l(\xi)\frac{d}{d\xi}\;|\;
l(\xi) = \alpha \xi + \alpha_1 \xi^2 +\cdots, \quad \alpha \in \R{}\;\}
\end{eqnarray*}
 Then, we have the following result.
 \begin{lemma}
 \label{lem6.1}
 The exponential map
 \begin{eqnarray*}
 \exp \; : \; \underline{d} & \rightarrow & {\mathcal D} \\
 \phantom{\exp \; : \;{}} \underline{l} & \mapsto & \exp(\underline{l})
 \end{eqnarray*}
 \label{exponential}
 is surjective. Moreover, the exponential map induces an
 isomorphism
 $$
 \exp \; : \; \underline{d}_+^0 \simeq {\mathcal D}_+^0 .
 $$
 \end{lemma}

 Since,  for any integer $n$,  we have
 $$
 \exp(2\pi n \sqrt{-1} \xi \frac{d}{d\xi}) = id,
 $$
 the exponential mapping is not injective on $\underline{d}$.

For any element
$\underline{l} \in {\mathcal D}_+^0$ we define $\exp(T[\underline{l}])$ by
$$ \exp(T[\underline{l}])= \sum_{k=0}^\infty
\frac{1}{k!}T[\underline{l}]^k. $$ Then, $\exp(T[\underline{l}])$
operates on $ \mathcal{H}_\lambda$ from the left and on $ \mathcal{H}_\lambda^\dagger$ from the right.

By Lemma \ref{exponential},  for any automorphism $h \in {\mathcal D}_+^0$,
there exist a unique $\underline{l} \in \underline{d}_+^0$ with
$\exp(\underline{l}) = h$. Now for $h \in {\mathcal D}_+^0$ define
the operator $G[h]$ by
$$
  G[h] = \exp (- T[\underline{l}])
$$
where $\exp (\underline{l}) = h$. Then, $G[h]$ operates on $ \mathcal{H}_\lambda$
from the left and on $ \mathcal{H}_\lambda^\dagger$ from the right.
Then the following important theorems hold.
\begin{theorem}[{[19, Theorem 3.2.4]}]
\label{thm6.1}
For any $h \in \cD_+^0$, $f(\xi)d\xi \in \bC((\xi))d\xi$,
$g(\xi) \in \bC((\xi))$   and
$\underline{l} = l(\xi)\frac{d}{d\xi} \in \bC((\xi))\frac{d}{d\xi}$,
we have the following equalities as operators on $\cF$ and $\cFd$.
\begin{eqnarray*}
(1) &&G[h](\psi[f(\xi)d\xi])G[h]^{-1}= \psi[h^*(f(\xi)d\xi)]  =
\psi[f(h(\xi))h'(\xi)d\xi] \\
(2) && G[h](\ovpsi[g(\xi)])G[h]^{-1}= \ovpsi[h^*(g(\xi))]  =
\ovpsi[g(h(\xi))] \\
(3)&& G[h_1\circ h_2] = G[h_1]G[h_2] \\
(4)&& G[h]T[\underline{l}]G[h]^{-1} = T[\ad(h)(\underline{l})]+
\frac16 \Res_{\xi=0}\big(\{h(\xi);\xi\}l(\xi)d\xi\big) .
\end{eqnarray*}
where $\{f(\xi); \xi\}$ is the Schwarzian derivative.
\end{theorem}

\begin{theorem}[{[19, Theorem 3.2.5]}]
\label{prop6.1}
For any $h_j \in {\mathcal D}_+^0$, $j=1,2, \ldots, N$ and a pointed
Riemann surface with formal neighbourhoods
$$\mathfrak{X}= (C; q_1,q_2, \ldots, q_N; \xi_1, \xi_2, \ldots, \xi_N)$$
put
$$
\mathfrak{X}_{(h)} =(C; q_1,q_2, \ldots, q_N;
h_1(\xi_1), h_2(\xi_2), \ldots, h_N(\xi_N)).
$$
Then, the isomorphism $G[h_1] \widehat{\otimes}\cdots
\widehat{\otimes}G[h_N]$
\begin{eqnarray*}
 \mathcal{H}_{\vec \lambda}^\dagger  & \rightarrow & \mathcal{H}_{\vec \lambda}^\dagger \\
\langle \phi_1\widehat{\otimes} \cdots \widehat{\otimes}\phi_N|
& \mapsto & \langle \phi_1G[h_1] \widehat{\otimes}\cdots
\widehat{\otimes}\phi_NG[h_N]|
\end{eqnarray*}
induces the canonical isomorphism
$$
G[\vh] = G[h_1] \widehat{\otimes}\cdots \widehat{\otimes}G[h_N]:
\Vdaglam(\mathfrak{X}) \rightarrow
\Vdaglam(\mathfrak{X}_{(h)})
$$
\end{theorem}

Let $\goX = (C;\vQ;\veta)$ be a Riemann surface with formal neighbourhoods and let
$q_{N+1}$ be a further point on
the curve $C$ and
$\eta_{N+1}$ a formal neighbourhood of $C$ at $q_{N+1}$. Put $ \vtQ  = (q_1, \dots ,
q_N, q_{N+1})$ and $ \vteta = (
  \eta_1, \dots , \eta_N, \eta_{N+1}).$ Let
$$
   {\widetilde {\goX}} = (C;\vtQ;\vteta).
$$

We have the canonical
isomorphism $P_{{\widetilde {\goX}},\goX}$ from $\mathcal V_{\veclam,0}^\dagger(\widetilde{\goX})$ to $
\mathcal V_{\veclam}^\dagger(\goX)$ as given in Theorem \ref{T3.2.1}.
Suppose now $\vxi$ is
another formal neighbourhood at $\vQ$ and that $\xi_{N+1}$ is a
formal neighbourhood at $q_{N+1}$. Let then $\vtxi= (\xi,
\xi_{N+1})$, $\goX' = (C;\vQ;\vxi)$ and ${\widetilde {\goX}} = (C;
   \vtQ;\vtxi).$  Let $\vh$ be the
   formal coordinate change $\vxi = \vh(\veta)$ and $\widetilde{\vh}$ the
   formal coordinate change $\vtxi = \vth(\vteta)$. Assume now $c(\goX,\vec \l) = c(\goX',\vec
   \l)$ and $c(\widetilde{\goX},\vec \l,0) = c(\widetilde{\goX}',\vec
   \l,0)$. Then $h_j\in \mathcal{D}^{p_j},$ $p_j\geq 1$ and ${\widetilde h}_j\in
   \mathcal{D}^{{\widetilde p}_j},$ ${\widetilde p}_j\geq 1$.
   We then get the
   following diagram:

\begin{equation}
\begin{CD}
  \mathcal V_{\veclam,0}^\dagger(\widetilde{\goX})  @>P_{{\widetilde \goX},\goX}>>
  \mathcal V_{\veclam}^\dagger(\goX)   \\
  @V G[\vh] VV @VV G[\vth] V\\
   \mathcal V_{\veclam,0}^\dagger(\widetilde{\goX}')
   @>P_{{\widetilde \goX}',\goX'}>>
                  \mathcal V_{\veclam}^\dagger(\goX')
\end{CD}\label{diagcocpro}
\end{equation}

\begin{proposition}\label{compcocpro}
The diagram (\ref{diagcocpro}) is commutative.
\end{proposition}

\proof
A simple explicit calculation shows that
\[L_k |0\rangle = 0\]
if $k> 0$ and
\[L_0 |0\rangle = \Delta_0 |0\rangle.\]
But by \eqref{2.2.8a} $\Delta_0=0 $.
From this we immediately get that
\[G[\tilde h_{N+1}] |0\rangle = |0\rangle,\]
which by the very construction of the propagation of vacua
isomorphism makes the above diagram commute. \eproof

\subsection{The definition of the space of vacua associated to a labeled marked Riemann surface.}\label{New3}
Let $\e C = (C,Q,W)$ be a pointed Riemann surface and let $\l$
be a labeling of $\e C$. We shall now define the space of vacua
attached to the pair $(\e C, \l)$. We do this by providing
canonical isomorphisms between the spaces of vacua associated to
all correspondingly label pointed Riemann surfaces with formal
neighbourhoods over $\e C$. First we treat the case where $\e C$ is saturated.

 We notice that the definition of
the space of vacua $\Vdaglam(\goX)$ associated to a pointed Riemann surface with formal
neighbourhoods $\goX = (C,\vQ,\veta)$
depends on the
ordering of the marked points $\vQ = (q_1, \dots, q_N)$.
Let $S$ be a permutation of $\{1, \ldots, N\}$. We then define $\goX_S = (C, S(vQ),
S(\veta))$. The permutation $S$ acting from $\Hlam$ to $\HSlam$, induces an isomorphism
\[S : \Vdaglam(\goX) \ra \VdagSlam(\goX_S).\]

Clearly, compositions of
permutations go to compositions of isomorphisms.

Let $\goX' = (C,\vQ,\veta)$ and $\goX' = (C,\vQ',\veta')$ be two pointed
Riemann surface  with formal
neighbourhoods such that $c(\goX) = \e C = c(\goX')$. Let $S$ be such that $S(\vQ) =
\vQ'$.
Let $\vh$ be the
formal change of coordinates from $S(\veta)$ to $\veta'$. Then
as discussed above we have the isomorphism
\begin{equation}
G[\vh]S : \Vdaglam(\goX) \ra \VdagSlam(\goX')\label{isospov}
\end{equation}
to identify the two spaces of vacua with.

\begin{definition}\label{Defspofv} Let $\e C = (C,Q,W)$ be a \fm marked
Riemann surface. Let $\l$ be a labeling of $\e C$ using the set
$P_\ell$.
The space of vacua associated to the labeled marked curve $(\e C, \l)$
is by definition
\[\Spofvlam(\e C) = \left. \coprod_{c(\goX,\vec \l) = (\e C,\l)}\Vdaglam(\goX)\right/ \sim,\]
where the disjoint union is over all labeled curves with formal neighbourhoods
with $(\e C,\l)$ as the underlying labeled marked curve, $\vec \lambda$ is compatible
with the labeling $\l$  and $\sim$ is the equivalence relation
generated by the preferred
isomorphisms  (\ref{isospov}).
\end{definition}

That the relation $\sim$ is an equivalence relation follows from  (1) and (2) in Theorem \ref{thm6.1}.
Further it is clear that

\begin{proposition}\label{Spofviso}
The natural quotient map from $\Vdaglam(\goX)$ to  $\Spofvlam(\e C)$
is an isomorphism
 for all labeled Riemann surfaces with formal neighbourhoods $(\goX, \vec \l)$
with $c(\goX,\vec \l) = (\e C,\l)$.
\end{proposition}

Suppose $(\e C_i,\l_i)$ are labeled marked Riemann surfaces and $\Phi : (\e C_1,\l_1) \ra
(\e C_2,\l_2)$ is a morphism of labeled marked Riemann Surfaces. Let $(\goX_2,\vec
\l)$ be a labeled pointed Riemann Surface with formal neighbourhoods such that $c(\goX_2,\vec \l_2) = (\e
C,\l_2)$. Let $\Phi^*\goX_2 = \goX_1$. Then $\Phi$ is a morphism
of labeled marked Riemann surfaces with formal neighbourhoods. We obviously
have that

\begin{proposition}\label{morphspofv}
The identity map on $\Hdaglam$ induces a linear isomorphism from
$\Vdaglam(\goX_1)$ to $\Vdaglam(\goX_2)$, which induces a well defined linear isomorphism
$\Vdag(\Phi)$ from $\Spofvlamone(\e C_1)$ to $\Spofvlamtwo(\e
C_2)$. Compositions of morphisms of labeled marked Riemann Surfaces
go to compositions of the induced linear isomorphisms.
\end{proposition}

Let $(\e C, \l)$ be a labeled marked Riemann surface which might not be \fmp\newline
Consider all \fm
labeled marked Riemann surfaces $(\e C', \l')$ obtained from $(\e C, \l)$ by
adding points labeled with the trivial label $0\in P_\ell$.

\begin{definition}\label{Defspofvnc}
The space of vacua associated to the labeled marked Riemann surface $(\e C, \l)$
is by definition
\[\Spofvlam(\e C) = \left. \coprod_{(\e C',\l')}\Spofvlamp(\e C')\right/ \sim,\]
where the disjoint union is over all labeled marked Riemann surfaces $(\e C', \l')$ discussed
above and $\sim$ is the equivalence relation
generated by the propagation of vacua isomorphisms given in Theorem
\ref{T3.2.1}, the permutations of the order of the marked points and the
change of formal coordinate isomorphisms given in Theorem \ref{prop6.1}.
\end{definition}

By Theorem \ref{thm6.1}, Proposition \ref{compcocpro} and 2,
this $\sim$ is also an equivalence relation. We also
remark that if we apply definition \ref{Defspofvnc} to a labeled marked Riemann surface $(\e C, \l)$
which is already saturated, then we obtain a space of vacua which
is naturally isomorphic to the space of vacua definition
\ref{Defspofv} produces.

\section{The bundle of vacua over Teichm\"{u}ller space}\label{New4}

\subsection{Definition of the sheaf of vacua}\label{Sheafofv}

 Let
$\mathfrak{F}  = ( \pi : \mathcal{C}  \rightarrow \mathcal{B}  ; s_1, \ldots, s_N;
\eta_1, \ldots,\eta_N)$ be a family of pointed saturated
Riemann Surfaces of genus $g$ with formal neighbourhoods.

The sheaf $\widehat{\gog}_N (\mathcal{B} )$ of  affine
Lie algebra over $\mathcal{B} $ is the sheaf of $\Ob$-module
$$
 \widehat{\gog}_N(\mathcal{B} )=
{\gog}\otimes_{\bolc }
 (\bojn \Ob((\xi_j))) \oplus \Ob
   \cdot c
$$
with the following commutation relation, which is $\mathcal{O} _{\mathcal{B} }$-bilinear.
\begin{eqnarray*}
     \lefteqn{[( X_1 \otimes f_1,  \ldots,  X_N \otimes f_N ),
  (Y_1 \otimes g_1, \ldots, Y_N \otimes g_N)] } \\
       &&  =     ([X_1,Y_1] \otimes  (f_1g_1),\ldots,
  [X_N,Y_N] \otimes (f_Ng_N))
             \oplus
          c \cdot \sumjn (X_j, Y_j) \Res_{{\xi_j}=0}(g_jdf_j)
\end{eqnarray*}
where $
X_j, \, Y_j \in {\mathfrak{g} }, \quad f_j, \, g_j \in \Ob(({\xi_j}))$ and
we require $c$ to be central.  Put
$$
 \widehat{\gog}(\goF) =
{\mathfrak{g} } \otimes_{\bolc }\pi_{*} (\Oc (*S))
$$
where we define
\begin{align*}
   S & = \sumjn s_j({\mathcal{B} }) \\
   \pi_{*} (\Oc (*S)) & =
%     \varinjlim_k \pi_{*} (\Oc (kS))\,.
\varinjlim_k \pi_{*} (\Oc (kS))\,.
\end{align*}

  Laurent expansion
using the formal neighbourhoods $\eta_j$'s gives an inclusion:
$$
 \tilde{t} : \pi_*(\Ob(*S)) \rightarrow
    \bojn \Ob(({\xi_j}))
$$
and we may regard $\widehat{\gog}(\goF)$ as a Lie subalgebra of
$\widehat{\gog}_N(\mathcal{B} )$.
For any $\veclam = (\lambda_1,
\ldots , \lambda_N) \in (P_\ell)^N$,  put
\begin{align*}
 \Hlam(\mathcal{B} ) & = \Ob \otimes_{\bolc }{\mathcal{H}}_\veclam \, ,   \\
 \Hdaglam(\mathcal{B} ) & =
\underline{\hbox{\rm Hom}}_{\Ob}(\Hlam(\mathcal{B} ), \Ob)
= \mathcal{O} _{\mathcal{B} } \otimes_{\bolc} \Hdaglam.
\end{align*}
The pairing \eqref{3.1-8} induces  an $\Ob$-bilinear pairing
\begin{equation} \label{perfectpairsheaf}
  \langle \phantom{X} | \phantom{X} \rangle :
       \Hdaglam(\mathcal{B} ) \times \Hlam(\mathcal{B} ) \rightarrow
  \Ob.
\end{equation}
 The sheaf of affine Lie algebra $\widehat{\gog}_N(\mathcal{B} )$ acts on $\Hlam(\mathcal{B} )$ and
$\Hdaglam(\mathcal{B} )$ by
\begin{align*}
        ( (X_1 \otimes \sum_{n \in {\mathbf{Z}}} a_n^{(1)}\xi_1^n),
  \ldots , (X_N & \otimes \sum_{n \in {\mathbf{Z}}} a_n^{(N)}\xi_N^n))
 (F \otimes |\Psi \rangle )   \\
  =  &  \sumjn \sum_{n \in {\mathbf{Z}}}(a_n^{(j)} F)
  \otimes \rho_j(X_j(n))|\Psi\rangle
\end{align*}
The action of $\widehat{\mathfrak g}_N(\mathcal{B} )$ on
${\mathcal{H}}_{\veclam}^\dagger (\mathcal{B} )$ is the dual action of
${\mathcal{H}}_{\veclam}(\mathcal{B} )$, that is,
$$
   \langle \Psi a | \Phi \rangle = \langle \Psi| a  \Phi \rangle
\quad \text{for any } \; a \in \widehat{\mathfrak{g} }_N .
$$
\begin{definition}\label{ShofVdef}
For the family $\goF$ of pointed Riemann surfaces with formal
neighbourhoods, we define the sheaves  of $\Ob$-modules on $\mathcal{B}$
\begin{align*}
 \Vlam(\goF) & = \Hlam(\mathcal{B} ) / \widehat{\gog}(\goF) \Hlam(\mathcal{B} ) \\
\Vdaglam(\goF) & = \underline{\hbox{\rm Hom}}_{\Ob}(\Vlam(\mathcal{B} ), \Ob).
\end{align*}
These are the
{\it sheaf of covacua\/} and the {\it sheaf  of vacua\/} attached to the family
$\goF$.
\end{definition}

Note that we have
$$
\Vdaglam(\goF)  = \{ \, \langle \Psi| \in \Hdaglam(\goF) \, | \;\;
      \langle \Psi | a = 0 \quad \text{for any }\; a \in \widehat{\gog}(\goF) \, \} .
$$
The pairing \eqref{perfectpairsheaf}) induces an $\Ob$-bilinear
pairing
$$
  \langle \phantom{X} |\phantom{X} \rangle :
     \Vdaglam(\goF) \times \Vlam(\goF) \rightarrow  \Ob .
$$

By Theorem 4.1.6 in \cite{Ue2} and Corollary 4.2.4 in \cite{Ue2} we have
the following highly non-trivial theorem.
Note that the theorem is true even for a family of nodal curves
(see Theorem \ref{stabllocalfreeness}).

\begin{theorem}\label{localfreeness}
The sheaves $ \Vdaglam(\goF)$ and $\Vlam(\goF)$ are locally free
sheaves of $\Ob$-modules of finite rank over $\mathcal B$. They are dual to each other.
\end{theorem}

Hence $\Vdaglam(\goF)$ is a holomorphic vector bundles over $\mathcal
B$.

\subsection{Properties of the sheaf of vacua.}\label{New4Props}

To construct the bundle of vacua over Teichm\"{u}ller space, we further need
the following obvious property of the sheaf of vacua construction.

\begin{lemma}\label{Lfamiso}
Let $\goF_i$ be two families of saturated pointed Riemann surfaces with formal neighbourhoods
over the same base $\mathcal B$.
Let $\Phi :
(\goF_1,\l_1) \ra (\goF_2,\l_2)$ be an isomorphism of labeled families,
which induces the identity map
 on the base. Then the identity map on $\Hdaglam(\mathcal B)$
 induces a canonical isomorphism
\begin{equation}\Vdag(\Phi) : \Vdaglamone(\goF_1) \ra
\Vdaglamtwo(\goF_2).\label{famiso}
\end{equation}
\end{lemma}

Suppose that we have two families of pointed Riemann surfaces with formal
neighbourhoods $\goF_i$, $i=1,2$
on a pointed saturated and stable surface $(\Si,P)$,
with the property that they have the same image $\Psi_{\goF_1}({\mathcal B}_1) =
\Psi_{\goF_2}({\mathcal B}_2)$ in Teichm\"{u}ller space
$\cT_{(\Si,P)}$ and that $\goF_2$ is a \good family.
For such a pair of families there exists by Proposition \ref{famequivalence} a unique
fiber preserving
biholomorphism $\Phi_{12} : \mathcal C_1\ra \mathcal C_2$ covering
$\Psi^{-1}_{\goF_2}\Psi_{\goF_1}$ such that $\Phi^{-1}_{\goF_2}
\Phi_{12} \Phi_{\goF_1} : (Y, P) \ra (Y,P)$ is isotopic to
$\Psi^{-1}_{\goF_2}\Psi_{\goF_1}\times \id$ through such fiber
preserving maps inducing the identity on the first order neighbourhood of $P$.

We note that $(\Phi_{12}^* (\veta_2))^{(1)} = (S\veta_1)^{(1)}$, where $S$ is some
permutation of $\{1, \ldots, N\}$, i.e. $\Phi_{12}^* (\veta_2)$
induce
the same first order formal neighbourhoods as $S\veta_1$ does. Let
$\vh$ be the formal change of coordinates from $S\veta_1$
to $\Phi_{12}^* (\veta_2)$. Let $\widetilde{\goF_1} = \Phi_{12}^*(\goF_2)$.
Then $\Phi_{12}$ induces an isomorphism of families from $\widetilde{\goF_1}$ to
$(\Psi^{-1}_{\goF_2}\Psi_{\goF_1})^{*}(\goF_2)$. Choose labelings
$\vec \l_1$ of $\goF_1$ and $\vec \l_2$ of $\goF_2$, which are compatible
under the above isomorphism.

\begin{proposition}\label{vbtransf}
The action of $S$ on $\Hdaglamone(\mathcal B_1)$ induces an
isomorphism from $\Vdaglamone(\goF_1)$ to $\Vdaglamtwo(S\goF_1)$, where
S acts on the family $\goF_1$ by permuting the numbering of the formal neighbourhoods
and the sections. Furthermore
$G[\vh] : \Hdaglamtwo(\mathcal B_1) \rightarrow \Hdaglamtwo(\mathcal B_1)$
induces an isomorphism from $\Vdaglamtwo(S\goF_1)$ to
$\Vdaglamtwo(\widetilde{\goF_1})$.
The natural pull back isomorphisms provided by Lemma
4.1.3 in \cite{Ue2} and the families isomorphism (\ref{famiso}) provided by Lemma
\ref{Lfamiso}
 induces an isomorphism from
$\Vdaglamtwo(\widetilde{\goF_1})$ to $(\Psi^{-1}_{\goF_2}\Psi_{\goF_1})^*\Vdaglamtwo(\goF_2)$.  The composite of these
three
isomorphism gives the transformation isomorphism
\begin{equation}G_{12} : \Vdaglamone(\goF_1) \ra
(\Psi^{-1}_{\goF_2}\Psi_{\goF_1})^*\Vdaglamtwo(\goF_2).\label{overlapiso}
\end{equation}
The isomorphisms satisfies the cocycle condition
\[G_{12} (\Psi^{-1}_{\goF_2}\Psi_{\goF_1})^*G_{23} = G_{13}.\]
\end{proposition}

\proof {The cocycle condition follows, since $G_{12}$ is induced from
the permutation $S$ and the isomorphism $G[\vec h]$ on the $\mathcal H$-level, combined
with formula (2) of Theorem 3.2.4 in \cite{Ue2}.}\eproof

Let $\goF = ( \pi : \mathcal{C}  \rightarrow \mathcal{B}  ; s_1, \ldots, s_N;
\eta_1, \ldots,\eta_N)$ be a family of pointed Riemann surfaces with formal
neighbourhoods. Let $s_{N+1}$ be a further section of $\pi$
disjoint from the $s_i$ and $\eta_{N+1}$ a formal neighbourhood
along $s_{N+1}$ and set $\widetilde{\goF} = ( \pi : \mathcal{C}  \rightarrow \mathcal{B}  ;
s_1, \ldots, s_{N+1};
\eta_1, \ldots,\eta_{N+1})$

\begin{theorem} \label{Povforfamilies}
The inclusion
\begin{eqnarray*}
     \Hlam(\mathcal B) &\longrightarrow &\Hlam(\mathcal B) \otimes {\mathcal H}_0(\mathcal B) \\
     |v \rangle &\longrightarrow & |v \rangle \otimes |0 \rangle
\end{eqnarray*}
induces the propagation of vacua isomorphism
$$
 P_{\widetilde{\goth F},\goth F} : \mathcal V_{\veclam,0}^\dagger(\widetilde{\goth F})
 \ra
\mathcal V_{\veclam}^\dagger(\goth F).
$$
\end{theorem}

\subsection{The definition of the bundle of vacua over Teichm\"{u}ller space}\label{New4D}

Let $\Si$ be a closed oriented smooth surface and let $P$ be a finite set of
marked points on
$\Si$. Let $\l$ be a labeling of $(\Si, P)$ and assume $(\Si, P)$
is \stable and \fmp
We now define a holomorphic vector bundle $\Spofvlam = \Spofvlam(\Si,P)$ over Teichm\"{u}ller space
$\cT_{(\Si,P)}$ using the cover $\{\Psi_{\goF}(\mathcal B)\}$, where $\goF$ runs
over the good families of complex structures on $(\Si,P)$.
\begin{definition}\label{dvbovac}
A holomorphic vector bundle $\Spofvlam = \Spofvlam(\Si,P)$ over Teichm\"{u}ller space
$\cT_{(\Si,P)}$ is specify to be the bundle $(\Psi_\goF^{-1})^*\Vdaglam(\goF)$ over $\Psi_{\goF}(\mathcal
B)$ for any good families of complex structures on $(\Si,P)$. On
overlaps of the image of two \good families, we use
the glueing isomorphism (\ref{overlapiso}) to glue the corresponding bundles together.
\end{definition}

Proposition \ref{vbtransf} implies that $\Spofvlam(\Si,P)$ is a vector bundle over $\cT_{(\Si,P)}$
with the following property.

\begin{proposition}\label{Teichpullback}
For any \stable and \fm family $\tgoF$ of pointed curves with
formal neighbourhoods on $(\Si,P)$ we have a preferred isomorphism
\[\Upsilon_{\tgoF} : \Vdaglam(\tgoF) \ra \Psi_{\tgoF}^*\Spofvlam(\Si,P)\]
induced by the transformation isomorphism between $\Vdaglam(\tgoF)$ and
$\Vdaglam(\goF)$, for \good families $\goF$ of complex structures on
$(\Si,P)$ such that $\Psi_\goF(\mathcal B)$ intersect $\Psi_{\tgoF}({\mathcal B}')$
nonempty.
\end{proposition}

Let $(\Si,P,\lambda)$ be a general labeled pointed surface, i.e.
$(\Si,P)$ might not be \stable nor \fmp Let
$(\Si,P_i,\lambda_i)$, $i=1,2$ be
any labeled marked surfaces obtained from $(\Si,P,\lambda)$ by
labeling further points by $0\in P_\ell$.
Assume that $(\Si,P_i)$ are \stable and \fm pointed surfaces. Let $\bar P = P_1 \cup
P_2$ and $\bar \l$ be the induced labeling of $\bar P$. Note that
$(\Si,\bar P)$ is also \stable and \fmp

We get holomorphic projection maps $\pi_i : \cT_ {(\Si,\bar P)}\ra \cT_
{(\Si,P_i)}$. As a direct consequence of Proposition
\ref{compcocpro} we get the following.

\begin{proposition}\label{propvaciso}
Iterations of the propagation of vacua isomorphism given in Theorem \ref{Povforfamilies} induces
natural isomorphisms of bundles
\[\Spofvlambar (\Sigma,\bar P) \cong \pi_1^*\Spofvlamone (\Sigma,P_1)
\cong \pi_2^*\Spofvlamtwo (\Sigma,P_2),\]
which satisfies associativity.
\end{proposition}

Suppose now $f : (\Si_1, P_1) \ra (\Si_2,P_2)$ is a morphism of \stable and \fm pointed
surfaces. Then of course $f$ induces a
morphism $f^*$ from $\cT_{(\Si_1,P_1)}$ to $\cT_{(\Si_2,P_2)}$.
Let $\l_1$ be a labeling of $(\Si_1,P_1)$ and let $\l_2$ be the
induced labeling on $(\Si_2,P_2)$ such that $f : (\Si_1,P_1,\l_1) \ra
(\Si_2,P_2,\l_2)$ is a morphism of labeled pointed surfaces.
Let now $\goF_1$ be a
\good family of stable pointed curves with formal
neighbourhoods of $(\Si_1,P_1)$. Then by composing with $f^{-1}\times
\id$ we get a \good family $\goF_2$ of stable pointed curves with formal
neighbourhoods on $(\Si_2,P_2)$ over the same base $\mathcal
B_1$. The identity morphism on $\Hdaglam(\mathcal B_1)$ then induces a
morphism $\Vdag(f) : \Vdaglam(\goF_1) \ra
\Vdaglam(\goF_2)$ which covers the identity on the base. This is
precisely the morphism induced from the morphism of families $\Phi_f =  f\times \id
: \goF_1 \ra
\goF_2$ by Lemma \ref{Lfamiso}. This intern then induces a
morphism $\Vdag(f) : (\Psi_{\goF_1}^{-1})^*(\Vdaglamone(\goF_1)) \ra
(\Psi_{\goF_2}^{-1})^*(\Vdaglamtwo(\goF_2))$ which covers
$f^* : \Psi_{\goF_1}(B_1) \ra \Psi_{\goF_2}(B_1)$.

\begin{proposition}\label{comptransf}
The above construction provides a well defined lift of $f^* : \cT_{(\Si_2,P_2)} \ra
\cT_{(\Si_1,P_1)}$ to a morphism $\Spofv(f) : \Spofvlamone(\Si_1,P_1)
\ra \Spofvlamtwo(\Si_2,P_2)$ which behaves well under compositions.
\end{proposition}

\proof {Let $f'$ be a diffeomorphism whose first order
neighbourhood from $P_1$ to $P_2$ is the same as $f$'s and such
that $f'$ is isotopic to $f$ among such. Let $\goF_{2'}$ be
obtained from $\goF_1$ by composing with $(f')^{-1}\times \id$.
Then
\[\Phi = \Phi_{\goF_1}\circ (((f')^{-1}\circ f) \times \id)
\circ \Phi^{-1}_{\goF_1} : \mathcal C_1 \ra \mathcal C_1\]
is the unique biholomorphism from Proposition \ref{famequivalence}.
Hence we get a commutative diagram
$$
\begin{CD}
  (\Psi_{\goF_1})^*(\Vdaglamone(\goF_1))  @> \Vdag(f) >>
    (\Psi_{\goF_2})^*(\Vdaglamtwo(\goF_2))  \\
  @V = VV @VV G_{22'} V\\
   (\Psi_{\goF_1})^*(\Vdaglamone(\goF_1))
   @> \Vdag(f') >>
(\Psi_{\goF_{2'}})^*(\Vdaglamtwo(\goF_{2'}))
\end{CD}
$$
which shows that $\Vdag(f) : \Vdaglamone(\Sigma_1,P_1) \ra
\Vdaglamtwo(\Sigma_2,P_2)$ is well defined. It is obvious that $\Vdag(fg) =
\Vdag(f)\Vdag(g)$.}

\eproof

Assume that $(\Si_i,P_i, \l_i)$ are labeled pointed surfaces,
which need not be neither \stable nor \fmp Let $f : (\Si_1,P_1,
\l_1) \ra (\Si_2,P_2,\l_2)$ be an orientation preserving
diffeomorphism of labeled pointed surfaces. Let $(\Si_i,
P'_i,\l'_i)$ be labeled pointed surfaces obtained from
$(\Si_i,P_i, \l_i)$ by labeling further points by $0\in P_\ell$
such that $(\Si_i, P'_i,\l'_i)$ are \stable and \fm labeled
pointed surfaces such that $f:(\Si_1,P'_1,\l'_1) \ra (\Si_2, P'_2,\l'_2)$ is a
morphism of labeled pointed surfaces. We obviously have the
following result.

\begin{proposition}\label{morphcomp}
The lift of $f^* : \cT_{(\Si_1,P'_1)} \ra \cT_{(\Si_2,P'_2)}$ to a
morphism $\Spofv(f) : \Spofvlampone(\Si_1,P'_1) \ra
\Spofvlamptwo(\Si_2,P'_2)$ as given by Theorem \ref{comptransf}
is compatible with the isomorphisms given
in Proposition \ref{propvaciso}.
\end{proposition}

\section{The connection in the bundle of vacua over Teichm\"{u}ller space}\label{New5}

\subsection{Twisted first order differential operators acting on the sheaf of vacua}

Let $\mathfrak{F}  = (\pi : \mathcal{C}  \rightarrow \mathcal{B} ; s_1, \ldots, s_N ; \eta_1, \ldots,
\eta_N)$ is a family of saturated pointed Riemann surfaces of genus $g$ with formal
neighbourhoods. We will assume that
$\goF^{(0)}= (\pi : \mathcal{C}  \rightarrow \mathcal{B} ; s_1, \ldots, s_N)$
be a {\it versal family\/}
 of pointed stable curves of genus $g$ in the sense of definition 1.2.2. in \cite{Ue2}.
We consider the divisors
$$
S_j = s_j(\mathcal{B} ), \qquad  S = \sum_{j=1}^N S_j.
$$

There is an exact sequence
$$
0 \rightarrow \Theta_{\mathcal{C} /\mathcal{B} } \rightarrow
 \Theta_{\mathcal{C} }  \overset{d\pi}{\rightarrow}
 \pi^*\Theta_{\mathcal{B} }     \rightarrow  0
$$
where $\Theta_{\mathcal{C} /\mathcal{B} }$ is a sheaf of
vector fields tangent to the fibres of $\pi$. Put
$$
\Theta'_{\mathcal{C} ,\pi} = d\pi^{-1}(\pi^{-1}\Theta_{\mathcal{B} }).
$$
Hence, $\Theta'_{\mathcal{C} ,\pi}$ is a sheaf of
vector field on $\mathcal{C} $
whose vertical components are constant along
the fibers of $\pi$. That is,  in a neighbourhood of a
point of a  fiber $\Theta'_{\mathcal{C} ,\pi}$
consists of germs of holomorphic vector fields of the form
$$
 a(z,u) \frac{\partial}{\partial z} +
  \sum_{i=1}^n b_i(u) \frac{\partial}{\partial u_i}
$$
where $(z,u_1 ,\ldots, u_n)$ is a system of local
coordinates such that the mapping $\pi$ is expressed as the projection
$$
 \pi (z,u_1 ,\ldots, u_n) =(u_1, \ldots, u_n).
$$

More generally, we can define a sheaf $\Theta'_{\mathcal{C} }(mS)_\pi$
as the one consisting of germs of meromorphic vector fields of the form
$$
 A(z,u) \frac{\partial}{\partial z} +
\sum_{i=1}^n B_i(u)  \frac{\partial}{\partial u_i}
$$
where $A(z,u)$ has the poles of order at most $m$
along $S$.

We have an exact sequence
$$
 0 \rightarrow \Theta_{\mathcal{C} /\mathcal{B} }(mS) \rightarrow
 \Theta'_{\mathcal{C} }(mS)_\pi  \overset{d\pi}{\rightarrow }
 \pi^{-1}\Theta_{\mathcal{B} }    \rightarrow  0 .
$$
Note that $\Theta'_{\mathcal{C} }(mS)_\pi $ has the structure of  a
sheaf of Lie algebras by the usual bracket operation on  vector
fields and the above exact sequence is one of sheaves of Lie algebras.

For $m > \dfrac1N(2g - 2)$  we have an exact sequence of
$\mathcal{O} _\mathcal{B} $-modules.
$$
 0 \rightarrow \pi_* \Theta_{\mathcal{C} /\mathcal{B} }(mS)  \rightarrow
\pi_*  \Theta'_{\mathcal{C} }(mS)_\pi  \overset{d\pi}{\rightarrow }
\Theta_{\mathcal{B} }    \rightarrow  0
$$
which is also an exact sequence of sheaves of Lie algebras.
Taking $m \rightarrow \infty$ we obtain the exact sequence
$$
 0 \rightarrow \pi_* \Theta_{\mathcal{C} /\mathcal{B} }(*S)  \rightarrow
\pi_*  \Theta'_{\mathcal{C} }(*S)_\pi  \overset{d\pi}{\rightarrow }
\Theta_{\mathcal{B} }   \rightarrow  0 .
$$

Recall that
we have  the following exact sequence of
${\mathcal{O} }_{\mathcal{B} }$-modules.

$$
0 \rightarrow  \Theta_{\mathcal{C} /\mathcal{B} }( - S))
\rightarrow \Theta_{\mathcal{C} /\mathcal{B} }( mS))
\rightarrow
 \bigoplus_{j=1}^N \bigoplus_{k=0}^m
{\mathcal{O} }_{\mathcal{B} }\xi_j^{-k}
\displaystyle{\frac{d}{d\xi_j}}
 \rightarrow 0
$$
Which for any positive integer $m \geq 4g -3$ gives the
following exact sequence
$$
 0 \rightarrow \pi_*(\Theta_{\mathcal{C} /\mathcal{B} }( mS))
\overset{b_m}{\rightarrow} \bigoplus_{j=1}^N \bigoplus_{k=0}^m
{\mathcal{O} }_{\mathcal{B} }\xi_j^{-k}
\displaystyle{\frac{d}{d\xi_j}}
\overset{\vartheta_m}{\rightarrow} R^1\pi_* \Theta_{\mathcal{C} /\mathcal{B} }(-S)
 \rightarrow 0
$$
From which we deduce the following exact sequence of
$\Ob$-modules
$$
 0 \rightarrow \pi_*(\Theta_{\mathcal{C} /\mathcal{B} }( * S))
\overset{b}{\rightarrow} \bigoplus_{j=1}^N {\mathcal{O} }_{\mathcal{B} }[\,\xi_j^{-1}\,]
\displaystyle{\frac{d}{d\xi_j}}
\overset\vartheta{\rightarrow} R^1\pi_* \Theta_{\mathcal{C} /\mathcal{B} }(-S)
 \rightarrow 0.
$$
Note that the mappings  $b$ and $b_m$ correspond to the Laurent expansions
with respect to $\xi_j$ up to zero-th order.

By the Kodaira-Spencer theory (see \cite{Ue2} for details) we have the
following commutative diagram.
$$\begin{matrix}
 0 \rightarrow & \pi_* \Theta_{\mathcal{C} /\mathcal{B} }(*S)  & \rightarrow  &
\pi_*  \Theta'_{\mathcal{C} }(*S)_\pi & \overset{d\pi}{\rightarrow }&
\Theta_{\mathcal{B} }  &  \rightarrow  0 \\
& & & & &&\\
& \Vert & \phantom{p} & \;\;\downarrow \;p &
\phantom{\rho} &\;\; \downarrow \;\rho &  \\
& & & & &&\\
0 \rightarrow & \pi_* \Theta_{\mathcal{C} /\mathcal{B} }(*S)  & \rightarrow &
\bigoplus_{j=1}^N \mathcal{O} _\mathcal{B} [\; \xi_j^{-1}] \dfrac d{d \xi_j}
 & \overset{\vartheta}{\rightarrow }&
R^1\pi_* \Theta_{\mathcal{C} /\mathcal{B} }( -S)   &  \rightarrow  0
\end{matrix}
$$
where $\rho$ is the Kodaira-Spencer mapping
of the family $\mathfrak{F}^{(0)}$ and $p$ is
given by taking the non-positive part of
the $\dfrac d {d\xi_j}$ part of the Laurent expansions
of the vector fields in $\pi_*\Theta_\mathcal{C} (mS)_\pi$
at $s_j(\mathcal{B})$ . Since our family  $\mathfrak{F}^{(0)}$
is versal, it follows from Proposition 1.2.6. in \cite{Ue2}, that
the Kodaira-Spencer mapping $\rho$ is  an
isomorphism of $\mathcal{O}_\mathcal{B} $-modules. Therefore, $p$
is an isomorphism.   Put
$$
\mathcal{L}(\mathfrak{F}) :=
\bigoplus_{j=1}^N \mathcal{O} _\mathcal{B} [\; \xi_j^{-1}] \frac d{d \xi_j}.
$$
Then, we have the following exact sequence
\begin{equation}
0 \rightarrow \pi_* \Theta_{\mathcal{C} /\mathcal{B} }(*S)  \rightarrow
\mathcal{L}(\mathfrak{F} )  \overset{\theta}{\rightarrow }
\Theta_{\mathcal{B} }    \rightarrow  0 \label{Lfs}
\end{equation}
of $\mathcal{O} _\mathcal{B} $-modules.  The Lie
bracket  $[\phantom{X}, \phantom{X}]_d$ on
$\mathcal{L}(\mathfrak{F} )$ is determined by requiring $p$ to be
a Lie algebra isomorphism.
Thus, for $\vec \ell$, $\vec m \in \mathcal{L}(\mathfrak{F})$ we have
\begin{equation}
\label{35f}
[\vec \ell, \vec m]_d = [\vec \ell, \vec m]_0 +
\theta(\vec \ell)(\vec m) - \theta(\vec m)(\vec \ell)
\end{equation}
where $[\phantom{X},\phantom{X}]_0$ is the usual bracket of
formal vector fields and the action of $\theta(\vec \ell)$ on
$$
\vec m = (m_1{d\over{d\xi_1}}, \ldots,
m_N{d\over{d\xi_N}})
$$
is defined by
$$
      \theta(\vec \ell)(\vec m) = ( \theta(\vec \ell)(m_1){d\over{d\xi_1}}, \ldots,
           \theta(\vec \ell)(m_N){d\over{d\xi_N}} ).
$$
Then, the exact sequence (\ref{Lfs}) is an exact sequence of
sheaves of Lie algebras.

For $\vec \ell = (\underline l_1,
\ldots, \underline l_N) \in \mathcal{L}(\goF)$,
the action $D(\vec \ell)$ on $\mathcal{H}_\veclam(\mathcal{B} )$ is defined by
$$
     D(\vec \ell) (F \otimes |\Phi \rangle)
    = \theta(\vec \ell)(F) \otimes |\Phi\rangle -
  F \cdot (\sum_{j=1}^N \rho_j(T[\underline l_j] ) |\Phi\rangle
$$
where
$$
    F \in {\mathcal{O} }_{\mathcal{B} }, \quad |\Phi\rangle \in {\mathcal{H}}_{\veclam}.
$$

We have the following propositions.

\begin{proposition}[{[19, Proposition 4.2.2]}]
The action  $D(\vec \ell)$
 of  $\vec \ell \in \mathcal{L}(\goF)$ on $\Hlam(\mathcal{B} )$ defined  above
has the
following properties.

{\rm 1)} \quad For any $f \in {\mathcal{O} }_{\mathcal{B} }$ we have
$$
 D(f\vec \ell) = f D(\vec \ell).
$$

{\rm 2)} \quad For $\vec \ell$, $\vec  m \in  \mathcal{L}(\goF)$ we have
$$
[\,D(\vec \ell), D(\vec m)\,] = D([\,\vec \ell, \vec m\,]_d)
 +\frac{c_v}{12} \sumjn \Res_{\xi_j=0}
\left( \frac{d^3 \ell_j}{d \xi_j^3} m_j d\xi_j \right)\cdot {\text id}.
$$

{\rm 3)} \quad For $f \in {\mathcal{O} }_{\mathcal{B} }$ and $|\phi\rangle \in
\Hlam(\mathcal{B} )$ we have
$$
  D(\vec \ell)(f|\phi\rangle) = (\theta(\vec \ell)(f))|\phi\rangle
 + f D(\vec \ell)|\phi \rangle.
$$
Namely, $D(\vec \ell)$ is a first order  differential operator, if
$\theta(\vec \ell ) \neq 0$.
\end{proposition}

We define the dual action of $\mathcal{L}(\goF)$ on $\Hdaglam(\mathcal{B} )$
by
$$
  D(\vec \ell) (F \otimes \langle \Psi|) =
   (\theta(\vec \ell)F) \otimes \langle \Psi|  +
   \sumjn F  \cdot \langle \Psi| \rho_j(T[\underline l_j]).
$$
where
$$
 F \in \Ob, \qquad \langle \Psi | \in \Hdaglam(\mathcal{B}).
$$
Then, for any $|\widetilde \Phi\rangle \in \Hlam(\mathcal{B})$ and
$\langle \widetilde\Psi| \in \Hdaglam(\mathcal{B})$, we have
$$
 \{D(\vec \ell) \langle \widetilde\Psi| \}|\widetilde\Phi\rangle +
 \langle \widetilde\Psi|\{D(\vec \ell)|\widetilde\Phi\rangle\}
=  \theta(\vec \ell)\langle \widetilde\Psi|\widetilde\Phi\rangle .
$$

\begin{proposition}[{[19, Proposition 4.2.3]}]
For any $\vec \ell \in \mathcal{L}(\goF)$
we have
$$
D(\vec \ell)(\widehat{\mathfrak{g} }(\mathfrak{F})\Hlam(\mathcal{B} )) \subset
   \widehat{\mathfrak{g} }(\mathfrak{F} )\Hlam(\mathcal{B} ).
$$
Hence, $D(\vec \ell)$ operates on $\Vlam(\goF)$. Moreover, it is a
first order differential operator, if $\theta(\vec \ell) \neq 0$.
\end{proposition}

\begin{proposition}[{[19, Proposition 4.2.8]}]
For each element
$\vec \ell \in
\mathcal{L}(\goF)$, $D(\vec \ell)$ acts on $\Vdaglam(\goF)$.  Moreover, if
$\theta(\vec \ell) \neq 0$, then $D(\vec \ell)$ acts on  $\Vdaglam(\goF)$
as a first order differential operator.
\end{proposition}

Note that for the natural bilinear pairing $\langle \phantom{X}|\phantom{X}\rangle
: \Vdaglam(\goF) \times \Vlam(\goF) \rightarrow \Ob$, we have the equality
$$
 \{ D(\vec \ell) \langle\Psi |\}|\Phi\rangle +
  \langle\Psi |\{D(\vec \ell)|\Phi\rangle  \}= \theta (\vec \ell)
 (\langle\Psi |\Phi\rangle).
$$

\subsection{The projectively flat connections on the on the bundle of vacua}

Let $\Si$ be a closed oriented smooth surface and let $P$ be a finite set
of marked points on $\Si$. Assume that $(\Si,P)$ is \stable
and \fm pointed surface.

 Let
$\goF = ( \pi \: \mathcal C \rightarrow \mathcal B; \vs; \veta )$
be a family of pointed Riemann Surfaces with formal
neighbourhoods on $(\Si, P)$. Recall the discussion of symmetric
bidifferentials from section 1.4 in \cite{Ue2}. We now introduce
the notion of a normalized symmetric bidifferential.

\begin{definition}\label{nsymbidif}
A symmetric bidifferential
$\omega \in H^0(\mathcal C \times_{\mathcal B} \mathcal C, \omega_{\mathcal C
\times_{\mathcal B} \mathcal C/\mathcal B}(2\Delta))$ with
$$\omega= \left(\frac{1}{(x-y)^2} + \text{holomorphic}\right)dxdy$$
in a neighbourhood of the diagonal of ${\mathcal C \times_{\mathcal B} \mathcal C}$ is called a
{\em normalized} symmetric bidifferential for the family $\goF$.
\end{definition}

For a Riemann surface $R$ we let $(\vec{\alpha}, \vec{\beta})= (\alpha_1, \ldots, \alpha_g,\beta_1,\ldots,\beta_g)$ be
a symplectic basis of
$H_1(R,\Z{})$. We can find a basis $\{\omega_1, \ldots, \omega_g\}$ of
holomorphic one forms of $R$ with
\begin{equation}
\label{betaone}
\int_{\beta_i}\omega_j = \delta_{i j}, \quad 1 \le i,j \le g.
\end{equation}
The matrix
\begin{equation}\label{period}
 \tau = (\tau_{ij}), \quad \tau_{ij} =
\int_{\alpha_i}\omega_j
\end{equation}
is then called the period matrix of the Riemann surface $R$.
The complex torus
$$
J(R) = {\mathbb C}^g/(\tau, I_g)
$$
is called a Jacobian variety. If we chose a point $P$ on $R$ we can define a holomorphic
mapping
\begin{alignat*}{2}
j : &R  &\quad \rightarrow & \quad J(R) \\
   & Q & \quad \mapsto &\left( \int_P^Q \omega_1, \ldots, \int_P^Q \omega_g \right).
\end{alignat*}
If a family of pointed Riemann surfaces is given, we can construct a family of Jacobian varieties
and a family of holomorphic mappings.

We have the following lemma as a consequence of the construction in
Section 1.4 in \cite{Ue2}.

\begin{lemma}\label{nsymbidifexists}
For any family of pointed Riemann surfaces with formal neighbourhoods
 $\goF$ on $(\Si,P)$ and any symplectic basis
$(\vec{\alpha}, \vec{\beta})=(\alpha_1, \ldots, \alpha_g,\beta_1,\ldots,\beta_g)$ of
$H_1(\Si, \Z{})$, there is a unique normalized symmetric bidifferential
$\omega \in H^0(\mathcal C \times_{\mathcal B} \mathcal C,
\omega_{\mathcal C \times_{\mathcal B} \mathcal C/\mathcal
B}(2\Delta))$ determined by formula
\begin{equation}\label{primeform}
\omega(x,y)dxdy = \frac{\partial^2 \log E(x,y)}{\partial x \partial y}dx dy
\end{equation}
where $E(x,y)(\sqrt{dx})^{-1}(\sqrt{dy})^{-1}$ is the prime form associated to the
symplectic basis $(\vec{\alpha}, \vec{\beta})$ for each Riemann surface.
\end{lemma}

For a prime form see Chapter II of \cite{Fa}.
Please do note that $\alpha$ and $\beta$ play the reverse roles
in \cite{Ue2}, but the same as in \cite{AU1}.

Note that a prime form of a Riemann surface $R$ (hence, also a normalized
symmetric bidifferential) is uniquely determined by a symplectic
basis $(\vec{\alpha}, \vec{\beta})$ of $H_1(R,\Z{})$ .
If $(\vec{\widehat{\alpha}}, \vec{\widehat{\beta}})=
(\widehat{\alpha}_1, \ldots, \widehat{\alpha}_g, \widehat{\beta}_1,\ldots,\widehat{\beta}_g)$
is another symplectic basis , there exists a symplectic matrix
$$
\Lambda = \left( \begin{array}{cc} A&B\\C&D \end{array} \right) \in Sp(g, \Z{})
$$
such that
\begin{equation} \label{action1}
\left(\begin{array}{c} \widehat{\alpha}_1\\ \vdots\\ \widehat{\alpha}_g \\ \widehat{\beta}_1\\
\vdots \\ \widehat{\beta}_g \end{array}\right)=
\left(\begin{array}{cc}A & B\\C & D\end{array}\right)
 \left(\begin{array}{c} \alpha_1\\ \vdots\\ \alpha_g \\ \beta_1\\
\vdots \\ \beta_g \end{array}\right).
\end{equation}
Also for any element $\Lambda \in Sp(g, \Z{})$,  by \eqref{action1}
we can define a new symplectic basis $\Lambda (\vec{\alpha},\vec{\beta})
= (\widehat{\alpha}_1, \ldots, \widehat{\alpha}_g, \widehat{\beta}_1,
\ldots, \widehat{\beta}_g)$. Then the normalized symmetric bidifferential
$\widehat{\omega}(x,y)dx dy$ associated to the symplectic basis $\Lambda (\vec{\alpha},\vec{\beta})$
and the normalized symmetric bidifferential $\omega(x,y)dx dy$ have a relation:
\begin{eqnarray} \label{bidiffrelation}
\widehat{\omega}(x,y)dxdy &=& \omega(x,y)dxdy \\
&& -\frac12 \sum_{i\le j} \left\{ v_i(x) v_j(y) +v_j(x)v_i(y)\right\}
\frac{\partial}{\partial \tau_{ij}}\log \det(C \tau +D)
\nonumber
\end{eqnarray}
where $\{\omega_1 =v_1(x)dx, \ldots ,\omega_g=v_g(x)dx\}$ is a basis of holomorphic one-forms on the
Riemann surface $R$ which satisfies \eqref{betaone}, $\tau_{ij}$ is defined by \eqref{period} and
$$
\Lambda = \left( \begin{array}{cc} A&B\\C&D \end{array} \right) .
$$
For details see \cite{Fa}, Chapter II.

Let $\omega$ be a normalized symmetric bidifferential for $\goF$.
In the formal neighbourhood $\eta_j$ we define the quadratic
differential ({\em projective connection}\/)
\begin{equation} \label{projectiveconnection}
S_{\omega}d\eta_j^2 = 6 \lim_{\overset{\eta \ra \eta_j}{\xi \ra \eta_j}}
\left( \omega(\eta,\xi) - \frac{d\eta d\xi}{(\eta -\xi)^2}\right) .
\end{equation}

We define
$$
a_\omega :
\mathcal{L}(\goF) \rightarrow \Ob
$$
as an $\Ob$-module homomorphism by
\[a_\omega(\vec \ell) = - \frac{c_v}{12} \sum_{j=1}^N \Res_{\eta_j =0} (\ell_j(\eta_j)
S_\omega(\eta_j) d\eta_j)\]
for all $\vec \ell \in \mathcal{L}(\goF)$.

\begin{definition}\label{conomega}
For each $b\in  \mathcal B$ and each element $X \in (\Theta_{\mathcal{B}})_b$,
there is an
element $\vec \ell \in  \mathcal{L}(\goF)_b$ with $\theta(\vec \ell) = X$. Define
an operator $\nabla_X^{(\omega)}$ acting on $\Vdaglam(\goF)_b$ from the left by
$$
\nabla_X^{(\omega)}(\langle \Phi|) = D(\vec \ell)([\langle\Phi|)
   + a_\omega(\vec \ell)([\langle \Phi|).
$$
\end{definition}

\begin{theorem}\label{holomorphicconnection}
$\nabla^{(\omega)}$ is a well-defined holomorphic connection in
$\Vdaglam(\goF)$.
\end{theorem}

In \cite{Ue2} we introduced the holomorphic connection on $\Vlam(\goF)_b$ and
our connection $\nabla^{(\omega)}$ is its dual connection (see Proposition 5.1.4 in \cite{Ue2}).

\begin{definition}
We define a connection in the vector bundle
$\Vdaglam(\goF)$ over $\mathcal B$ by letting its $(1,0)$-part be
given by the holomorphic connection $\nabla^{(\omega)}$ just
defined and its $(0,1)$-part be given by the $\bar \partial$-operator
determined by the holomorphic structure on $\Vdaglam(\goF)$. We also
denote this connection $\nabla^{(\omega)}$.
\end{definition}

\begin{theorem}\label{curvaturecom}
The curvature of the connection $\nabla^{(\omega)}$ is given by
\[R^\omega(X,Y) = \left\{ -a_\omega(\vec{n}) +  X(a_\omega(\vec{m})) -
Y(a_\omega(\vec{l})) - \frac{c_v}{12}\sum_{j=1}^N
\Res_{\xi_j=0}\left( \frac{d^3 l_j}{d\xi_j^3} m_jd\xi_j \right)\right\}
\otimes \id\]
where $\vec \ell$  and $\vec m$ are liftings of $X$ and $Y$ to
$\mathcal{L}(\goF)$, i.e.
$\theta(\vec \ell) = X$ and $\theta(\vec m) = Y$, and $\vec{n} =
[\vec{l}, \vec{m}]_d$ (see \eqref{35f}).
Hence we see that the connection is projectively flat and the
curvature is of type $(2,0)$.
\end{theorem}

\proof
It follows from the definition of the connection
in the above definition that the
$(1,1)$ and $(0,2)$-part of the curvature vanishes.

\eproof

\subsection{The definition of the connection in the bundle of vacua over Teichm\"{u}ller space}\label{New5.3}

Suppose we have two \good families $\goF_i$, $i=1,2$ of pointed Riemann Surfaces with for mal neighbourhoods,
with
the property that they have the same image $\Psi_{\goF_1}
({\mathcal B}_1) = \Psi_{\goF_2}({\mathcal B}_2)$ in Teichm\"{u}ller
space $\cT_{(\Si,P)}$.

\begin{lemma}\label{contransf}
Let $\nabla^{(\omega)}_i$ be the connection in $\Vdaglam(\goF_i)$. Then we have that
\[G_{12}^*(\nabla^{(\omega)}_2) = \nabla^{(\omega)}_1.\]
\end{lemma}

\proof Since the connection is descended from the $\mathcal H$-level and $G_{12}$ is also descended from this level, we just need
to check the transformation rule on this level. Up on the $\mathcal H$-level it follows
straight from Theorem \ref{thm6.1} (3).
 \eproof

\begin{theorem}\label{conTeich}
Let $(\Si,P,\l)$ be a closed oriented stable and \fm marked
surface and let $(\vec \alpha, \vec \beta) = (\alpha_1, \ldots,
\alpha_g,\beta_1,\ldots,\beta_g)$ be a symplectic basis of
$H_1(\Si, \Z{})$. There is a unique connection $\nabla^{(\vec
\alpha, \vec \beta)} = \nabla^{(\vec
\alpha, \vec \beta)}(\Si,P)$ in the bundle $\Spofvlam(\Si,P)$ over
$\cT_{(\Si,P)}$ with the property that for any
\good family
$\goF$ of stable pointed Riemann surfaces with formal neighbourhoods on $(\Si,P)$ we have that
\[\Psi_{\goF}^*(\nabla^{(\vec \alpha, \vec \beta)}) = \nabla^{(\omega)}.\]
In particular the connection is holomorphic and projectively flat with $(2,0)$-curvature
as described in Theorem \ref{curvaturecom}.
If we act on the symplectic basis $(\vec \alpha,\vec \beta)$ by an element $\Lambda =
\left(\begin{array}{cc}A & B\\C & D\end{array}\right) \in
\text{Sp} (g,\Z{})$ by \eqref{action1}, then we have
\begin{equation}
\nabla^{\Lambda(\vec \alpha, \vec \beta)} - \nabla^{(\vec \alpha, \vec \beta)} =
- \frac{c_v}{2} \Pi^*( d \log \det (C \tau + D )),\label{contchbasis}
\end{equation}
where $\Pi$
is the period mapping of holomorphic one-forms form the base space of  $\mathfrak{F}$ to the
Siegel upper-half plane of degree $g$. If $f : (\Si_1,P_1,\l_1) \ra
(\Si_2,P_2,\l_1)$ is an orientation preserving diffeomorphism of
labeled pointed surfaces which maps the symplectic basis $(\vec
\alpha^{(1)},\vec \beta^{(1)})$ of $H_1(\Si_1,\Z{})$ to the
symplectic basis $(\vec \alpha^{(2)},\vec \beta^{(2)})$ of
$H_1(\Si_2,\Z{})$ then we have that
\[\Spofv(f)^*(\nabla^{(\vec \alpha^{(2)},\vec \beta^{(2)})}) = \nabla^{(\vec \alpha^{(1)},
\vec \beta^{(1)})}.\]
\end{theorem}

\proof The existence of the connection is a consequence of Lemma \ref{contransf}.
The transformation law (\ref{contchbasis}) is proved in section
5.2 in \cite{Ue2}.\eproof

\begin{proposition}\label{Teichpullbackcon}
For any \stable and \fm family $\tgoF$ of pointed Riemann surfaces with
formal neighbourhoods on $(\Si,P)$ the preferred isomorphism
\[\Upsilon_{\tgoF} : \Vdaglam(\tgoF) \ra \Psi_{\tgoF}^*\Spofvlam(\Si,P)\]
given in Proposition \ref{Teichpullback} preserves connections.
\end{proposition}

This follows directly from Lemma \ref{contransf}. \eproof

\begin{proposition}\label{Rplusact}
The $\RPP$-action on $\cT_{(\Si,P)}$ lifts by the use of the
connection $\nabla^{(\vec
\alpha, \vec \beta)}$ to an action $\Spofvlam$ of $\RPP$ on
$\Spofvlam(\Si,P)$ which preserves the connection $\nabla^{(\vec
\alpha, \vec \beta)}$.
\end{proposition}

\proof Since the $\RPP$ action on the formal coordinates is just obtained by scaling
the coordinates by positive scalars, we get a well defined homomorphism from
$\RPP$ to the group of formal coordinates changes for any family of $N$-pointed Riemann surfaces with
formal coordinates, hence by composing with the group homomorphism $G$ we get an action
of $\RPP$ on $\Spofvlam(\Si,P)$. Note that we have here used Theorem 3.2.4.
(2) for $p=0$ in \cite{Ue2}, but
only for these special real coordinates changes. This action preserved the connection
by Theorem 3.2.4. (3) in \cite{Ue2}.
\eproof

Let now $\Sib_i$, $i=1,2$ be marked surface and let $f_j : (\Si_1,P_1) \ra
(\Si_2,P_2)$, $j=1,2$ be diffeomorphisms of pointed surfaces,
which induce the same morphism of marked surfaces from $\Sib_1$ to
$\Sib_2$. Then there exists a unique $v\in \RPP$ such that $v \cdot d_{P_1}f_1 =
d_{P_2}f_2$.

\begin{lemma}\label{momswd}
We have a commutative diagram
\begin{equation*}
\begin{CD}
\Spofvlamone(\Si_1,P_1) @> \Spofv(f_1) >> \Spofvlamtwo(\Si_2,P_2)\\
@V = VV               @VV\Spofvlamtwo(v)V\\
\Spofvlamone(\Si_1,P_1) @> \Spofv(f_2) >> \Spofvlamtwo(\Si_2,P_2),
\end{CD}
\end{equation*}
\end{lemma}

\proof By the construction of $\Spofvlami(\Si_i,P_i)$ we just need to
check the commutativity
on the $\mathcal H$-level, where this just amounts to $G$ being a homomorphism,
which again is the content of Theorem 3.2.4.
(2) in \cite{Ue2}.
 \eproof

Let $(\Si, P,\l)$ be a general labeled pointed surface which might
not be \stable or \fmp Let now $(\tilde{P},\tilde{\l})$ be
obtained from $(P,\l)$ by labeling further points not in $P$ by
$0\in P_\ell$ such that $(\tilde{\Si},\tilde{P},\tilde{\l})$ is
a \stable and \fm labeled pointed surface. Let $\tilde{\pi} : \cT_{(\Si,\tilde{P})}\ra
\cT_{(\Si,P)}$ be the natural projection map.

\begin{proposition}\label{propvaccon}
The connection $\nabla^{(\vec
\alpha, \vec \beta)} = \nabla^{(\vec
\alpha, \vec \beta)}(\Si,P)$ is flat with trivial holonomy when
restricted to any of the fibers of $\tilde{\pi} : \cT_{(\Si,\tilde{P})}\ra
\cT_{(\Si,P)}$.
\end{proposition}

\proof The propagation of vacua isomorphism applied along the
fibers of $\tilde \pi$ is compatible with the connection.
In the notation of chapter 5 of \cite{Ue2} this is seen
by the following calculation. Let $\tilde N = |\tilde P|$.
Let $\goF$ be a family of $\tilde N$-pointed curves
with formal neighbourhoods on $(\Sigma,\tilde
P)$ such that $\Psi_{\goF}(\mathcal B)$ is contained in a fiber of
$\tilde \pi$. Further we can assume that when we consider $\goF$
over $(\Sigma,P)$, then it is simply a product family of a fixed
$N$-pointed curve with formal neighbourhoods crossed with the base
$\mathcal B$, i.e. we only vary the coordinates and the points in
$\tilde P - P$. Then for any tangent field $X$ on $\mathcal B$, we
choose a corresponding $\vec \ell = (\ell_1, \ldots, \ell_N, \ell_{N+1}, \ldots, \ell_{\tilde
N})\in \mathcal L(\goF)$,
such that $\ell_j = 0$, $j=1, \ldots, N$ and $\ell_j \in \C{}[[\xi]]$, $j=N+1, \ldots, \tilde N$.
Then for a section of $\mathcal
H_{\tilde \lambda}(\mathcal B)$ of the form $F\otimes|\Phi\rangle\otimes |0\rangle\otimes \ldots \otimes |0\rangle$
we
compute that
\begin{eqnarray*}
\nabla^{(\omega)}_X(F |\Phi\rangle\otimes |0\rangle\otimes \ldots \otimes |0\rangle) &=&
X(F) |\Phi\rangle\otimes |0\rangle\otimes \ldots \otimes |0\rangle\\
&& - F \sum_{j=N+1}^{\tilde N} \rho_{j}(T[\ell_j]) (|\Phi\rangle\otimes |0\rangle\otimes
\ldots \otimes |0\rangle) \\
&& - a_\omega(\vec \l)F |\Phi\rangle\otimes |0\rangle\otimes \ldots \otimes
|0\rangle\\
&=&
X(F) |\Phi\rangle)\otimes |0\rangle\otimes \ldots \otimes
|0\rangle.
\end{eqnarray*}
Here we have used that $\rho_{j}(T[\ell_j])(|0\rangle)=0$ because $\ell_j\in \C{}[[\xi]]$
and $a_\omega(\vec \l)=0$ for the same reason.
Hence we see that $\nabla^{(\omega)}$ is the trivial connection
in $\mathcal H_{\tilde \lambda}(\mathcal B)$. Hence the connection
is flat with trivial holonomy on the subbundle $\Vdagtlam(\goF) \subset \Hdagtlam(\mathcal B)$.
\eproof

Let $(\Si,P,\lambda)$ be a general labeled pointed surface, i.e.
$(\Si,P)$ might not be \stable nor \fmp Let
$(\Si,P',\lambda')$ and $(\Si,P'',\lambda'')$ be
labeled marked surfaces obtained from $(\Si,P,\lambda)$ by
labeling further points not in $P$ by $0\in P_\ell$.
Assume that $(\Si,P')$ and $(\Si,P'')$ are \stable and \fm pointed surfaces.
Let $\bar P = P' \cup
P''$ and $\bar \l$ be the induced labeling of $\bar P$. Note that
$(\Si,\bar P)$ is also \stable and \fmp

\begin{proposition}\label{pullbackcon}
Let $(\vec \alpha, \vec \beta)$ be a symplectic basis of $H_1(\Si,\Z{})$. The
isomorphisms given in Theorem \ref{propvaciso} satisfies
\[ \nabla^{(\vec
\alpha, \vec \beta)}(\Si,\bar P) = (\pi')^*\nabla^{(\vec
\alpha, \vec \beta)}(\Si,P') = (\pi'')^*\nabla^{(\vec
\alpha, \vec \beta)}(\Si,P'').\]
\end{proposition}

\proof
We only have to consider the case of adding one point to a \stable
and \fm curve, i.e. say $\bar P = P'$ and $P'$ is obtained from $P''$ by adding
one more point. Let $\bar N$ be the number of points in $\bar P$ and
$\goF$ be a family of $\bar N$-pointed curves
with formal neighbourhoods on $(\Sigma, \bar P)$. For any tangent
field $X$ on $\mathcal B$, we choose a corresponding $\vec \ell =
(\ell_1, \ldots, \ell_{\bar N}, \ell_{\bar N+1})\in \mathcal L(\goF)$,
such that  $\ell_{\bar N+1} \in \C{}[[\xi]]$. Then the same computation as
above shows that
\[\nabla^{(\omega)}(F |\Phi\rangle\otimes |0\rangle ) =
\nabla^{(\omega)}(F |\Phi\rangle)\otimes |0\rangle.\]
The Proposition follows directly from this.
\eproof

Let now $(\Sigma_i,P'_i,\lambda'_i)$ and
$(\Sigma_i,P''_i,\lambda''_i)$ be \stable and \fm labeled pointed surfaces,
obtained from the labeled pointed surfaces
$(\Sigma_i,P_i,\lambda_i)$, by labeling further points with the
zero-label. Assume that $f' : (\Si_1,P_1',\l_1') \ra (\Si_2,P_2',\l_2')$ and
$f'' : (\Si_1,P_1'',\l_1'') \ra (\Si_2,P_2'',\l_2'')$ are
diffeomorphisms which induce isotopic maps from $(\Si_1,P_1)$ to
$(\Si_2,P_2)$, where the isotopy is through maps which induces the
same map from $PT_{P_1}\Si_1$ to $PT_{P_2}\Si_2$.

\begin{proposition}\label{propvmorph}
With respect to the propagation of vacua isomorphisms, we
get that
\[(\pi')^* \Spofv(f') = (\pi'')^* \Spofv(f'').\]
\end{proposition}

This follows directly from the way the morphisms for $f'$ and $f''$ are defined.

\section{Definition of the space of \abelian vacua associated to a Riemann surface.}\label{New3ab}

\subsection{Fermion Fock space}

Let $\hZ$ be the set of all half integers. Namely
$$
\hZ = \{ n + 1/2\, | \, n \in \bZ \,\}.
$$
Let $\cWd$ be an infinite-dimensional vector space over $\bC$
with a filtration $\{F^m\cWd\}_{m \in \bZ}$ which satisfies the following conditions.
\begin{enumerate}
\item  The filtration $\{F^m\cWd\}$ is decreasing;
\item  $\bigcup_{m \in \bZ} F^m\cWd = \cWd$, \quad
$\bigcap_{m \in \bZ} F^m\cWd = \{0\}$;
\item $\dim_\bC F^m\cWd/ F^{m+1}\cWd = 1$;
\item The vector space $\cWd$ is complete with respect to the uniform topology such that
$\{F^m\cWd\}$ is a basis of open neighbourhoods of 0.
\end{enumerate}
We introduce a  basis $\{e^\nu \}_{\nu \in \hZ}$ of $\cWd$ in such a way that
$$
e^{m+1/2} \in F^m\cWd \setminus F^{m+1}\cWd.
$$
Then, each element $ u \in \cWd$ can uniquely be expressed in the form
$$
u= \sum_{\nu> n_0, \nu \in \hZ}^\infty a_\nu  e^\nu
$$
for some $n_0$ and with respect to this basis the filtration is given  by
$$
F^m\cWd = \left\{ u \in \cWd \, \left|
\, u=  \sum_{\nu> m, \nu \in \hZ}^\infty a_\nu  e^\nu\, \right. \right\}.
$$
We fix the basis $\{e^\nu\}_{\nu \in \hZ}$ throughout the present
paper.

Let $\bC((\xi))$ be a field of formal Laurent series over the complex number field.
Then the basis gives us a filtration preserving linear isomorphism
\begin{eqnarray*}
\bC((\xi)) & \cong & \cWd \\
\xi^n  &\mapsto &e^{n+1/2}.
\end{eqnarray*}
By mapping $\xi^n d\xi$ to $e^{n+1/2}$ we of course also get a filtration
preserving linear isomorphism between $\bC((\xi))d\xi$ and $\cWd$.

We let $\{\oe_\nu\}_{\nu \in \hZ}$ be the dual basis of $\{e^\nu\}_{\nu \in \hZ}$.
Then,  put
$$
\cW = \bigoplus_{\nu \in \hZ} \bC\oe_\nu
$$
Then $\cW$ is the topological dual of the vector space $\cWd$.
There is a natural pairing $(\phantom{X}|\phantom{X}) : \cWd\times \cW \rightarrow \bC$
defined by
$$
(e^\nu|\oe_\mu) = \delta^\nu_\mu.
$$
In other word we have
$$
(u|v) = v(u).
$$

Let us introduce the semi-infinite exterior product of the vector spaces
$\cW$ and $\cWd$.  For that purpose we first introduce the notion of
a Maya diagram.

\begin{definition}{\rm
A Maya diagram $M$ of the charge $p$, $p \in \bZ$ is a set
$$
M = \left\{ \mu(p-1/2), \mu(p-3/2), \mu(p-5/2), \ldots\right\},
$$
where $\mu$ is an increasing function
$$
\mu : \hZ_{<p} = \{ \nu \in \hZ\, | \, \nu < p\,\} \rightarrow \hZ
$$
such that there exists an integer $n_0$ such that
$$
\mu(\nu) = \nu
$$
for all $\nu < n_0$.

The function $\mu$ is called the {\it characteristic function} of the
Maya diagram $M$.
The set of Maya diagrams of charge $p$ is written as ${\mathcal M}_p$.
}
\end{definition}

For a Maya diagram M we have $\mu(\nu) = \nu$ for almost all $\nu$. Therefore
the set
$$
\{ \mu(\nu) - \nu \,|\, \nu \in \hZ, \, \mu(\nu) - \nu >0\,\}
$$
is finite and the number
$$
d(M) = \sum_{\nu \in \hZ} (\mu(\nu) - \nu)
$$
is finite. The number $d(M)$ is also written as $d(\mu)$ and it is
called the {\it degree} of the Maya diagram $M$ with characteristic function  $\mu$.
The finite set of Maya diagrams of degree $d$ and change $p$ is
denoted ${\mathcal M}_p^d$. Clearly ${\mathcal M}_p = \coprod_{d} {\mathcal
M}_p^d$.

For a Maya diagram $M$ of charge $p$ we define two semi-infinite
products
\begin{eqnarray*}
|M\rangle &= & \oe_{\mu(p-1/2)}\wedge \oe_{\mu(p-3/2)} \wedge \oe_{\mu(p-5/2)}
\wedge \cdots \\
\langle M| &=& \cdots \wedge
 e^{\mu(p-5/2)}\wedge  e^{\mu(p-3/2)} \wedge  e^{\mu(p-1/2)}
\end{eqnarray*}
Formally, these semi-infinite products is just another notation
for the corresponding Maya diagram. This notation is
particular convenient for the following discussion. However, by using
the basis $e_\nu$, we clear indicate the relation to the vector
spaces $\cW$ and $\cWd$.

For any integer $p$ put
\begin{eqnarray*}
|p \rangle &= & \oe_{p-1/2}\wedge \oe_{p-3/2} \wedge \oe_{p-5/2}
\wedge \cdots \\
\langle  p| &=& \cdots \wedge
 e^{p-5/2}\wedge  e^{p-3/2} \wedge  e^{p-1/2}
\end{eqnarray*}

Now the {\it fermion Fock space} $\cFd(p)$
of {\it charge} $p$ and the {\it dual fermion Fock space} $\cF(p)$ of
{\it charge} $p$ are defined by
\begin{eqnarray*}
\cF(p) &=&  \bigoplus_{M \in {\mathcal M}_p} \bC |M\rangle \\
\cFd(p) &=& \prod_{M \in {\mathcal M}_p} \bC \langle M|
\end{eqnarray*}
We observe that
\[\cF(p) = \bigoplus_{d\geq 0} \cF_d(p),\]
where
\[\cF_d(p) =  \bigoplus_{M\in {\mathcal M}_p^d} \bC |M\rangle. \]

The dual pairing
$$
\langle \cdot | \cdot \rangle : \cFd(p) \times \cF(p) \rightarrow \bC
$$
is given by
$$
\langle M| N \rangle = \delta_{M,N}, \quad M,  N \in {\mathcal M}_p
$$
Put also
\begin{eqnarray*}
\cF      &=&  \bigoplus_{p \in \bZ} \cF(p) \\
\cFd   &=& \bigoplus_{p \in \bZ} \cFd(p)
\end{eqnarray*}
The vector space $\cFd$ is called the {\it fermion Fock space} and $\cF$ is
called the {\it dual fermion Fock space}.
These are the
semi-infinite exterior products of the vector spaces $\cWd$ and
$\cW$ respectively, which we shall be interested in. We only
define the fermion Fock space by using the basis $e_\nu$,
since we are fixing this basis throughout.

The above pairing can be extended to the one on $\cFd \times \cF$ by
assuming that the paring is zero on $\cFd(p) \times \cF(p')$ if $p \ne p'$.

Let us introduce the {\it fermion operators} $\psi_\nu$ and $\ovpsi_\nu$  for
all half integers $\nu \in \hZ$ which act on $\cF$ from the left  and on
$\cFd$ from the right.
\begin{eqnarray}
\hbox{\rm Left action on $\cF$} &\quad & \psi_\nu = i(\oe_\nu), \quad
\ovpsi_\nu = \oe_{-\nu} \wedge \\
\hbox{\rm Right action on $\cFd$} &\quad & \psi_\nu = \wedge  e^\nu, \quad
\ovpsi_\nu = i(e^{-\nu})
\end{eqnarray}
where $i(\cdot)$ is the interior product. For example we have
\begin{eqnarray*}
\psi_{-3/2}|0\rangle &=& i(\oe_{-3/2}) \oe_{-1/2}\wedge \oe_{-3/2}\wedge \cdots =
- \oe_{-1/2}\wedge\oe_{-5/2}\wedge \oe_{-7/2}\wedge \cdots , \\
\langle 0 | \ovpsi_{5/2} &=& \cdots \wedge e^{-5/2}\wedge e^{-3/2} \wedge e^{-1/2} i(e^{-5/2})
=\cdots \wedge e^{-7/2}\wedge e^{-3/2} \wedge e^{-1/2}
\end{eqnarray*}
Note that $\psi_\nu$ maps $\cF(p)$ to $\cF(p-1)$, hence decreases the charge by one,
and $\ovpsi_\nu$ maps $\cF(p)$ to $\cF(p+1)$, hence increase the charge by one.
Similarly the right action of $\psi_\nu$ maps $\cFd(p)$ to $\cFd(p+1)$ and
$\ovpsi_\nu$ maps  $\cFd(p)$ to $\cFd(p-1)$. It is easy to show that
for any $\langle u| \in \cF$ and $|v\rangle \in \cFd$ we have
$$
\langle u | \psi_\nu v\rangle = \langle u  \psi_\nu | v\rangle,
\quad \langle u | \ovpsi_\nu v\rangle = \langle u  \ovpsi_\nu | v\rangle .
$$

The fermion operators have the following anti-commutation relations as operators
on $\cF$ and $\cFd$.
\begin{eqnarray}
{[\psi_\nu, \psi_\mu ]}_+  &=&  0, \label{anticomm1}\\
{[\ovpsi_{\nu}, \ovpsi_{\mu} ]}_+ &=& 0,  \label{anticomm2}\\
{[\psi_\nu, \ovpsi_{\mu} ]}_+  &=& \delta_{\nu +\mu, 0}, \label{anticomm3}
\end{eqnarray}
where we define
$$
{[A,B]}_+ = AB+BA.
$$

Note that for each Maya diagram $M$ of charge $p$ we can find
non-negative half integers
$$
\mu_1<\mu_2<\cdots<\mu_r<0, \quad
 \nu_1<\nu_2<\cdots<\nu_s<0, \quad r\ge 0, \,s\ge 0
 $$
with $r- s=p$ and $\mu_i \neq \nu_j$ such that
\begin{equation}
\label{maya}
 |M\rangle = (-1)^{\sum_{i=1}^s\nu_i  +s/2}
  \ovpsi_{\mu_1}\ovpsi_{\mu_2}\cdots\ovpsi_{\mu_r}
 \psi_{\nu_s}\psi_{\nu_{s-1}}\cdots\psi_{\nu_1}|0\rangle.
\end{equation}
The negative half integers $\mu_i$'s and $\nu_j$'s are uniquely determined by
the Maya diagram $M$.

The {\it normal ordering} $\normalord\phantom{X}\normalord$ of the
fermion operators are defined as follows.
$$
\normalord A_\nu B_\mu \normalord = \left\{
\begin{array}{ll}
- B_\mu A_\nu & \hbox{\rm if $\mu <0$ and $\nu>0$}, \\
A_\nu B_\mu & \hbox{\rm otherwise,}
\end{array}
\right.
$$
where $A$ and $B$ is $\psi$ or $\ovpsi$. By \eqref{anticomm1}, \eqref{anticomm2} and
\eqref{anticomm3} the normal ordering is non-trivial if and only if $\mu <0$ and
$A_\nu ={\bar B}_{-\mu} $.

The field operators $\psi(z)$ and $\ovpsi(z)$ are defined by
\begin{eqnarray}
\psi(z)&=& \sum_\mu \psi_\mu z^{-\mu -1/2} , \\
\overline{\psi}(z) &=& \sum_\mu \overline{\psi}_\mu z^{-\mu -1/2} .
\end{eqnarray}
The current operator $J(z)$ is defined by
\begin{equation}
J(z) = \normalord \overline{\psi}(z)\psi(z) \normalord = \sum_{n \in \bZ} J_n z^{-n-1}
\end{equation}

Note that thanks to the normal ordering, the operator $J_n$ can operate on $\cF$ and $\cFd$
even though $J_n$ is an infinite sum of  operators.

The energy-momentum tensor $T(z)$ is defined
by
\begin{equation}\label{emab}
T(z) = \normalord  \frac{d \psi(z)}{dz} \ovpsi(z)  \normalord = \sum_{n \in \bZ} L_n z^{-n-2}.
\end{equation}
Again due to the normal ordering, the coefficients $L_n$ operates on $\cF$ and $\cFd$.
The set $\{L_n\}_{n \in \bZ}$ forms  the {\it Virasoro} algebra with central charge $c=-2 $.

The field operators $\psi(z)$ and $\ovpsi(z)$, the current operator $J(z)$
and the energy-momentum tensor $T(z)$ form the so-called
spin $j=0$ $bc$-system or ghost system in the physics literature.

\subsection{Abelian Vacua}

For a positive integer $N$ put
\begin{eqnarray*}
\cF_N &=& \bigoplus_{p_1, \ldots,p_N \in\bZ} \cF(p_1) \otimes \cdots \otimes \cF(p_N) ,\\
\cFd_N &=& \bigoplus_{p_1, \ldots, p_N \in \bZ}
{\cFd}(p_1)  \hat{\otimes}  \cdots \hat{\otimes} {\cFd}(p_N),
\end{eqnarray*}
where $\hat{\otimes}$ means the complete tensor product.

\begin{definition}
\label{dfn3.1}  Let ${\mathfrak X} = (C; q_1,q_2,
\ldots, q_N; \xi_1, \xi_2, \ldots, \xi_N)$ be a \fm pointed
Riemann Surface with formal neighbourhoods .
The abelian vacua (ghost vacua in [2])  $\Vdagab(\gX)$
of the spin $j=0$ ghost system associated to $\mathfrak{X}$ is a linear subspace of
$\cFd_N$ consisting of elements $\langle \Phi|$ satisfying the following conditions:
\begin{enumerate}
\item  For all $|v\rangle \in \cF_N$, there exists a meromorphic function
$f \in H^0(C,\cO_C(*\sum_{j=1}^Nq_j))$
such that
$\langle \Phi| \rho_j(\psi(\xi_j))|v\rangle$ is the Laurent expansion of $f$
 at the point $q_j$ with respect to the formal coordinate $\xi_j$;
\item  For all $|v\rangle \in \cF_N$, there exists  a meromorphic one-form
$\omega \in H^0(C,\omega_C(*\sum_{j=1}^Nq_j))$ such that
$\langle \Phi| \rho_j(\ovpsi(\xi_j))|v\rangle d\xi_j$ is the Laurent expansion
of $\omega$
at the point $q_j$ with respect to the coordinates $\xi_j$,
\end{enumerate}
where $\rho_j(A)$ means that the operator $A$ acts on the $j$-th component of
$\cF_N$  as
$$
\rho_j(A)|u_1 \otimes u_2 \otimes \cdots \otimes u_N\rangle =
(-1)^{p_1+\cdots+p_{j-1}} |u_1 \otimes \cdots \otimes u_{j-1}\otimes
Au_j \otimes u_{j+1} \otimes \cdots \otimes u_N\rangle.
$$
\end{definition}
We will reformulate the above two conditions into gauge conditions. For that purpose we
introduce the following notation.

For a meromorphic one-form $\omega \in H^0(C,\omega_C(*\sum_{j=1}^Nq_jj))$ we let
$$
\omega_j = (\sum_{n=-n_0}^\infty a_n \xi_j^n )d\xi_j
$$
be the Laurent expansion at $q_j$ with respect to the coordinate $\xi_j$.
Then, for the field operator $\psi(z)$ let us define $\psi[\omega_j]$ by
$$
\psi[\omega_j]= \Res_{\xi_j=0}(\psi(\xi_j)\omega_j) =
\sum_{n= -n_0}^\infty a_n\psi_{n+1/2}.
$$
Similarly we can define $\ovpsi[\omega_j]$.
For  a meromorphic function $f \in H^0(C,\cO_C(*\sum_{j=1}^Nq_j))$ we let $f_j(\xi_j)$
be
the Laurent expansion of $f$ at $q_j$ with respect to the coordinate $\xi_j$.
For the field  operator $\psi(z)$ define $\psi[f_j]$ by
$$
\psi[f_j] = \Res_{\xi_j=0}(\psi(\xi_j)f_j(\xi_j)d\xi_j)
$$
Put
\begin{eqnarray*}
\psi[\omega] = (\psi[\omega_1], \ldots, \psi[\omega_N]) , &&
\ovpsi[\omega] = (\ovpsi[\omega_1], \ldots, \ovpsi[\omega_N]) \\
\psi[f]= (\psi[f_1], \ldots, \psi[f_N]), &&
\ovpsi[f]= (\ovpsi[f_1], \ldots, \ovpsi[f_N]).
\end{eqnarray*}
Then, these operate  on $\cF_N$ from the left
and on $\cFd_N$ from the right. For example,
$\ovpsi[f]$ operates on $\cF_N$ from the left by
\begin{eqnarray*}
\ovpsi[f] |u_1\otimes \cdots \otimes u_N \rangle  & = &
\sum_{j=1}^N\rho_j(\ovpsi[f_j])|u_1\otimes \cdots \otimes u_N \rangle \\
&= & \sum_{j=1}^N(-1)^{p_1+\cdots+p_{j-1}}|u_1 \otimes \cdots \otimes u_{j-1}
\otimes \ovpsi[f_j]u_j \otimes u_{j+1}\otimes \cdots \otimes u_N\rangle
\end{eqnarray*}
for $|u_j\rangle \in \cFd(p_j)$ and  operates on $\cFd_N$ from the right by
\begin{eqnarray*}
\langle v_N \otimes \cdots v_1|\ovpsi[f]   & = &
\sum_{j=1}^N \langle v_N \otimes \cdots v_1|\rho_j(\ovpsi[f_j]) \\
&= & \sum_{j=1}^N(-1)^{p_1+\cdots+p_{j-1}} \langle v_N \otimes
\cdots \otimes v_{j+1}
\otimes v_j \ovpsi[f_j] \otimes v_{j-1}\otimes \cdots \otimes v_1 |
\end{eqnarray*}
for $\langle v_j| \in \cFd(p_j)$.
\begin{theorem}[{[2, Theorem 3.1]}]
\label{thm3.1}
The element $\langle \Phi| \in \cFd_N$ belongs to
the space of abelian vacua  $\Vdagab(\gX)$
of the $j=0$ ghost system if and only if $\langle \Phi|$ satisfies
the following two conditions.
\begin{enumerate}
\item $\langle\Phi|\psi[\omega] =0$ for any meromorphic one-form $\omega \in
H^0(C, \omega_C(* \sum_{j=1}^Nq_j))$.
\item  $\langle\Phi|\ovpsi[f] =0$
for any meromorphic function $f \in H^0(C,\cO_C(*\sum_{j=1}^Nq_j))$.
\end{enumerate}
\end{theorem}

The first (resp. second) condition in the above theorem is called the
first (resp. second) gauge condition. The first and second gauge conditions can be
rewritten in the following form:
\begin{enumerate}
\item  $\sum_{j=1}^N(-1)^{p_1+\cdots+p_{j-1}}
\langle \Phi|   u_1 \otimes \cdots
\otimes \cdots \otimes u_{j-1}
\otimes \psi[\omega_j]u_j \otimes u_{j+1}\otimes \cdots \otimes u_N\rangle =0$
for any $\omega \in H^0(C, \omega_C(*\sum_{j=1}^Nq_j))$ and $|u_j\rangle \in \cF(p_j)$,
$j=1,2,\ldots,N$.
\item  $\sum_{j=1}^N(-1)^{p_1+\cdots+p_{j-1}}
\langle \Phi|  u_1 \otimes \cdots \otimes u_{j-1}
\otimes \ovpsi[f_j]u_j \otimes
u_{j+1}\otimes \cdots \otimes u_N\rangle =0$
for any $f \in H^0(C, \cO_C(*\sum_{j=1}^Nq_j))$ and $|u_j\rangle \in \cF(p_j)$,
$j=1,2,\ldots,N$.
\end{enumerate}

It is easy to show that the abelian vacua  $\Vdagab(\gX)$ is a
finite dimensional vector space. More strongly we can prove the
following theorem.
\begin{theorem}[{[2, Theorem 3.2]}]
\label{thm3.2}
 For any pointed Riemann surface $\gX= (C;Q_1,\ldots,Q_N;\xi_1,\ldots
\xi_N)$ with formal neighbourhoods we have
$$
\dim_\bC \Vdagab(\gX) =1.
$$
\end{theorem}

In the later application we need to consider a disconnected Riemann surface.
The following proposition is an immediate consequence of the definition.
\begin{proposition}[{[2, Proposition 3.1]}]
\label{prop3.1}
Let
$$\gX_1=(C_1;q_1, \ldots,q_M;\xi_1, \ldots, \xi_M)$$
and
$$\gX_2=(C_2;q_{M+1}, \ldots,q_M;\xi_{M+1}, \ldots, \xi_N)$$
be \fm pointed  Riemann surfaces  with formal neighbourhoods. Let $C$ be the disjoint union
$C_1 \sqcup C_2$ of the Riemann surfaces $C_1$, $C_2$.  Put
$$
\gX=(C; q_1, \ldots, q_N; \xi_1, \ldots,\xi_N).
$$
Then we have
$$
\Vdagab(\gX) = \Vdagab(\gX_1) \otimes
\Vdagab(\gX_2).
$$
\end{proposition}

Now we can introduce the dual  abelian vacua.
\begin{definition}
\label{dfn3.2}
Let $\cF_{\ab}(\gX)$ be the subspace of $\cF_N$ spanned by
$\psi[\omega] \cF_N$,
$\omega \in H^0(C,\omega_C(*\sum_{j=1}^Nq_j))$
and $\ovpsi[f] \cF_N$,  $f \in H^0(C,\cO(*\sum_{j=1}^Nq_j))$.
Put
$$
\cV_{\ab}(\gX) = \cF_N/\cF_{\ab}(\gX).
$$
The quotients space $\cV_{\ab}(\gX)$ is called
 the space of {\it dual abelian vacua} or {\it dual ghost vacua} of the $j=0$ ghost system.
\end{definition}

Since
$\Vdagab(\gX)$ is finite dimensional,  $\cV_{\ab}(\gX)$ is
dual to $\Vdagab(\gX)$.

\begin{theorem}[{[2, Theorem 3.2.]}]
\label{thm3.3}
The space of abelian vacua
$\Vdagab(\gX)$ is isomorphic to the determinant of the canonical bundle
$\omega_C$
$$\Vdagab(\gX) \cong \det(H^0(C,\omega_C)). $$
\end{theorem}

Let $\gX=(C;q_1,\ldots,q_N;\xi_1,\ldots,\xi_N)$  be an $N$-pointed Riemann surface with
formal neighbourhoods. Let $q_{N+1}$ be a non-singular point and choose a
formal coordinate $\xi_{N+1}$ of $C$ with center $q_{N+1}$.  Put
$$
\widetilde{\gX} = (C;q_1,\ldots,q_N,q_{N+1} ;\xi_1,\ldots,\xi_N, \xi_{N+1}).
$$
Then the canonical linear mapping
\begin{eqnarray*}
\iota : \cF_N & \rightarrow & \cF_{N+1} \\
|v\rangle & \mapsto & |v\rangle \otimes |0\rangle
\end{eqnarray*}
induces  the canonical mapping
$$
\iota^* :\cFd_{N+1}  \rightarrow \cFd_N.
$$
\begin{theorem}[{[2, Theorem 3.4]}]
\label{thm3.4}
The canonical mapping $\iota^*$ induces
an isomorphism
$$
\Vdagab(\widetilde{\gX}) \cong \Vdagab(\gX) .
$$
\end{theorem}

This isomorphism is denoted the ``Propagation of abelian vacua''
isomorphism.

Let us consider change of formal neighbourhoods. We use the same notation as in \S4.5.
For any automorphism $h \in {\mathcal D}_+^0$ we can define the action $G[h]$ on the fermion
Fock space $\mathcal F$ by using the energy-momentum tensor \eqref{emab} as in \S4.5.

\begin{proposition}[{[2, Proposition 6.1]}]
\label{prop6.1a}
For any $h_j \in {\mathcal D}_+^0$, $j=1,2, \ldots, N$ and $N$-pointed curve
$$\mathfrak{X}= (C; Q_1,Q_2, \ldots, Q_n; \xi_1, \xi_2, \ldots, \xi_N)$$
with formal neighbourhoods, put
$$
\mathfrak{X}_{(h)} =(C; Q_1,Q_2, \ldots, Q_N;
h_1(\xi_1), h_2(\xi_2), \ldots, h_N(\xi_N)).
$$
Then, the isomorphism $G[h_1] \widehat{\otimes}\cdots
\widehat{\otimes}G[h_N]$
\begin{eqnarray*}
\cFd_N  & \rightarrow & \cFd_N \\
\langle \phi_1\widehat{\otimes} \cdots \widehat{\otimes}\phi_N|
& \mapsto & \langle \phi_1G[h_1] \widehat{\otimes}\cdots
\widehat{\otimes}\phi_NG[h_N]|
\end{eqnarray*}
induces the canonical isomorphism
$$
G[h_1] \widehat{\otimes}\cdots \widehat{\otimes}G[h_N]:
\Vdagab(\gX) \rightarrow
\Vdagab(\mathfrak{X}_{(h)})
$$
\end{proposition}

\subsection{The space of \abelian vacua associated to a Riemann surface}
Let $\e C$ be a compact Riemann surface. For a \stable
and \fm pointed Riemann surface with formal neighbourhoods $\goX$ we denote
by $\tc(\goX)$ the underlying Riemann surface.
Suppose we now have two pointed Riemann surfaces with formal
neighbourhoods $\goX_i$ such that $\tc(\goX_1) = \tc(\goX_2)=\e C$.
Choose for each component of $\e C$ a point with a formal neighbourhood, which is not
a point with formal neighbourhoods of $\goX_i$, i=1,2. Let $\goX_0$
be the resulting \stable
and \fm pointed Riemann surface with formal neighbourhoods (if $\goX_0$ is not stable, then add further
points with formal neighbourhoods). Then iterations of the
propagation of vacua isomorphism determined by the inclusion of $\mathcal
F_N$ into $\mathcal F_{N+1}$ given by $|v\rangle \mapsto |v\rangle\otimes
|0\rangle$, induces isomorphisms from $\Vdagab(\goX_0)$ to
$\Vdagab(\goX_i)$, $i=1,2$. It is elementary to check that the
resulting isomorphism from $\Vdagab(\goX_1)$ to $\Vdagab(\goX_2)$
is independent of $\goX_0$.

Furthermore, we get from the commutativity of the following
diagram %(see \S4.5  for the definition of $G[h]$)
\begin{equation*}
\begin{CD}
{\mathcal F}_1 @> = >> {\mathcal F}_1\\
@V \id \otimes |0\rangle VV               @V \id \otimes |0\rangle VV\\
{\mathcal F}_2 @> \id\otimes G[h] >> {\mathcal F}_2,
\end{CD}
\end{equation*}
which  follows from the fact that $G[h] |0\rangle =
|0\rangle$,
that these isomorphisms are also compatible with the change of
formal coordinates isomorphism induced by $G[h]$.

\begin{definition}\label{Defspofvab}
The space of \abelian vacua associated to the Riemann surface $\e C$
is by definition
\[\Vdagab(\e C) = \coprod_{\tc(\goX) = \e C}\Vdagab(\goX)/\sim,\]
where the disjoint union is over all Riemann surfaces with formal neighbourhoods
with $\e C$ as the underlying Riemann surface and $\sim$ is the equivalence relation
generated by the isomorphisms discussed above.
\end{definition}

It is obvious that
\begin{proposition}\label{Spofvisoab}
The natural quotient map from $\Vdagab(\goX)$ to  $\Vdagab(\e C)$
is an isomorphism
 for all Riemann surfaces with formal neighbourhoods $\goX$
with $\tc(\goX) = \e C$.
\end{proposition}

Suppose $\e C_i$, $i=1,2$ are Riemann surfaces and $\Phi : \e C_1 \ra
\e C_2$ is a morphism of labeled marked Riemann surfaces. Let $\goX_2$ be a
Riemann surface with formal neighbourhoods such that $\tc(\goX_2) = \e
C_2$. Let $\Phi^*\goX_2 = \goX_1$. Then $\Phi$ is a morphism
of Riemann surfaces with formal neighbourhoods. We clearly have that

\begin{proposition}\label{morphspofvab}
The identity map on $\mathcal F_N$ induces a linear isomorphism from
$\Vdagab(\goX_1)$ to $\Vdagab(\goX_2)$, which induces a well defined linear isomorphism
$\Vdagab(\Phi)$ from $\Vdagab(\e C_1)$ to $\Vdagab(\e
C_2)$. Compositions of morphisms of labeled marked Riemann surfaces
go to compositions of the induced linear isomorphisms.
\end{proposition}

\section{Definition of the line bundle of \abelian vacua over Teichm\"{u}ller space}\label{New4ab}

\subsection{Sheaf of abelian vacua}

Let
$$
\gF = (\pi : \cC \rightarrow \cB; s_1, \ldots, s_N; \xi_1, \ldots, \xi_N)
$$
be a family of $N$-pointed semi-stable curves with formal neighbourhoods. That is
$\cC$ and $\cB$ are  complex manifolds, $\pi$ is a proper holomorphic
mapping, and for each point
$b \in \cB$,  $\gF(b)= (C_b =\pi^{-1}(b); s_1(b),
\ldots, s_N(b); \xi_1, \ldots,\xi_N)$
is an $N$-pointed semi-stable curve with formal neighbourhoods.
We let $\Sigma$
be the locus of double points of the fibers of $\gF$ and let $D$ be
$\pi(\Sigma)$. Note that $\Sigma$ is a non-singular submanifold of
codimension two in $\cC$,  and $D$ is a divisor
in $\cB$ whose irreducible components
$D_i$, $i = 1, 2, \dots ,m'$ are non-singular.

In this section  we use the following notation freely.
$$
S_j = s_j(\cB), \quad   S = \sum_{j=1}^N S_j.
$$
Put
$$
\cF_N(\cB)  = \cF_N \otimes_\bC \cO_\cB, \quad
\cFd_N(\cB)= \cO_\cB  \otimes_\bC    \cFd_N.
$$
\begin{definition}
\label{dfn4.1}
%{\rm
We define the subsheaf $\Vdagab(\gF)$ of $\cFd_N(\cB)$ by the gauge conditions:
\begin{eqnarray*}
&&\sum_{j=1}^N \langle \Phi | \psi[\omega_j] =  0 ,
\quad \hbox{\rm for all $\omega  \in \pi_*(\omega_{\cC/\cB}(*S))$}, \\
&&\sum_{j=1}^N \langle \Phi | \ovpsi[f_j] =  0,
\quad \hbox{\rm for all $f  \in \pi_*\cO_\cC(*S)$} .
\end{eqnarray*}
where $\omega_j$ and  $f_j$ are the Laurent expansion of
$\omega$ and $f$  along $S_j$ with respect to
the coordinate $\xi_j$.

The sheaf $\Vdagab(\gF)$ is called the sheaf of  ($j=0$) abelian vacua
or the sheaf of abelian vacua of the family $\gF$. Similarly the
sheaf $\cV_{\ab}(\gF)$ of  ($j=0$)  dual
abelian vacua of the family is defined by
$$
\cV_{\ab}(\gF)= \cF_N(\cB)/\cF_{\ab}(\gF).
$$
where  $\cF_{\ab}(\gF)$ is the $\cO_\cB$-submodule of $\cF_N(\cB)$ given by
$\cF_{\ab}(\gF) = \cF^0_{\ab}(\gF) + \cF^1_{\ab}(\gF)$,
where $\cF^0_{\ab}(\gF)$ is the span of $\ovpsi[f]\cF_N(\cB)$ for
all $f \in \pi_*\cO_\cC(*S)$ and $\cF^1_{\ab}(\gF)$ is the span of
$\psi[\omega]\cF_N(\cB)$ for all
$\omega \in \pi_*\omega_{\cC/\cB}(*S)$.
%}
\end{definition}

Note that we have
$$
\Vdagab(\gF) = \underline{\Hom}_{\cO_\cB}(\cV_{\ab}(\gF),
\cO_\cB).
$$
Moreover, by the right exactness of the tensor product we have
that
\begin{equation}
\label{4.1}
 \cV_{\ab}(\gF)\otimes_{\cO_\cB} \cO_{\cB, b}/\mathfrak{m}_b
 \cong \cV_{\ab}(\gF(b)).
\end{equation}

\begin{theorem}[{[2, Theorem 5.2]}]
\label{thm5.2a}
The sheaves $\cV_{\ab}(\gF)$ and $\Vdagab(\gF)$ are invertible
$\cO_\cB$-modules.  They are  dual to each other.
\end{theorem}

\subsection{The line bundle of \abelian vacua over Teichm\"{u}ller space}

Let $\goF_i$, $i=1,2$  be two families of \stable and \fm pointed Riemann surfaces with formal
neighbourhoods. Assume we have a morphism of
families (not necessarily preserving sections nor formal
coordinates)
\begin{equation*}
\begin{CD}
{\mathcal C}_1 @> \Phi >> {\mathcal C}_2\\
@V  VV               @VVV\\
{\mathcal B}_1 @> \Psi >> {\mathcal B}_2,
\end{CD}
\end{equation*}
which is a fiberwise biholomorphism.

Let now $\goF_0 = (\mathcal C_1 \ra \mathcal B_1;\vec s_0,\vec
\eta_0)$ be obtained from $\goF_1$, by replacing $(\vec s_1,\vec
\eta_1)$ by $(\vec s_0,\vec
\eta_0)$ such that $\vec s_0(\mathcal B_1)$ is disjoint from $\vec s_1(\mathcal
B_1)$ and from $\vec {\tilde s}_2(\mathcal
B_1)$, where $\Phi\vec  {\tilde s}_2 = \vec s_2 \Psi$. The propagation of vacua
isomorphism induces an isomorphism between $\Vdagab(\goF_0)$ and
$\Vdagab(\goF_1)$. Furthermore the propagation of vacua induces an
isomorphism between $\Vdagab(\goF_0)$ and $\Vdagab(\tilde
\goF_2)$, where $\tilde \goF_2 = (\mathcal C_1 \ra \mathcal B_1; \vec {\tilde s}_2,
\vec {\tilde \eta}_2)$ and $\Phi \vec {\tilde \eta}_2= \Psi^*\vec
\eta_2$. The identity on $\mathcal F_N(\mathcal B_1)$ induces an isomorphism
between
$\Vdagab(\tilde{\mathcal F}_2)$ and $\Vdagab(\Psi^*(\mathcal F_2))$.
Composing these with the pull back isomorphism just as in the
non-abelian case, we arrive at the following proposition.

\begin{proposition}\label{Tfamisoab}
We get an induced bundle morphism
\begin{equation}
\begin{CD}
\Vdagab(\goF_1) @> \Vdagab(\Phi) >> \Vdagab(\goF_2)\\
@V  VV               @VVV\\
{\mathcal B}_1 @> \Psi >> {\mathcal B}_2,
\end{CD}\label{}
\end{equation}
determined as above.
Moreover, composition of such family morphisms goes to composition
of the induced bundle morphisms.
\end{proposition}

Suppose now that we have two families $\goF_i$, $i=1,2$ over
$\Si$
with the property that they have the same image $\tPsi_{\goF_1}({\mathcal B}_1) =
\tPsi_{\goF_2}({\mathcal B}_2)$ in Teichm\"{u}ller space
$\cT_{\Si}$ and that $\goF_2$ is a \good family with respect to $\Si$.
For such a pair of families there exists by Proposition \ref{famequivalence} a unique
fiber preserving
biholomorphism $\Phi_{12} : \mathcal C_1\ra \mathcal C_2$ covering
$\Psi^{-1}_{\goF_2}\Psi_{\goF_1}$ such that $\Phi^{-1}_{\goF_2}
\Phi_{12} \Phi_{\goF_1} : Y \ra Y$ is isotopic to
$\Psi^{-1}_{\goF_2}\Psi_{\goF_1}\times \id$ through such fiber
preserving maps.

By Theorem \ref{thm5.2a} we have that
$\Vdagab(\goF_i)$ are holomorphic line bundles over $\mathcal
B_i$. By Proposition \ref{Tfamisoab} we get induced a glueing
isomorphism
\begin{equation}
\Vdagab(\Phi_{12}) : \Vdagab(\goF_1) \ra \Vdagab(\goF_2). \label{overlapisoab}
\end{equation}

\begin{definition}\label{dvbovacab}
Let $\Si$ be a closed oriented smooth surface.
We now define a line bundle $\Vdagab = \Vdagab(\Si,P)$ over Teichm\"{u}ller space
$\cT_{\Si}$ using the cover $\{\tPsi_{\goF}(\mathcal B)\}$, where $\goF$ runs
over the \stable and \fm  \good families of pointed Riemann surfaces with formal neighbourhoods over $\Si$.
Over $\tPsi_{\goF}(\mathcal
B)$ we specify the line bundle as $(\tPsi_\goF^{-1})^*\Vdagab(\goF)$. On
overlaps of the image of two \good families, we use
the glueing isomorphism $\Vdagab(\Phi_{12})$  to
glue the corresponding bundles together.
\end{definition}

We obviously have the following
\begin{proposition}\label{Teichpullbackab}
For any \stable and \fm family $\tgoF$ of pointed Riemann surfaces with
formal neighbourhoods over $\Si$ we have a preferred isomorphism
\[{\tilde \Upsilon}_{\tgoF} : \Vdagab(\tgoF) \ra \tPsi_{\tgoF}^*\Vdagab(\Si)\]
induced by the transformation isomorphism between $\Vdagab(\tgoF)$ and
$\Vdagab(\goF)$, for \good families $\goF$ of pointed Riemann surfaces with
formal neighbourhoods
over $\Si$ such that $\tPsi_\goF(\mathcal B)$ intersect $\tPsi_{\tgoF}({\mathcal B}')$
nonempty.
\end{proposition}

Suppose now $f : \Si_1 \ra \Si_2$ is a morphism of
surfaces. Then of course $f$ induces a
biholomorphism $f^*$ from $\cT_{\Si_1}$ to $\cT_{\Si_2}$.
Let now $\goF_1$ be a
\good family of stable pointed Riemann surfaces with formal
neighbourhoods over $\Si_1$. Then by composing with $f^{-1}\times
\id$ we get a \good family $\goF_2$ of stable pointed Riemann surfaces with formal
neighbourhoods over $\Si_2$ over the same base $\mathcal
B_1$. The identity morphism on $\mathcal F_N(\mathcal B_1)$ then induces a
morphism $\Vdagab(f) : \Vdagab(\goF_1) \ra
\Vdagab(\goF_2)$ which covers the identity on the base. This is
precisely the morphism induced from the morphism of families $\Phi_f =  f\times \id
: \goF_1 \ra
\goF_2$ by Proposition \ref{Tfamisoab}. This in turn induces a
morphism $\Vdagab(f) : (\tPsi_{\goF_1}^{-1})^*(\Vdagab(\goF_1)) \ra
(\tPsi_{\goF_2}^{-1})^*(\Vdagab(\goF_2))$ which covers
$f^* : \tPsi_{\goF_1}(B_1) \ra \tPsi_{\goF_2}(B_1)$.

\begin{proposition}\label{comptransfab}
The above construction provides a well defined lift of $f^* : \cT_{\Si_2} \ra
\cT_{\Si_1}$ to a morphism $\Vdagab(f) : \Vdagab(\Si_1)
\ra \Vdagab(\Si_2)$ which behaves well under compositions.
\end{proposition}

The proof is exactly the same as the proof of Proposition \ref{comptransf}.

\section{The connection in the line bundle of \abelian vacua over Teichm\"{u}ller space.}\label{New5ab}

Let us use the same notation as in \S6.1.
Let $\goF$ be a family of stable and saturated pointed Riemann surfaces with formal neighbourhoods
 on $(\Si,P)$.
For an element ${\vec \ell} = ({\ell}_1, \dots, {\ell}_N) $ in $$
\mathcal{L}(\mathfrak{F}) := \bigoplus_{j=1}^N \mathcal{O}
_\mathcal{B} [\; \xi_j^{-1}] \frac d{d \xi_j}, $$ the action
$D({\vec \ell})$ on $\goF$ is defined by
\begin{equation}
\label{4.2.8}
     D(\vec \ell) (F \otimes | u  \rangle)
    = \theta(\vec \ell)(F) \otimes | u \rangle -
F \cdot \Bigl(\sum_{j=1}^N \rho_j(T[{\ell}_j] )\Bigr) | u \rangle,
\end{equation}
where
$$
    F \in {\cO}_{\cB}, \quad | u \rangle \in \cF_N,
$$
and
$$
T[\ell] = \Res_{z=0}(T(z) \ell(z) dz),.
$$
Here $T(z)$ is the energy-momentum tensor T of spin $j=0$ bc ghost system.
Then the action has the similar properties as those of Proposition 6.1.
We define the dual action of $\cL(\gF)$ on $\cFd_N(\cB)$
by
\begin{equation}
\label{4.2.9}
  D(\vec \ell) (F \otimes \langle \Phi|) =
   \theta(\vec \ell)(F) \otimes \langle \Phi|  +
   \sum_{j=1}^N F  \cdot \langle \Phi| \rho_j(T[{\ell}_j]).
\end{equation}
where
$$
 F \in \cO_{\cB}, \quad \langle \Phi | \in \cFd_N(\cB).
$$
  Then, for any $| u \rangle \in \cF_N(\cB)$ and
$\langle\Phi | \in \cFd_N(\cB)$, we have
\begin{equation}
\label{4.2.10}
 \{D(\vec \ell) \langle \Phi| \}|\widetilde{\Phi}\rangle +
 \langle \Phi|\{D(\vec \ell)|\widetilde{\Phi}\rangle\}
=  \theta(\vec \ell)\langle \Phi|\widetilde{\Phi}\rangle .
\end{equation}

Now  the  operator $D(\vec \ell)$ acts on $\cV_{\ab}(\gF)$.
\begin{proposition}[{[2, Proposition 4.2]}]
\label{prop4.3a}
For any $\vec \ell \in \cL(\gF)$
we have
$$
D(\vec \ell)(\cF_{\ab}(\gF))  \subset
   \cF_{\ab}(\gF).
$$
Hence, $D(\vec \ell)$ operates on $\cV_{\ab}(\gF)$. Moreover, it is a
first order differential operator, if $\theta(\vec \ell) \neq 0$.
\end{proposition}

Now choose a meromorphic bidifferential
$$
\omega \in H^0(\cC\times_{\cB}\cC,
\omega_{\cC\times_{\cB}\cC}(2\Delta))
$$
defined by \eqref{primeform}.
Put
\begin{equation}
\label{aomega}
b_\omega(\vec{\ell}) =\sum_{j=1}^N \Res_{\xi_j=0}
\Big(\ell_j(\xi_j) S_\omega(\xi_j)d\xi_j\Bigr)
\end{equation}
where $S_\omega$ is the projective connection defined by \eqref{projectiveconnection}.
Then this defines an $\cO_\cB$-module homomorphism
$$
b_\omega : \cL(\gF) \rightarrow \cO_\cB,
$$
and if $\theta(\vec{\ell})=0$ then we have that
$$
\langle \Phi |\{D(\vec{\ell})| u \rangle \} =\frac16 b_\omega(\vec{\ell})
\langle \Phi | u \rangle .
$$
For a vector field $X$ on $\cB$
choose $\vec{\ell} \in \cL(\gF)$ such that $\theta(\vec{\ell})=X$.
Then the connection on $\Vdagab(\gF)$ is defined by
\begin{equation}
\label{4.19a}
\nabla_X^{(\omega)}(\langle\Phi|) =  D(\vec{\ell})(\langle\Phi|)
+ \frac16 b_\omega(\vec{\ell})    (\langle\Phi|),
\end{equation}
for $\langle\Phi| \in \Vdagab(\gF)$.
This is well-defined
and \eqref{4.19a} is independent of the choice of $\vec{\ell} \in\cL(\gF)$
with $\theta(\vec{\ell})=X$.
Just like for the non-abelian conformal field theory (see for example \cite{Ue2}, section 5)
we can prove the following theorem.
\begin{theorem}[{[2, Theorem 4.2]}]
\label{thm4.2a}
The operator $\nabla^{(\omega)}$ defines a
projectively flat holomorphic connection of
the sheaves $\cV_{\ab}(\gF)$ and $\Vdagab(\gF)$.
Moreover, the connection has a regular singularity along the locus
$D \subset \cB$ which is the locus of the singular curves.
The connection $\nabla^{(\omega)}$ depends on the choice of
bidifferential $\omega$ and if we choose another bidifferential
$\omega'$ then there exists a holomorphic one-form $\phi_{\omega,\omega'}$
on $\cB$ such that
\begin{equation}
\label{4.20a}
\nabla_X^{(\omega)} - \nabla_X^{(\omega')}= \frac16\langle
\phi_{\omega,\omega'}, \, X\rangle.
\end{equation}
Moreover, the curvature form $R$ of $\nabla_X^{(\omega)}$ is given by
\begin{equation}
\label{4.21a}
R(X,Y) = \frac16 \Big\{b_\omega(\vec{n}) - X(b_\omega(\vec{m}) )
+Y(b_\omega(\vec{\ell}) )- \sum_{j=1}^N
\Res_{\xi_j=0}\big(\frac{d^3\ell_j}{d\xi_j}m_jd \xi_j\bigr)\Bigr\},
\end{equation}
where $X$,$Y \in \Theta_{\cC/\cB}(*S))$,
$\vec{\ell}$, $\vec{m} \in \cL(\gF)$ with $X=\theta(\vec{\ell})$,
$Y= \theta(\vec{m})$, and $\vec{n} \in \cL(\gF)$ is defined by
$\vec{n} = [\vec{\ell}, \vec{m}]_d$  $($see \eqref{35f}$)$.
\end{theorem}

\begin{corollary} \label{comparison}
Let $\mathfrak{F}$ be a family of stable and saturated pointed Riemann surfaces
with formal neighbourhoods on $(\Si, P)$.
If we use the same bidifferential $\omega$ to define the connections on the bundle of vacua
and the line bundle of abelian vacua on $\mathfrak{F}$, then we have
\begin{equation}\label{relation}
R^\omega(X,Y) = \frac{c_v}{2} R(X,Y)\otimes \id .
\end{equation}
\end{corollary}

\begin{proposition}\label{confamab}
Let $\goF$ be a family of stable and saturated pointed Riemann surfaces with formal neighbourhoods
 on $(\Si,P)$ and choose a symplectic basis
$(\alpha_1, \ldots, \alpha_g,\beta_1,\ldots,\beta_g)$ of $H_1(\Si,
\Z{})$. Let $\omega\in H^0(\mathcal C \times_{\mathcal B} \mathcal
C, \omega_{\mathcal C \times_{\mathcal B} \mathcal C/\mathcal
B}(2\Delta))$ be the normalized symmetric bidifferential
determined by this data. Then there is the connection
$\tnabla^{(\omega)}$ in the bundle $\Vdagab(\goF)$ whose $(1,0)$-part is given by
formula (4.23) in \cite{AU1} and whose $(0,1)$-part is just the $\overline{\partial}$-operator
in this holomorphic line-bundle. The curvature of this connection is
given by the formula \eqref{4.21a}.
\end{proposition}

\proof {By the definition of $\tnabla^{(\omega)}$ and the
definition \eqref{aomega} of $b_\omega$, we see
that the $(1,1)$ and $(0,2)$-part of the curvature is zero.}
\eproof

Suppose now that we have two \good families $\goF_i$, $i=1,2$ with
the property that they have the same image $\tPsi_{\goF_1}
({\mathcal B}_1) = \tPsi_{\goF_2}({\mathcal B}_2)$ in Teichm\"{u}ller
space $\cT_{\Si}$. For such a pair of families there exists by
Proposition \ref{famequivalence} a unique fiber preserving
biholomorphism $\Phi_{12} : \mathcal C_1\ra \mathcal C_2$ covering
$\Psi^{-1}_{\goF_2}\Psi_{\goF_1}$ such that $\Phi^{-1}_{\goF_2}
\Phi_{12} \Phi_{\goF_1} : (Y, P) \ra (Y,P)$ is isotopic to
$\Psi^{-1}_{\goF_2}\Psi_{\goF_1}\times \id$.

\begin{lemma}\label{contransfab}
Let $\tnabla^{(\omega)}_i$ be the connection in $\Vdagab(\goF_i)$
described in Proposition \ref{confamab}. Then we have that
\[\Vdagab(\Phi_{12})^*(\tnabla^{(\omega)}_2) = \tnabla^{(\omega)}_1.\]
\end{lemma}

This follows from Theorem \ref{thm6.1} above, by the same argument
as in the non-abelian case.

\begin{theorem}\label{conTeichab}
Let $\Si$ be a closed oriented
surface and let $(\vec \alpha, \vec \beta) = (\alpha_1, \ldots,
\alpha_g,\beta_1,\ldots,\beta_g)$ be a symplectic basis of
$H_1(\Si, \Z{})$. There is a unique connection $\tnabla^{(\vec
\alpha, \vec \beta)} = \tnabla^{(\vec
\alpha, \vec \beta)}(\Si,P)$ in the bundle $\Vdagab(\Si,P)$ over
$\cT_{\Si}$ with the property that for any
\good family
$\goF$ of stable pointed Riemann surfaces with formal neighbourhoods over $\Si$ we have that
\[\tPsi_{\goF}^*(\tnabla^{(\vec \alpha, \vec \beta)}) = \tnabla^{(\omega)}.\]
In particular the connection is compatible with the holomorphic line-bundle structure
on $\Vdagab(\Si,P)$. The curvature is of type $(2,0)$ as stated in Proposition
\ref{confamab}.

If we act on the symplectic basis $(\vec \alpha,\vec \beta)$ by an element $\Lambda =
\left(\begin{array}{cc}A & B\\C & D\end{array}\right) \in
\text{Sp} (g,\Z{})$ so as to obtain $\Lambda(\vec \alpha, \vec \beta)$,
as defined by \eqref{action1},
then
\begin{equation}
\tnabla^{\Lambda(\vec \alpha, \vec \beta)} - \tnabla^{(\vec \alpha, \vec \beta)} =
\frac{1}{2} \Pi^*( d \log \det (C \tau + D )),\label{contransformab}
\end{equation}
where $\Pi$
is the period mapping of holomorphic one-forms form the base space of  $\mathfrak{F}$ to the
Siegel upper-half plane of degree $g$. If $f : \Si_1 \ra
\Si_2$ is an orientation preserving diffeomorphism of
surfaces which maps the symplectic basis $(\vec
\alpha^{(1)},\vec \beta^{(1)})$ of $H_1(\Si_1,\Z{})$ to the
symplectic basis $(\vec \alpha^{(2)},\vec \beta^{(2)})$ of
$H_1(\Si_2,\Z{})$ then we have that
\[\Vdagab(f)^*(\tnabla^{(\vec \alpha^{(2)},\vec \beta^{(2)})}) = \tnabla^{(\vec \alpha^{(1)},
\vec \beta^{(1)})}.\]
\end{theorem}

\proof {The existence is a consequence of Lemma \ref{contransfab}.
The transformation law \eqref{contransformab} follows from the above Theorem \ref{thm4.2a}.} \eproof

\begin{proposition}\label{Teichpullbackconab}
For any \stable and \fm family $\tgoF$ of pointed Riemann surfaces with
formal neighbourhoods over $\Si$ the preferred isomorphism
\[{\tilde \Upsilon}_{\tgoF} : \Vdagab(\tgoF) \ra \tPsi_{\tgoF}^*\Vdagab(\Si)\]
given by Proposition \ref{Teichpullbackab} preserves connections and is compatible
with the lift $\Vdagab(f)$.
\end{proposition}

This follows directly from Lemma \ref{contransfab} and Theorem \ref{conTeichab}.

\section{The preferred non-vanishing section of the bundle of \abelian vacua.}\label{prefsec}

Let $\gX=(C;Q;\xi)$ be a one-pointed smooth curve of genus $g$
with a formal neighbourhood. We shall show that if we fix a
symplectic basis $(\vec{\alpha}, \vec{\beta})=(\alpha_1, \ldots, \alpha_g, \beta_1, \ldots,
\beta_g)$ of $H_1(C, \bZ)$, then there is a canonical preferred non-zero
vector $\langle \omega(\gX,(\{\alpha,\beta\})| \in
\Vdagab(\gX)$. Let us choose a normalized basis
$\{\omega_1, \ldots, \omega_g\}$ of holomorphic one-forms on $C$
which is characterized by
\begin{equation}
\label{betaone2}
\int_{\beta_i}\omega_j = \delta_{i j}, \quad 1 \le i,j \le g.
\end{equation}
The period matrix is given by $$ \tau = (\tau_{ij}), \quad \tau_{ij} =
\int_{\alpha_i}\omega_j. $$ Now the numbers $I_n^i$, $n = 1,2, \ldots$, $i=1,\ldots
g$ are defined by
$$ \omega_i = (\sum_{n=1}^\infty I_n^i
\xi^{n-1})d\xi. $$
Note that the numbers $I_n^i$ depend on the
symplectic basis $(\vec{\alpha}, \vec{\beta})$ and the formal neighbourhood $\xi$.

For a positive integer $n \ge 1$ let $\omega_Q^{(n)}$ be a meromorphic
one-form on $C$ which has a pole of order $n+1$ at $Q$ and holomorphic elsewhere
such that
\begin{eqnarray}
\label{omegaQ1}
\int_{\alpha_i} \omega_Q^{(n)} &= &
-\frac{2 \pi \sqrt{-1}I_n^i}{n},
\quad \int_{\beta_i} \omega_Q^{(n)} = 0, \quad 1\le i \le g  \\
\label{omegaQ}
\omega_Q^{(n)} &=& \bigl( \frac{1}{\xi^{n+1}} +
\sum_{m=1}^\infty q_{n,m}\xi^{m-1}\bigr)d\xi.
\end{eqnarray}
These conditions uniquely determine $\omega_Q^{(n)}$. Note that the second equality
of \eqref{omegaQ1} and \eqref{omegaQ} imply the first equality of \eqref{omegaQ1}.
The preferred  element $\langle \omega(\gX,\{\alpha, \beta\})| \in \Vdagab(\gX)$
is defined by
$$
\langle \omega(\gX,\{\alpha, \beta\}) |  = \cdots e(\omega_{g+2} )\wedge e(\omega_{g+1})
\wedge e(\omega_g) \wedge  \cdots \wedge e(\omega_1),
$$
where
$$
\omega_{g+n} = \omega_Q^{(n)}.
$$
For details see Lemma 3.1 and its proof of \cite{AU1}.
We call $\{\omega_n\}$, $n=1,2,\ldots$ a normalized basis for $\gX$.
Note that the normalized basis depends on the choice of a symplectic
basis of $H_1(C,\bZ)$ and the coordinate $\xi$.

\begin{theorem}[{[2, Theorem 6.2]}]
\label{thm6.2}
For $h(\xi) \in \cD_+^0$ put $\gX_h = \{ C;Q;\eta= h(\xi)\}$.
Then
$$
\langle \omega(\gX,\{\alpha, \beta\})| G[h] = \langle \omega(\gX_h,\{\alpha, \beta\})|,
$$
where $G[h] : \Vdagab(\gX) \rightarrow \Vdagab(\gX_h)$ is  the canonical isomorphism
given in Proposition \ref{prop6.1a}
\end{theorem}
\begin{theorem}[{[2, Theorem 6.3]}]
\label{thm6.3}
Let $(\alpha_1, \ldots \alpha_g,\beta_1, \ldots, \beta_g )$
and $(\widetilde{\alpha}_1, \ldots \widetilde{\alpha}_g,
\widetilde{\beta}_1, \ldots, \widetilde{\beta}_g )$ be symplectic
bases of $H^1(C, \bZ)$ of the non-singular curve $C$. Assume that
$\{\beta_1, \ldots, \beta_g\}$ and
$\{\widetilde{\beta}_1, \ldots, \widetilde{\beta}_g \}$ span the same
Lagrangian sublattice in $H^1(C, \bZ)$. Then
$$
\langle \omega(\gX, \{\alpha,\beta\}) |  = \det U
\langle \omega(\gX, \{\widetilde{\alpha}, \widetilde{\beta}\}) | ,
$$
where $U \in GL(g,\bZ)$ is defined by
$$
\left( \begin{array}{c}
\widetilde{\beta}_1 \\ \vdots \\ \widetilde{\beta}_g
\end{array} \right)
= U
\left( \begin{array}{c}
\beta_1 \\ \vdots \\ \beta_g \end{array} \right) .
$$
\end{theorem}

Let $\{p, q\}$ be two smooth points on the curve $C$ with formal neighbourhoods
$\xi$, $\eta$, respectively. Put
$\gX_0 = (C; p,q; \xi , \eta)$, $\gX_1= (C; p; \xi)$, $\gX_2=(C;q;\eta)$.
Then the natural imbeddings
\begin{eqnarray*}
\iota_1&:&   \cF \hookrightarrow \cF_2 \\
  &&  |u\rangle \mapsto |u \rangle \otimes |0\rangle \\
\iota_2&:&   \cF \hookrightarrow \cF_2 \\
  &&  |u\rangle \mapsto |0\rangle  \otimes  |u \rangle\\
\end{eqnarray*}
induce canonical isomorphisms
$$
\begin{array}{ccccc}
&&\Vdagab(\gX_0)& &\\
&&  &&   \\
& {}^{\iota_1^*} \swarrow&& \searrow^{\iota_2^*} & \\
&&  &&   \\
&\Vdagab(\gX_1) && \Vdagab(\gX_2) &
\end{array}
$$ by Theorem \ref{thm3.4}.

\begin{theorem}[{[2, Theorem 6.4]}]
\label{thm6.4}
Under the above notation we have
$$
\iota_2^*\circ (\iota^*_1)^{-1}(\langle \omega(\gX_1,\{\alpha, \beta\})| )
= \langle\omega(\gX_2,\{\alpha, \beta\})|.
$$
\end{theorem}

Let $\Si$ be a closed oriented surface. Assume first that $\Si$ is
connected. As described above, the
choice of a symplectic basis gives a preferred section in the line
bundle of abelian vacua associated to any family of stable and
saturated pointed Riemann surfaces with formal neighbourhoods over $\Si$. We have
that

\begin{theorem}\label{prefsecfam}
Let $\goF$ be a family of stable and saturated pointed Riemann surfaces with formal neighbourhoods
 over $\Si$ and choose a symplectic basis
$(\vec \alpha, \vec \beta)$ of $H_1(\Si,
\Z{})$. Then there is a preferred non-vanishing holomorphic section
$s^{(\vec \alpha, \vec \beta)}_{\goF}$ in the bundle $\Vdagab(\goF)^{\otimes 2}$ given by
$$
s^{(\vec \alpha, \vec \beta)}_{\goF}(t) = (\langle \omega(\gX_t, \{\alpha(t), \beta(t)\}|)^{\otimes 2}.
$$
If we act on the symplectic basis $(\vec \alpha, \vec \beta)$ by
an element $\Lambda =
\left(\begin{array}{cc}(U^t)^{-1} & B\\0 & U\end{array}\right) \in
\text{Sp} (g,\Z{})$ in order to obtain $\Lambda(\vec \alpha, \vec \beta)$ as described
 in \eqref{action1}, then
\begin{equation}
s^{\Lambda(\vec \alpha, \vec \beta)} = s^{(\vec \alpha, \vec
\beta)}.\label{sectiontransformfam}
\end{equation}
\end{theorem}

\proof This is clear from Theorem \ref{thm6.3}, since we have $\det U = \pm 1$.
\eproof

Thus the section $s^{(\vec \alpha, \vec \beta)}_{\goF}$ only really
depends on the Lagrangian subspace $L = \spane \{\beta_i\}$ and we
therefore denote it $s_\goF (L)$.

Suppose now that $\Si$ is not connected and that $\Si =\coprod_i\Si_i$ is the decomposition
of $\Si$ into
its connected components $\Si_i$. Let $\goF$ be a family of stable pointed Riemann surfaces with formal
neighbourhoods over $\Si$. Let $\goF_i$
be the restriction of $\goF$ to $\Si_i$. Let $N_i$ be the number of sections of $\goF_i$
and $N = \sum_i N_i$ the number of sections
of $\goF$. We obviously have the following lemma.

\begin{lemma}\label{DJDlemma}
The isomorphism $\mathcal F_N \cong \otimes_i \mathcal F_{N_i}$ induces an isomorphism of
holomorphic line bundles
$$\Vdagab(\goF) \cong \otimes_i \Vdagab(\goF_i),$$
which is compatible with the connections.
\end{lemma}

Suppose now that $(\vec\alpha, \vec\beta) = ((\vec\alpha_i, \vec\beta_i))$ is a symplectic
basis of $H_1(\Si,\Z)$. We then define the preferred
section to be
$$s_\goF^{(\vec\alpha, \vec\beta)} = \otimes_i s_{\goF_i}^{(\vec\alpha_i, \vec\beta_i)}.$$

For the rest of this section $\Si$ is just any closed oriented
surface.

Suppose now that we have two \good families $\goF_i$, $i=1,2$ with
the property that they have the same image $\tPsi_{\goF_1}
({\mathcal B}_1) = \tPsi_{\goF_2}({\mathcal B}_2)$ in Teichm\"{u}ller
space $\cT_{\Si}$. For such a pair of families there exists by
Theorem \ref{famequivalence} a unique fiber preserving
biholomorphism $\Phi_{12} : \mathcal C_1\ra \mathcal C_2$ covering
$\Psi^{-1}_{\goF_2}\Psi_{\goF_1}$ such that $\Phi^{-1}_{\goF_2}
\Phi_{12} \Phi_{\goF_1} : Y \ra Y$ is isotopic to
$\Psi^{-1}_{\goF_2}\Psi_{\goF_1}\times \id$ through fiber preserving diffeomorphisms.

\begin{lemma}\label{prefsecfamch}
Let $s^{(\vec \alpha, \vec \beta)}_{\goF_i}$ be the preferred sections of $\Vdagab(\goF_i)^{\otimes 2}$
described in Theorem \ref{prefsecfam}. Then we have that
\[\Vdagab(\Phi_{12})^{\otimes 2}
(s^{(\vec \alpha, \vec \beta)}_{\goF_1}) = s^{(\vec \alpha, \vec \beta)}_{\goF_2}.\]
\end{lemma}

\proof {Since $\Vdagab(\Phi_{12})$ is induced by the propagation of vacua isomorphism and the coordinate
change isomorphism, this lemma follows from Theorem \ref{thm6.2} and \ref{thm6.4}.} \eproof

\begin{theorem}\label{secTeichab}
Let $\Si$ be a closed oriented
surface and let $(\vec \alpha, \vec \beta) = (\alpha_1, \ldots,
\alpha_g,\beta_1,\ldots,\beta_g)$ be a symplectic basis of
$H_1(\Si, \Z{})$. Then there is a unique non-vanishing holomorphic section
$s^{(\vec \alpha, \vec \beta)}= s^{(\vec \alpha, \vec \beta)}_{\Si}$
in the bundle $\Vdagab(\Si)^{\otimes 2}$ over
$\cT_{\Si}$ with the property that for any
\good family
$\goF$ of stable pointed Riemann surfaces with formal neighbourhoods over $\Si$ we have that
\[(\tPsi_{\goF}^*)^{\otimes 2}(s^{(\vec \alpha, \vec \beta)}) = s^{(\vec \alpha, \vec \beta)}_{\goF}.\]
The sections transforms according to the transformation rule (\ref{sectiontransformfam}.
 If $f : \Si_1 \ra \Si_2$ is an orientation
preserving diffeomorphism of surfaces which maps the symplectic
basis $(\vec \alpha^{(1)},\vec \beta^{(1)})$ of $H_1(\Si_1,\Z{})$
to the symplectic basis $(\vec \alpha^{(2)},\vec \beta^{(2)})$ of
$H_1(\Si_2,\Z{})$ then we have that
\[(\Vdagab(f)^*)^{\otimes 2}(s^{(\vec \alpha^{(2)},\vec \beta^{(2)})}_{\Si_2}) = s^{(\vec \alpha^{(1)},
\vec \beta^{(1)})}_{\Si_1}.\]
\end{theorem}

\proof {The first part of this theorem follows from the above.
Since the choice of basis of meromorphic 1-forms with the relevant properties
is natural w.r.t. morphisms of Riemann surfaces, we get that the preferred section
transforms just like the bases does.} \eproof

Likewise, we see that the section only depends on the Lagrangian
subspace and we denote it therefore $s(L) = s_\Si(L)$.

\begin{proposition}\label{Teichpullbackconabsec}
For any \stable and \fm family $\tgoF$ of pointed Riemann surfaces with
formal neighbourhoods over $\Si$ the preferred isomorphism
\[{\tilde \Upsilon}_{\tgoF}^{\otimes 2} : \Vdagab(\tgoF)^{\otimes 2} \ra
\tPsi_{\tgoF}^*\Vdagab(\Si)^{\otimes 2}\]
given by Proposition \ref{Teichpullbackab} preserves the preferred sections.
\end{proposition}

This Follows from Lemma \ref{prefsecfamch}.

\section{The geometric construction of the modular functor.}\label{construction}

For the convenience of the reader,
let us summarize the results of the sheaf of vacua constructions
over Teichm\"{u}ller spaces of pointed surfaces obtained in non-abelian case in sections \ref{New3} to \ref{New5}
 and in the abelian case in sections \ref{New3ab} to \ref{prefsec}.

\begin{theorem}\label{mainconstT}
Let $(\Si, P, \l)$ be a \stable and \fm labeled pointed surface.

\begin{itemize}
\item The sheaf of vacua construction (see Definition \ref{dvbovac}) yields
a vector bundle $\Spofvlam = \Spofvlam(\Si, P)$ over the
Teichm\"{u}ller space $\cT_{(\Si,P)}$ of $(\Si,P)$ whose fiber at a
complex structure $\e C$ on $(\Si, P)$ is identified (via the
isomorphism given in Proposition \ref{Teichpullback}) with the
space of vacua $\Spofvlam (\e C)$ as defined in Definition
\ref{Defspofv}.
\item For each symplectic basis  of
$H_1(\Si,\Z)$, we get induced a
connection in
$\Spofvlam$ over $\cT_{(\Si,P)}$. Any two of these connections
differ by a global scalar-value 1-form on $\cT_{(\Si,P)}$. See
Theorem \ref{conTeich}.
\item Each of these connections is projectively flat and their curvature
are described in details in Theorem \ref{curvaturecom} and \ref{thm4.2a}.

\item There
is a natural lift of morphisms of pointed surfaces to these
bundles covering induced biholomorphisms between Teichm\"{u}ller
spaces, which preserves compositions. See Proposition
\ref{comptransf}.
\item A morphism of pointed surfaces transforms these connections according to
the way it transforms symplectic bases of the first homology.
See Theorem \ref{contransf}.
\end{itemize}
\end{theorem}

\begin{remark}\label{remconbaslag}{\em If we choose a Lagrangian subspace
$L$ of $H_1(\Si,\Z)$ and
constrain the symplectic basis $(\alpha_i,\beta_i)$ of
$H_1(\Si,\Z)$ such that $L = \spane \{\beta_i\}$ then we see from
the transformation laws in Theorem \ref{contransf}, that we
get a connection in $\Spofvlam$ which depends only on $L$.}
\end{remark}

Since the connections in the vector bundle $\Spofvlam$ are only
projectively flat, we need a 1-dimensional theory with
connections, whose curvature
 after taking tensor products, can cancel this curvature and result in a bundle with a flat
connection. There are obstructions to doing this mapping class
group equivariantly, so we expect to see central extension of the
mapping class groups occurring. As we shall see below, this is
exactly what happens, when one extracts the necessary root of the
abelian theory treated in \cite{AU1}, so as to get the right
scaling of the curvature. Again, the following theorem summarizes
the results and constructions, now in the abelian case treated in
section \ref{New3ab} to \ref{prefsec}.

\begin{theorem}\label{mainconstTab}
Let $\Si$ be a closed oriented surface.

\begin{itemize}
\item The sheaf of abelian vacua construction (see Definition \ref{dvbovacab})
yields a line bundle $\Vdagab = \Vdagab(\Si)$ over the Teichm\"{u}ller
space $\cT_{\Si}$ of $\Si$, whose fiber at a complex structure $\e
C$ on $\Si$ is identified (via the isomorphism given in
Proposition
\ref{Teichpullbackab}) with the space of abelian vacua $\Vdagab
(\e C)$ as defined in Definition \ref{Defspofvab}.
\item For each symplectic basis  of
$H_1(\Si,\Z)$, we get induced a holomorphic connection in
$\Vdagab$ over $\cT_{\Si}$ (see Theorem \ref{conTeichab}). The
difference between the connections associated to two different
basis's is the global scalar-value 1-form on $\cT_{\Si}$ given in
(\ref{contransformab}).
\item The curvature of each of these connections are described in Proposition \ref{confamab}.
\item For each symplectic basis  of
$H_1(\Si,\Z)$, we also get a preferred non-vanishing section of
$(\Vdagab)^{\otimes 2}$ as specified in Theorem \ref{secTeichab}. The
transformation formula (\ref{sectiontransformfam}) states how the
preferred sections transforms under change of the symplectic basis
of $H_1(\Si,\Z)$.
\item There
is a natural lift of morphisms of surfaces to these bundles
covering induced biholomorphisms between Teichm\"{u}ller space, which
preserves compositions. See Proposition \ref{comptransfab}.
\item A morphism of surfaces transforms these connections and the preferred sections
according to the way it transforms symplectic bases of the first
homology. See Theorem \ref{contransfab} and \ref{secTeichab}.
\end{itemize}
\end{theorem}

\begin{remark}{\em If we choose a Lagrangian subspace $L$ of $H_1(\Si,\Z)$ and
constrain the symplectic basis $(\alpha_i,\beta_i)$ of
$H_1(\Si,\Z)$ such that $L = \spane \{\beta_i\}$ then we see from
the transformation laws in Theorem \ref{secTeichab}, that we
get a preferred non-vanishing section $s=s(L)$ and a connection in
$(\Vdagab)^{\otimes 2}$ which only depends on $L$. }
\end{remark}

From the discussion of the curvatures of the connections in
$\Spofvlam$ and $\Vdagab$, i.e. by comparing the curvature formula \eqref{relation},
it is clear that the root of $(\Vdagab)^{\otimes 2}$
we are seeking is $c_{\cv}$. The following theorem
provided us with such a root.

\begin{theorem}\label{fracpowerab}
For any marked surface $\Sib = (\Si, L)$ there exists a line
bundle, which we denoted $(\Vdagab)^{-\frac{1}{2}c_{\cv}}(L)=
(\Vdagab)^{-\frac{1}{2}c_{\cv}}(\Sib)$, over $\cT_\Si$ that
satisfies the following:
\begin{itemize}
\item $(\Vdagab)^{-\frac{1}{2}c_{\cv}}$ is a
functor from the category of marked surfaces to the category of
line bundles over Teichm\"{u}ller spaces of closed oriented surfaces.
\item If we choose a symplectic basis of $H_1(\Si,\Z)$ for a marked
surface $\Sib$ then we get induced a connection in
$(\Vdagab)^{-\frac{1}{2}c_{\cv}}(L)$, whose curvature is
$-\frac{1}{2}c_{\cv}$ times the curvature of the corresponding
connection in $\Vdagab$. The difference between the connections
associated to two different bases is $-\frac{1}{2}c_{\cv}$ times
the global scalar-value 1-form on $\cT_{\Si}$ given in
(\ref{contransformab}).
\end{itemize}
\end{theorem}

\proof Let $((\Vdagab)^{\otimes 2})^*$ be the complement of the zero section of
$(\Vdagab)^{\otimes 2}$. Let $\widetilde{({\Vdagab})^{\otimes 2}}$ be the fiberwise universal cover of
$((\Vdagab)^{\otimes 2})^*$ based at the section $s(L)$. This is a completely
functorial construction on pairs of line bundles and non-vanishing
sections. There is a unique lift of the $\C{}^*$-action on
$((\Vdagab)^{\otimes 2})^*$ to a $\C{}$-action on $\widetilde{({\Vdagab})^{\otimes 2}}$ with respect to the
covering map $\exp$ from $\C{}$ to $\C{}^*$. For any $\alpha\in \C{}^*$ we can now
functorially define a line bundle $((\Vdagab)^{\otimes 2})^\alpha(L)$ as
follows:
\[((\Vdagab)^{\otimes 2})^\alpha(L) = \widetilde{({\Vdagab})^{\otimes 2}} \times_{\rho_\alpha}\C{},\]
where $\rho_\alpha(z) : \C{} \ra \C{}$ is the linear map given by
multiplication by $\exp(\alpha z)$ for all $z\in \C{}$. We
emphasis the dependence of this bundle on the section $s(L)$ and
hence on $L$ in the notation for this bundle. Here we choose
$\alpha =- c_v/4$ to define $(\Vdagab)^{-\frac{1}{2}c_{\cv}}$.

It is clear from the construction of $((\Vdagab)^{\otimes 2})^\alpha(L)$, that a
connection in $\Vdagab$ will induce a connection in
$((\Vdagab)^{\otimes 2})^\alpha(L)$, whose curvature two-form is $\alpha/2$ times
the curvature two-form of that connection in $\Vdagab$. For the
construction of the functor on the morphisms of marked surfaces,
we refer to \cite{Walker} and \cite{A}. \eproof

By pulling $(\Vdagab)^{-\frac{1}{2}c_{\cv}}(L)$ with its
connection back to $\cT_{(\Si,P)}$ from $\cT_\Si$, we get a line
bundle with with a connection on $\cT_{(\Si,P)}$, which we also
denote $(\Vdagab)^{-\frac{1}{2}c_{\cv}}(L)$.

Let now $(\Sib, \l) = (\Si, P,V,L, \l)$ be a \stable and \fm
labeled marked surface. From the above Theorems \ref{mainconstT}
and \ref{fracpowerab}, we see that there is a well defined flat
connection in the vector bundle $\Spofvlam \otimes
(\Vdagab)^{-\frac{1}{2}c_{\cv}}(L)$ over $\cT_{(\Si, P)}$ gotten
by taking the tensor product connection of the two connections
induced by any symplectic basis $(\alpha_i,\beta_i)$ of
$H_1(\Si,\Z)$. Now $\cT_{(\Si, P)}$ forms a $\RPP$-principal
bundle over the reduced Teichm\"{u}ller space $\ctT_{(\Si, P)}$. Hence
we can use the flat connection to push forward this bundle to
obtain a bundle with a flat connection over the reduced
Teichm\"{u}ller space.

\begin{definition}\label{Flatbdloredteich}
For the \stable and \fm labeled marked surface $(\Sib,\l)$ we
define the vector bundle $\Spofvlam(\Sib)$ with its flat
connection $\nabla(\Sib,\l)$ as the push forward of the bundle
$\Spofvlam \otimes (\Vdagab)^{-\frac{1}{2}c_{\cv}}(L)$ to the
reduced Teichm\"{u}ller space $\ctT_{(\Si, P)}$ followed by
restriction to the fiber $\cT_{\Sib}$.
\end{definition}

\begin{remark}{\em By Theorem \ref{mainconstT} and \ref{mainconstTab} and Lemma
\ref{momswd} we see that morphism of \stable and \fm marked
surfaces induces isomorphisms of flat vector bundles covering
corresponding diffeomorphisms of Teichm\"{u}ller spaces of the
corresponding marked surfaces.}
\end{remark}

However, for a labeled marked surface $(\Sib, \l) = (\Si, P,V,L,
\l)$, which is not \stable or not \fm we need to say a little
more. Namely, let $(\Sib',\l')$ be obtained from $(\Sib,\l)$ by
further labeling points not in $P$ by the zero label $0\in P_\ell$
and choose projective tangent vectors at these new labeled points,
such that $\Sib'$ is both \stable and \fmp. Let $\pi'$ be the
projection map from $\cT_{\Sib'}$ to $\cT_{\Sib}$.

\begin{proposition}\label{nonstabsat}
The connection $\nabla(\Sib',\l')$ has trivial holonomy along the
fibers of the projection map $\pi'$. The connection
$\nabla(\Sib',\l')$ induces a flat connection in the bundle over
$\cT_{\Sib}$ obtained by push forward $\Spofvlamp(\Sib')$ along
$\pi'$ using $\nabla(\Sib',\l')$. If $(\Sib'',\l'')$ is another
\stable and \fm labeled marked surface obtained from $(\Sib,\l)$
in the same way by adding zero-labeled point to $P$, then
iterations of the propagation of vacua isomorphisms given in
Proposition \ref{propvaciso} induces a connection preserving
isomorphism between the corresponding pair of bundles over
$\cT_{\Sib}$.
 \end{proposition}

This proposition follows directly from Proposition \ref{pullbackcon}
and the definition of the flat vector bundle
$\Spofvlamp(\Sib',P')$.

\begin{definition}\label{defnonstabsat}
We define the vector bundle with its flat connection
$(\Spofvlam(\Sib),\nabla(\Sib,\l))$ over $\cT_{\Sib}$ to be
$\pi'_*(\Spofvlamp(\Sib'), \nabla(\Sib',\l'))$ for any \stable and
\fm labeled marked surface $(\Sib',\l')$ obtained from $(\Sib,\l)$
by adding zero-labeled marked points to $P$.
\end{definition}

Suppose now $f : (\Sib_1,\l_1) \ra (\Sib_2,\l_2)$ is a morphism of
labeled marked surfaces and that $(\Sib'_i,\l'_i)$ is obtained as
above from $(\Sib_i,\l_i)$ by adding zero-labeled points and
further that $f': (\Sib'_1,\l'_1) \ra (\Sib'_2,\l'_2)$ is any
morphism of labeled marked surfaces, which induces $f$ when
restricted to $(\Sib_1,\l_1)$. We then have the following result as a direct
consequence of Proposition \ref{propvmorph}.

\begin{proposition}\label{mornonstabsat}
The induced morphism of flat vector bundles $\Spofv(f') :
\Spofvlampone(\Sib'_1) \ra \Spofvlamptwo(\Sib'_2)$ induces a
morphism of flat vector bundles from $\Spofvlamone(\Sib_1)$ to
$\Spofvlamtwo(\Sib_2)$ which only depends on $f$ and which behaves
well under compositions of morphism of labeled marked surface.
\end{proposition}

Let us now collect the thus fare obtained in the following
theorem.

\begin{theorem}
The construction given above gives a functor from the category of
labeled marked surfaces to the category of vector bundles with
flat connections over Teichm\"{u}ller spaces of marked surfaces.
\end{theorem}

The modular functor we seek is now simply just obtained by
composing with the functor which takes covariant constant sections
of vector bundles with connections.

\begin{definition}[The functor $V^{\fg}_{\ell}$]\label{def.main} Let $\ell$ be a positive
integer. Let $P_{\ell}$ be the finite set defined in
(\ref{labelset}) with the involution $\dagger$ as defined
by (\ref{involution}). Let $(\Sib, \l) = (\Si,P,V, L, \l)$ be a
labeled marked surface using the label set $P_{\ell}$. The functor
$V^{\fg}_{\ell}$ is by definition the composite of the functor,
which assigns to $(\Sib, \l)$ the flat vector bundle
$\Spofvlam(\Sib)$ over $\cT_{\Sib}$, and the functor, which takes
covariant constant sections.
\end{definition}

\begin{remark}{\em For a labeled marked surface $(\Sib,\l)$ and a
complex structure $\e C$ on it, we see that Proposition
\ref{Teichpullback} and \ref{Teichpullbackab} give an isomorphism
\[V^{\fg}_{\ell}(\Sib,\l) \cong \Spofvlam(\e C)\otimes (\Vdagab)^{-\frac{1}{2}c_{\cv}}(L)(\e C),\]
since $\cT_{\Sib}$ is contractible. Moreover, if ${\e f} :
(\Sib_1,\l_1)\ra (\Sib_2,\l_2)$ is a morphism of labeled marked
Riemann surfaces, which is realized by a morphism of labeled marked Riemann surfaces
$\Phi: {\e C}_1 \ra {\e C}_2$, such that $\Phi^*(L_2) = L_1$, then
we have the following commutative diagram
$$\begin{CD}
V^{\fg}_{\ell}(\Sib_1,\l_1) @>\cong>> \Spofvlam(\e C_1)\otimes
(\Vdagab)^{-\frac{1}{2}c_{\cv}}(L_1)(\e C_1)\\ @V
{V^{\fg}_{\ell}(\e f)} VV @V\Vdag(\Phi)\otimes
\Vdagab(\Phi)^{-\frac{1}{2}c_{\cv}} VV\\
V^{\fg}_{\ell}(\Sib_2,\l_2) @>\cong>>\Spofvlam(\e C_2)\otimes
(\Vdagab)^{-\frac{1}{2}c_{\cv}}(L_2)(\e C_2).
\end{CD}$$
}\end{remark}

\begin{remark}{\em Let  $(\Sib,\l)$ be a labeled marked surface and
suppose that $(\Sib',\l')$ is obtained from $(\Sib,\l)$ by
labeling further points by $0\in P_\ell$, then by Proposition
\ref{nonstabsat} we get induced an isomorphism
\[V^{\fg}_{\ell}(\Sib,\l) \cong V^{\fg}_{\ell}(\Sib',\l').\]
}\end{remark}

Let  $(\Sib,\l)$ be a labeled marked surface. Let $\goF = (\pi : \mC \ra \mathcal B,\vec s,\vec \eta)$
be a family of stable and saturated Riemann surfaces with formal neighbourhoods over $\Sib$.

\begin{definition}
We define $V^{\fg}_{\ell}(\goF,\l)$ to be the covariant constant sections
of the flat bundles $ \Spofvlam(\goF)\otimes (\Vdagab)^{-\frac{1}{2}c_{\cv}}(L)(\goF)$ over $\mathcal B$.
\end{definition}

Form this definition it is clear that we get an isomorphism
$$I_\goF : V^{\fg}_{\ell}(\goF,\l) \ra V^{\fg}_{\ell}(\Sib,\l).$$

In order for the functor $V^{\fg}_{\ell}$ to be modular, we need
to further construct the disjoint union isomorphism and the
glueing isomorphism and to check that the axioms of a modular
functor is satisfied. First we construct the disjoint union
isomorphism. The glueing isomorphism will be constructed in the
following section.

Let $(\Sib_i,\l_i) = (\Si_i,P_i,V_i,L_i, \l_i)$, $i=1,2$, be two
\stable and \fm labeled marked surfaces and let $(\Sib,\l) =
(\Sib_1,\l_1) \sqcup (\Sib_2,\l_2)$. Let $\l = \l_1 \sqcup \l_2$.
We have that $\cT_{\Sib} = \cT_{\Sib_1}\times\cT_{\Sib_2}$. Let
$\pi_i : \cT_{\Sib}\ra \cT_{\Sib_i}$ be the projection onto the
$i$'th factor. We clearly have that

\begin{proposition}\label{disjointuniiso}
The natural isomorphism $\Hlam\cong\Hlamone \otimes \Hlamtwo$ (for
any ordering $\vec \l, \vec \l_1$ and $\vec \l_2$ of $\l, \l_1$
and $\l_2$ respectively) induces an isomorphism of vector bundles
with connections
\begin{equation}
\Spofvlam(\Si,P) \cong \pi_1^*\Spofvlamone(\Si_1,P_1) \otimes
\pi_2^*\Spofvlamtwo(\Si_2,P_2),\label{isodisj}
\end{equation}
where we use the Lagrangian subspaces to fix the connections in
all three bundles. The isomorphism is compatible with isomorphism
induced by disjoint union of morphism of corresponding labeled
marked surfaces.
\end{proposition}

\begin{remark}{\em These disjoint union isomorphisms are clearly compatible
with the propagation of vacua isomorphisms given in Proposition
\ref{nonstabsat}.}
\end{remark}

Further it is easy to see that
\begin{proposition}\label{disjointuniisoab}
The isomorphism given in Lemma \ref{DJDlemma} induces an
isomorphism of line bundles with connections
\begin{equation}
\Vdagab(\Si) \cong \pi_1^*\Vdagab(\Si_1) \otimes
\pi_2^*\Vdagab(\Si_2),\label{isodisjab}
\end{equation}
where we use the Lagrangian subspaces to fix the connections in
all three bundles. The isomorphism is compatible with isomorphism
induced by disjoint union of morphism of corresponding labeled
marked surfaces. Moreover the preferred sections of the squares of  these bundles
specified by the given Lagrangian subspaces are compatible with
this isomorphism.
\end{proposition}

From this proposition it then follows that we get the corresponding
isomorphism of $-\frac{1}{2}c_{\cv}$-power of these bundles.
Combining this with (\ref{isodisj}) we now get induced a preferred
isomorphism of flat vector bundles
\[\Spofvlam(\Sib) \cong \pi_1^*\Spofvlamone(\Sib_1) \otimes
\pi_2^*\Spofvlamtwo(\Sib_2),\] which intern induces the required
isomorphism of the corresponding vector spaces of covariant
constant sections:
\[V^{\fg}_{\ell}(\Sib,\l) \cong V^{\fg}_{\ell}(\Sib_1,\l_1) \otimes V^{\fg}_{\ell}(\Sib_2,\l_2)\]
which is natural with respect to disjoint union of morphisms.

\section{Sheaf of vacua and gluing.}\label{shofvacandglue}

Let ${\Sib} = (\Si, \{p_-,p_+\}\sqcup P,\{v_-,v_+\}\sqcup V,L)$ be
a marked surface. Let $$c : P(T_{p_-}\Si) \ra P(T_{p_+}\Si)$$ be a
glueing map and $\Si_c$ the glueing of $\Si$ at the ordered pair
$((p_-,v_-),(p_+,v_+))$ with respect to $c$ as described in
section \ref{AxiomsMF}. We shall first assume that $(\Si_c,P)$ is
\stable and \fmp

Let $\goF = (\pi : \mC \ra \mathcal B; s_-,s_+,{\vec s};
\eta_-,\eta_+,\vec{\eta})$ be a family of pointed Riemann surfaces with
formal neighbourhoods on $\Sib$ over a simply-connected base
$\mathcal B$. Let $D$ be the unit disk in the complex plane.
Assume we have holomorphic functions $x_\pm : U_\pm \subset \mC
\ra D$ such that for each $b\in \mathcal B$ we have that
$x_\pm\mid_{\pi^{-1}(b)} : U_\pm\cap \pi^{-1}(b)\ra D$ are local
coordinates for $\pi^{-1}(b)$ centered at $p_\pm$ and further that
$x_\pm = \eta_\pm$ as formal neighbourhoods. Further we assume
that $d_{p_\pm}(x_\pm \mid_{\pi^{-1}(b)})(v_\pm) = 1$ and that $$c
= P(d_{p_+}(x_+\mid_{\pi^{-1}(b)}))^{-1} \circ P({\overline{\cdot}})
\circ P(d_{p_-}(x_-\mid_{\pi^{-1}(b)})) : P(T_{p_-}\Si)
 \ra P(T_{p_+}\Si)$$
 where $P(\overline{\cdot}) : P(T_0 D) \ra P(T_0 D)$ is induced by
 the the real linear map $z \mapsto {\overline z}$.
 Assume that $P \subset \Si - (U_- \cup U_+)$.

Set $\mathcal B_c = \mathcal B \times D$ and $\pi_D : \mathcal B_c
\ra D $ be the projection onto the second factor. Let us now
construct a \stable and \fm family of pointed curves with
formal neighbourhoods $\goF_{c} = (\pi_c : \mathcal C_c \ra
\mathcal B_c,\vec s, \vec \eta)$, in the sense of Definition
\ref{stablefamily} in the Appendix below,
 by applying the glueing construction pointwise over $\mathcal B$ to
 $\goF$:

 Let
$$\mathcal C^1 = \{(z,w,\tau)\in D^{\times 3} \mid zw = \tau\}$$
$$\mathcal C_c^1 = C^1 \times \mathcal B$$ and
\[\mathcal C_c^2 = \{(y,\tau) \in \mathcal C \times D \mid y \in U_\pm \Rightarrow |x_\pm(y)| > |\tau|\}\]
Let then
\[\mathcal C_c = \mathcal C_c^1 \cup_{\phi} \mathcal C_c^2,\]
where
\[\phi : ((U_- - p_-)\times D \cup (U_+ - p_+) \times D) \cap \mathcal C_c^2 \ra \mathcal C_c^1\]
is given by
\[\phi(y,\tau) = \left\{ \begin{array}{ll}
          (x_-(y), \tau/x_-(y),\tau,\pi(y)), & y \in U_- - p_-\\
          (\tau/x_+(y),x_+(y),\tau,\pi(y)), & y \in U_+ - p_+
          \end{array}\right. .\]
One easily checks that $\mathcal C_c$ is a smooth complex manifold
of dimension $\dim(\mathcal B) + 2$ and that we have an obvious
holomorphic projection map $\pi_c : \mathcal C_c \ra \mathcal
B_c$. Let $\goF_{c} = (\pi_c : \mathcal C_c \ra \mathcal B_c,\vec
s, \vec \eta)$.

The fibers over $\mathcal B \times \{0\}$ are nodal
 curves, hence $\goF_c$ is not a family of pointed Riemann
 surfaces with formal neighbourhoods, however it is a family of
pointed stable curves with formal neighbourhoods in the sense
of Definition \ref{stablefamily} in the Appendix below. If $D^* =
D \setminus \{0\}$, then the restricted family $\goF_c|_{\mathcal
B \times D^*}$ is however a family of pointed Riemann
 surfaces with formal neighbourhoods.

Set $\tD = \{ \zeta\in {\mathbb C} | \im(\zeta) > 0\}$. On $\tD$ we
now consider the real coordinates $(r,\theta)$ given by $r(\zeta)
= \exp(-2 \pi \im(\zeta))$ and $\theta(\zeta) = \re(\zeta)$. Let
$(r_\pm,\theta_\pm)$ be $x_\pm$ composed with polar coordinates.

Let ${\tilde {\mathcal B}}_c= \mathcal B \times {\tilde D}$ and
$p_c : {\tilde {\mathcal B}}_c \ra \mathcal B_c^* = \mathcal
B\times D^*$ be given by $p_c(b,\zeta) = (b, \exp(2 \pi i
\zeta))$. Then ${\tilde {\mathcal B}}_c$ is the universal cover of
${\mathcal B}_c$. Let ${\tilde {\mathcal C}}_c = p_c^*{\mathcal
C}_c$, ${\mathcal C}'_c = {\mathcal C}_c\mid_{{\mathcal B}_c^*}$,
$\tgoF_c = p^*_c\goF_c$, ${\tilde \pi}_c : {\tilde {\mathcal C}}_c
\ra {\tilde {\mathcal B}}_c$, ${\tilde \pi}_{\tilde D} : {\tilde
{\mathcal B}}_c \ra {\tilde D}$ and ${\tilde \pi}_{\mathcal B} :
{\tilde {\mathcal B}}_c \ra {\mathcal B}$.

Let
\[V_\pm = \Phi_{\goF}^{-1}(U_\pm)\]
and
\[\tx_\pm : V_\pm \ra {\mathcal B}\times D \]
be given by
\[\tx_\pm = (\pi_{\mathcal B}, x_\pm \circ \Phi_{\goF}).\]
Let us now define a fiber preserving diffeomorphism
\[f : {\tilde {\mathcal C}}_c \ra \Si_c\times {\tilde B}_c\]
by
\[f(y,r,\theta) = \left\{
\begin{array}{ll}
  (\tx_-^{-1}(\pi(y),\chi_r(r_-(y)),\theta_-(y) + \frac{1}{2}\frac{1-r_-(y)}{1 -
  r^{1/2}}\theta), r, \theta)
  & \mbox{if }y\in U_-, 1\geq r_-(y) \geq r^{1/2} \\
  (\tx_+^{-1}(\pi(y),-\chi_r(r_-(y)),-\theta_-(y) - \frac{1}{2}\frac{r_-(y)-r}{ r^{1/2}-r}\theta),
   r, \theta)
  & \mbox{if }y\in U_-, r^{1/2}\geq r_-(y) \geq r,
\end{array}
\right. \] and extend $f$ to all of ${\tilde {\mathcal C}}_c$ by
the map $\Phi_{\goF}^{-1}\times \id$ on $({\mathcal C}-(U_+\cup
U_-))\times \tD$.

Here $\chi_r$ is a smooth family of diffeomorphisms
\[\chi_r : [1,r] \ra [1,-1], \mbox{ }r \in (0,1),\]
with the properties that $\chi_r$ is the identity near $1$,
$\chi_r$ maps $\rho \mapsto -r/\rho$ near $r$ and $\chi_r(r^{1/2})
= 0$ for each $r\in (0,1)$. We will furthermore assume that for
all $\rho \in (0,1)$ we have that
\[\lim_{r\ra 0}\chi_r(\rho) = \rho\mbox{ and } \lim_{r\ra 0} \chi_r(r/\rho) = -\rho,\]
for all $\rho \in (0,1).$

The extra conditions on $\chi_r$ implies that the limit $\lim_{r
\ra 0}q \circ f(\cdot,r,0) : \Si' \ra \Si'$ exists and is equal to
$\id : \Si' \ra \Si'$. We observe that the monodromy
$(f|_{\{b\}\times \pi^{-1}(\zeta+1)}) \circ (f|_{\{b\}\times
\pi^{-1}(\zeta)})^{-1}$ is a Dehn twist in $P(T_{p_-}\Si)$.

Using $f^{-1} : \Si_c\times {\tilde {\mathcal B}}_c \ra {\tilde
{\mathcal C}}_c$ we see that $\tgoF_c$ is a family of \stable and
\fm curves with formal neighbourhoods on $\Si_c$.

The inclusion
\[\Hlam(\mathcal B) \hookrightarrow \bigoplus_{\mu \in P_\ell}
{\mathcal H}_{\mu, \mu^\dagger, \vec \lambda}(\mathcal B),\] given
by
\[|\phi\rangle \mapsto \bigoplus_{\mu \in P_\ell}
 |0_{\mu, \mu^\dagger} \otimes \phi\rangle\]
induces an isomorphism of vector bundles
\begin{equation}
\Vdaglam(\goF_c|_{{\mathcal B}\times \{0\}}) \cong \bigoplus_{\mu
\in P_\ell}\Vdagmulam(\goF).\label{factfam}
\end{equation}
This is the content of Theorem 4.4.9 in \cite{Ue2}.

The \abelian sheaf of vacua construction applied to the family
$\goF_{c}$ gives a holomorphic line bundle $\Vdagab(\goF_{c})$
over $\mathcal B_c$. This follows from Theorem \ref{thm5.2ab} below.
We get an isomorphism of vector bundles
\begin{equation}
\Vdagab(\goF_c)|_{{\mathcal B}\times \{0\}} \cong \Vdagab(\goF)\label{factfamab}
\end{equation}
induced by the isomorphism given in Theorem \ref{thm3.5} below. The
preferred section $s_{\goF_c}(L_c)$ is continuous over $\mathcal
B_c$. Over $\pi^{-1}_D(0)$ it is mapped via the above
isomorphism to the preferred section $s_{\goF}(L)$ of
$\Vdagab(\goF)$. This follows by Theorem \ref{thm6.5} and \ref{thm6.7}.
As discussed before the Lagrangian subspace $L$ determines connections in the bundles
$\oplus_\mu\Vdagmulam(\goF)$ and $\Vdagab(\goF)$.

\begin{proposition}\label{nodalfambidif}
The Lagrangian subspace $L_c \subset H_1(\Si_c,\Z{})$ determines a
unique normalized symmetric bidifferential $\omega_c \in
H^0({\mathcal C}_c\times_{{\mathcal B}_c}{\mathcal C}_c,
\omega_{{\mathcal C}_c\times_{{\mathcal B}_c}{\mathcal
C}_c}(2\Delta))$ specified by formula \eqref{primeform} for
any symplectic basis $(\vec \alpha,\vec \beta)$ of
$H_1(\Si_c,\Z{})$ such that $L_c = \spane \{\beta_i\}$.
\end{proposition}

\proof {Given any symplectic basis $(\vec \alpha,\vec \beta)$ of
$H_1(\Si_c,\Z{})$ such that $L_c = \spane \{\beta_i\}$, formula
\eqref{primeform} defines a normalized symmetric
bidifferential $\tilde\omega_c \in H^0({\tilde{\mathcal
C}}_c\times_{{\tilde {\mathcal B}}_c}{\tilde {\mathcal C}}_c,
\omega_{{\tilde {\mathcal C}}_c\times_{{\tilde {\mathcal B}}_c}{\tilde
{\mathcal C}}_c}(2\Delta))$. The monodromy of the fibration $\pi
: {\mathcal C}_c|_{{\tilde \pi}^{-1}_{\mathcal B}(b)} \ra D$
is the Dehn twist in the Riemann surface $P(T_{p_-}\Sigma)$, hence it
 preserves $L_c$. But then by applying the transformation law
given in (5.2.7) and (5.2.8) \cite{Ue2}, we see that $\tilde
\omega_c$ is invariant under this monodromy and therefore
$\tilde\omega_c$ is the pull back of a unique $$\omega'_c \in
H^0({\mathcal C}'_c\times_{{\mathcal B}_c}{\mathcal C}'_c,
\omega_{{\mathcal C}'_c\times_{{\mathcal B}_c}{\mathcal
C}'_c}(2\Delta)).$$ By the general theory for bidifferentials on
such families, see e.g. \cite{Fa} chapter III, we have that there
is a unique $\omega_c \in H^0({\mathcal C}_c\times_{{\mathcal
B}_c}{\mathcal C}_c, \omega_{{\mathcal C}_c\times_{{\mathcal
B}_c}{\mathcal C}_c}(2\Delta))$ such that $\omega_c\mid_{{\mathcal
C}'_c\times_{{\mathcal B}_c^*}{\mathcal C}'_c} = \omega'_c$.}
\eproof

By Definition \ref{conomega} and Theorem \ref{curvaturecom} we get that
$\omega_c$ determines a projectively flat connection in
$\Vdaglam(\goF_c)|_{{\mathcal B}_c^*}$ and by Theorem \ref{thm4.2a}
a connection in $\Vdagab(\goF_c)|_{{\mathcal
B}_c^*}$.

Let us now recall the conclusion of the glueing constructions on the sheaf of vacua
both in the non-abelian and abelian case applied to the family
$\goF_c$:

The explicit formula \eqref{formalsolution1} and Theorem \ref{thm5.3.4} below give an isomorphism between
sections of $\Vdaglam(\goF_c)|_{{\mathcal B}\times \{0\}}$ and sections of
$\Vdaglam(\tgoF_c)$,
which are covariant constant along the fibers of $\tilde \pi_{{\tD}}$.

The connection in $\Vdagab (\goF_c)|_{{\mathcal B}^*_c}$
determined by $\omega_c$ extends to a connection on all of
$\Vdagab (\goF_c)$, hence we get by parallel transport along the
fibers of $\pi_{{D}}$ an isomorphism between sections of
$\Vdagab(\goF_c)\mid_{{\mathcal B}\times \{0\}}$ and sections of
$\Vdagab (\goF_c)$, which are covariant constant along the fibers
of $\pi_{{D}}$. This follows from Theorem \ref{glueabcova} below and formula \eqref{formalsolution2} gives an explicit
formula for this isomorphism.

By applying the fractional power construction to the line bundle
$\Vdagab (\goF_c)^{\otimes 2}$ with the preferred section
$s_{\goF_c}(L_c)$, we get a line bundle over $\mathcal B_c$, which
we denote $\Vdagab (\goF_c)^{-\frac{1}{2}c_v}(L_c)$. By the very
construction of this bundle we see that $\Vdagab
(\goF_c)^{-\frac{1}{2}c_v}(L_c)\mid_{{\mathcal B}\times \{0\}}$ is
identified with $(\Vdagab)(\goF)^{-\frac{1}{2}c_v}(L)$. We get a
connection in this bundle from its construction and an isomorphism
from sections of $(\Vdagab)(\goF)^{-\frac{1}{2}c_v}(L)$ to
sections of $\Vdagab (\goF_c)^{-\frac{1}{2}c_v}(L_c)$ over
$\mathcal B_c$, which are covariant constant along the
fibers of $\pi_{{D}}$.

\begin{theorem}\label{glueconsistent}
The tensor product of these two glueing constructions gives an isomorphism $I_c(\goF,x_\pm)$ from covariant
constant sections of
$\bigoplus_{\mu\in P_\ell} \Vdagmulam(\goF) \otimes \Vdagab(\goF)^{-\frac{1}{2}c_v}(L)$
over $ {\mathcal B}$ to covariant constant sections of
$\Vdaglam(\tgoF_c) \otimes \Vdagab(\tgoF_c)^{-\frac{1}{2}c_v}(L_c)$
over ${\tilde {\mathcal B}}_c$:
$$I_c(\goF,x_\pm) : \bigoplus_{\mu\in P_\ell} V^{\fg}_{\ell}(\goF,\mu,\mu^\dagger,\l) \ra V^{\fg}_{\ell}(\tgoF_c,\l).$$
\end{theorem}

\proof {From an element of $\bigoplus_{\mu\in P_\ell}
V^{\fg}_{\ell}(\goF,\mu,\mu^\dagger,\l)$, we get a section of
$\Vdaglam(\goF_c) \otimes \Vdagab(\goF_c)^{-\frac{1}{2}c_v}(L_c)
\mid_{{\mathcal B}\times \{0\}}$. By Theorem 5.3 and Remark 5.1 in
\cite{AU1} we see that the covariant derivative of the section of
$\Vdaglam(\tgoF_c) \otimes
\Vdagab(\tgoF_c)^{-\frac{1}{2}c_v}(L_c)$ obtained by glueing
vanishes, since the contributions from the non-abelian factor
exactly cancels the contribution form the $-\frac{1}{2}c_v$-power
of the abelian factor. Hence glueing results in a covariant
constant section of $\Vdaglam(\tgoF_c) \otimes
\Vdagab(\tgoF_c)^{-\frac{1}{2}c_v}(L_c)$ over ${\tilde {\mathcal
B}}_c$. Since the two glueing constructions give isomorphisms, it
is clear that $I_c$ is an isomorphism.} \eproof

Let $\e C^{(i)}$, $i=1,2,$ be two complex structures on $\Si$ and
let $x_\pm^{(i)} : U_\pm^{(i)}\ra D$ be coordinates around $p_\pm$
with $d_{p_\pm}x_\pm^{(i)}(v_\pm) = 1$ such that $c = P((d_{p_+}x^{(i)}_+)^{-1}
\circ P({\overline{\cdot}}) \circ P(d_{p_-}x^{(i)}_-) : P(T_{p_-}\Si)
 \ra P(T_{p_+}\Si)$. Let $\eta_j^{(i)}$ be formal coordinates around $p_j\in C^{(i)}$.

\begin{theorem}\label{indepglue}
For such two pairs $(\e C^{(i)}, x^{(i)}_\pm)$, $i=1,2$,  of complex
structures and holomorphic coordinates on $(\Si, \{p_+,p_-\}\cup
P)$ we have that
\[I_c(\e C^{(1)}, x^{(1)}_\pm) = I_c(\e C^{(2)}, x^{(2)}_\pm).\]
\end{theorem}

This follows straight from Theorem \ref{glueconsistent}, since we clearly have the following

\begin{lemma}\label{contfamily}
There exists a family of pointed Riemann surfaces with formal neighbourhoods
$\goF = (\pi : \mathcal C \ra \mathcal B, \vec s, \vec \eta)$ on
$(\Si,P)$, holomorphic functions $x_\pm : U_\pm \subset \mathcal C \ra D$ and $b_i \in
\mathcal B$ $i=1,2$ such that the following holds
\begin{itemize}
\item The base $\mathcal B$ is simply-connected.
\item Restriction to the fiber
\[(\pi^{-1}(b_i),\vec \eta\mid_{\pi^{-1}(b_i)},x_\pm\mid_{\pi^{-1}(b_i)} :
U_\pm \cap \pi^{-1}(b_i)\ra D)\]
over $b_i$, $i=1,2$, is the same complex structure on $(\Si,P)$ as
\[(\e C^{(i)},\vec \eta^{(i)},x_\pm^{(i)} : U_\pm^{(i)} \ra D)\]
with the same formal coordinates and the same coordinates around
$p_\pm$.
\item For each $b\in \mathcal B$ we have that
$x_\pm\mid_{U_\pm\cap\pi^{-1}(b)} : U_\pm\cap\pi^{-1}(b) \ra D$ are
holomorphic coordinates around $p_\pm\in \pi^{-1}(b)$.
\end{itemize}
\end{lemma}

\begin{definition}\label{IC}
We define the glueing isomorphism
\[I_c=I_c(\Sib,\l) : \bigoplus_{\mu \in P_\ell}
V^{\fg}_{\ell}(\Sib,\mu,\mu^\dagger,\l) \ra V^{\fg}_{\ell}(\Sib_c,\l)\]
to be equal to $I_c(\e C, x_\pm)$ for any pair $(\e C, x_\pm)$ of
a complex structure and holomorphic coordinates on $(\Si,\{p_+,p_-\}\cup
P)$.
\end{definition}

Recall that it is assumed that $(\Si_c,P)$ is \stable and \fmp Let
now $(\Sib',\l')$ be a labeled marked surface obtained
from $(\Sib,\l)$ by labeling further points by $0\in P_\ell$.

\begin{proposition}\label{probofvacglue}
We get the following commutative diagram of isomorphisms:
\begin{equation}
\begin{CD}
\oplus_{\mu\in P_\ell}V^{\fg}_{\ell}(\Sib,\mu,\mu^\dagger,\l) @>I_c(\Sib,\l)>>
V^{\fg}_{\ell}(\Sib_c,\l)\\
@V\cong VV               @V\cong VV\\
\oplus_{\mu\in P_\ell}V^{\fg}_{\ell}(\Sib',\mu,\mu^\dagger,\l') @>I_c(\Sib',\l')>>
V^{\fg}_{\ell}(\Sib_c', \l'),
\end{CD}\label{ccprobofvacglue}
\end{equation}
where the vertical isomorphisms are the ones gotten from Theorem
\ref{nonstabsat}.
\end{proposition}

\proof {This follows since the glueing constructions both in the non-abelian case and the in the abelian case
commutes with the propagation of vacua isomorphism. This is clear from \eqref{formalsolution1}  and
 \eqref{formalsolution2}.}\eproof

\begin{definition}\label{gluenonstab}
In the cases where $(\Sib,\l)$ is not \stable or not \fmc we define
the glueing isomorphism, to be the unique isomorphism
\[I_c = I_c(\Sib,\l): \oplus_{\mu\in P_\ell}V^{\fg}_{\ell}(\Sib,\mu,\mu^\dagger,\l)
\ra V^{\fg}_{\ell}(\Sib_c,\l) \]
which makes the diagram (\ref{ccprobofvacglue}) commutative for any
\stable and \fm labeled marked surface $(\Sib_c', \l')$ obtained
from $(\Sib,\l)$ by labeling further points by $0\in P_\ell$.
\end{definition}

By the naturality of the glueing construction we have that

\begin{proposition}\label{gluemorph}
The glueing isomorphism are compatible with the isomorphisms
induced by glueing morphisms of marked surfaces. That is suppose
$$\e f :  (\Si^1, \{p^1_-,p^1_+\}\sqcup P^1,\{v^1_-,v^1_+\}\sqcup V^1,L^1)
     \ra (\Si^2, \{p^2_-,p^2_+\}\sqcup P^2,\{v^2_-,v^2_+\}\sqcup V^2,L^2)$$
is a morphism of marked surfaces and that there are glueing maps
\[c_i : P(T_{p^i_-}\Si^i) \ra P(T_{p^i_+}\Si^i), \]
such that $(d_{p^1_+}f)^{-1} c_2 d_{p^1_-}f = c_1$, then we get
the
following commutative diagram
\begin{equation*}
\begin{CD}
\oplus_{\mu\in P_\ell}V^{\fg}_{\ell}(\Sib^1,\mu,\mu^\dagger,\l^1) @>I_{c_1}(\Sib^1,\l^1)>>
V^{\fg}_{\ell}(\Sib^1_{c_1},\l^1)\\
@V {V^{\fg}_{\ell}(\e f)} VV               @V {V^{\fg}_{\ell}(\e f)} VV\\
\oplus_{\mu\in P_\ell}V^{\fg}_{\ell}(\Sib^2,\mu,\mu^\dagger,\l^2) @>I_{c_2}(\Sib^2,\l^2)>>
V^{\fg}_{\ell}(\Sib^2_{c_2},\l^2),
\end{CD}
\end{equation*}
for all labelings $\l^i$ of $P^i$ compatible with $\e f$.
\end{proposition}

We summarize the results on the glueing
construction.

\begin{theorem}\label{Glueingiso}
There is an isomorphism $I_c$ from $\oplus_{\mu \in P_\ell}
V^{\fg}_{\ell}(\Sib,\mu,\mu^\dagger,\l)$ to
$V^{\fg}_{\ell}(\Sib_c,\l)$ as specified in Definition \ref{IC} and \ref{gluenonstab},
which is independent of the glueing map $c$ in the following
sense:

If $c_i : P(T_{p_-}\Si) \ra P(T_{p_+}\Si)$, $i=1,2$, are glueing
maps and $f : \Sib_{c_1}\ra \Sib_{c_2}$ is a diffeomorphism as
described in Remark \ref{remarkglue2}, then we have that
\[I_{c_2} = V^{\fg}_{\ell}(f) I_{c_1}\]
as isomorphisms from $\oplus_{\mu \in P_\ell}
V^{\fg}_{\ell}(\Sib,\mu,\mu^\dagger,\l)$ to
$V^{\fg}_{\ell}(\Sib_{c_2},\l)$. Moreover the isomorphisms $I_c$
are compatible with glueing of morphisms of labeled marked
surfaces.
\end{theorem}

Let $(\Sib',\l')$ be another labeled marked surface, which is
\stable and \fmp Let $\Sib''$ be the disjoint union of
$\Sib$ and $\Sib'$. Let $\Sib_c''$ be the glueing of $\Sib''$
using the glueing map $c$. Then we clearly have that $\Sib_c'' = \Sib_c \sqcup
\Sib'$. It is trivial to check that

\begin{proposition}\label{duiglue}
The glueing isomorphism is compatible with the disjoint union
isomorphism, namely the following diagram is commutative
\begin{equation*}
\begin{CD}
\bigoplus_{\mu\in P_\ell}V^{\fg}_{\ell}(\Sib'',\mu,\mu^\dagger,\l,\l')
@>>>
\bigoplus_{\mu\in P_\ell}V^{\fg}_{\ell}(\Sib,\mu,\mu^\dagger,\l,)\otimes
V^{\fg}_{\ell}(\Sib'\l')\\
@VI_c(\Sib'',\l,\l') VV               @VI_c(\Sib,\l)\otimes \id VV\\
V^{\fg}_{\ell}(\Sib_c'',\l,\l')
@>>> V^{\fg}_{\ell}(\Sib_c,\l)\otimes V^{\fg}_{\ell}(\Sib',\l') .
\end{CD}
\end{equation*}
\end{proposition}

Let now ${\Sib} = (\Si, \{p^{(1)}_-,p^{(1)}_+,p^{(2)}_-,p^{(2)}_+\}
\sqcup P,\{v^{(1)}_-,v^{(1)}_+,v^{(2)}_-,v^{(2)}_+\}\sqcup V,L)$
be a marked surface. Let $$c^{(i)} : P(T_{p^{(i)}_-}\Si) \ra P(T_{p^{(i)}_+}\Si)$$ be
glueing maps and $\Si_{c^{(i)}}$ the glueing of $\Si$ at the ordered pair
$((p^{(i)}_-,v^{(i)}_-),(p^{(i)}_+,v^{(i)}_+))$
with respect to $c^{(i)}$. Let $\Si_{c^{(12)}}$ be the glueing
with respect to $c^{(12)} = c^{(1)}\sqcup c^{(2)}$.

\begin{theorem}\label{glueglue}
The glueing isomorphisms commute, meaning the following diagram is commutative
\[\begin{CD}
\bigoplus_{\mu_1,\mu_2\in P_\ell}V^{\fg}_{\ell}(\Sib,\mu_1,\mu_1^\dagger,\mu_2,\mu_2^\dagger,\l)
@>\bigoplus_{\mu_2\in P_\ell}I_{c^{(1)}}(\Sib,\mu_2,\mu_2^\dagger,\l)>>
\bigoplus_{\mu_2\in P_\ell}V^{\fg}_{\ell}(\Sib_{c^{(1)}},\mu_2,\mu_2^\dagger,\l)\\
@V\bigoplus_{\mu_1\in P_\ell}I_{c^{(2)}}(\Sib,\mu_1,\mu_1^\dagger,\l)VV
@V I_{c^{(2)}}(\Sib_{c^{(1)}},\l)VV\\
\bigoplus_{\mu_1\in P_\ell}V^{\fg}_{\ell}(\Sib_{c^{(2)}},\mu_1,\mu_1^\dagger,\l)
 @>I_{c^{(1)}}(\Sib_{c^{(2)}},\l)>>
V^{\fg}_{\ell}(\Sib_{c^{(12)}}, \l).
\end{CD}\]
\end{theorem}

\proof

Choose a complex structure on $\Si$ and let $\e C$ denote the resulting marked
Riemann surface. The complex structure $\e C$ gives a point in the Teichm\"{u}ller
space $\cT_{(\Si,\{p^{(1)}_-,p^{(1)}_+,p^{(2)}_-,p^{(2)}_+\}
\sqcup P)}$.
Choose centered coordinates
$x^{(i)}_\pm : U_\pm \ra D$ around $p^{(i)}_\pm$ with $d_{p^{(i)}_\pm}x^{(i)}_\pm(v^{(i)}_\pm)
= 1$ and such that
$c^{(i)} = P(d_{p^{(i)}_+}x^{(i)}_+)^{-1} \circ P({\overline{\cdot}}) \circ P(d_{p^{(i)}_-}x^{(i)}_-)
: P(T_{p^{(i)}_-}\Si)
 \ra P(T_{p^{(i)}_+}\Si)$.
The following construction of a smooth $3$-dimensional complex manifold $\mathcal C$
with a holomorphic map $\pi : \mathcal C \ra D\times D$ is the main ingredient in this proof :

Let
$$\mathcal C_1 = \{(z^{(i)}_1,w^{(i)}_2,\tau^{(1)},\tau^{(2)})\in (\C{2}\sqcup \C{2})\times D\times D \mid z^{(i)}w^{(i)} =
\tau^{(i)}, |z^{(i)}| < 1, |w^{(i)}| < 1, |\tau^{(i)}| <
1, i=1,2\}$$
and
\[\mathcal C_2 = \{(y,\tau) \in \Si\times D\times D \mid y \in U^{(i)}_\pm
\Rightarrow |x^{(i)}_\pm(y)| > |\tau^{(i)}|\}\]
Let then
\[\mathcal C = \mathcal C_1 \cup_{\phi} \mathcal C_2,\]
where
\[\phi : ( (U^{(i)}_- - p^{(i)}_-)\times D\times D \cup (U^{(i)}_+ - p^{(i)}_+)
\times D\times D )\cap \mathcal C_2 \ra \mathcal C_1\]
is given by
\[\phi(y,\tau) = \left\{ \begin{array}{ll}
          (x^{(i)}_-(y), \tau^{(i)}/x^{(i)}_-(y),\tau^{(1)},\tau^{(2)}),
          & y \in U^{(i)}_- - p^{(i)}_-\\
          (\tau^{(i)}/x^{(i)}_+(y),x^{(i)}_+(y),\tau^{(1)},\tau^{(2)}), & y \in U^{(i)}_+ - p^{(i)}_+
          \end{array}\right. .\]
One easily checks that $\mathcal C$ is a smooth complex manifold
of dimension $3$ and that we have an obvious holomorphic
projection map $\pi : \mathcal C \ra D\times D$. Choose formal
neighbourhoods $\vec \eta$ for the points $P$. We thus get a
family of stable pointed curves with formal
neighbourhoods $\goF_{c^{(12)}}$ over $D\times D$ obtained by
applying the glueing construction at the two pairs $(p_-^{(1)},
p_+^{(1)})$ and $(p_-^{(2)}, p_+^{(2)})$. Let $p = p_1\times p_2 :
\tD\times \tD \ra D^*\times D^*$, be the projection and let
$\tilde \goF_{c^{(12)}} = p^*(\goF_{c^{(12)}}|_{D^*\times D^*})$.
We denote the two projections onto the first factor by $p^{(1)}
: \tD \times \tD \ra \tD $ and onto the second by $p^{(2)}
: \tD \times \tD \ra \tD.$

Further let $\tilde \goF_{c^{(1)}}$ be the pull back under $p_1$
of the normalization of $\goF_{c^{(12)}}\mid_{D^*\times\{0\}}$ at
$[p_-^{(2)}] = [ p_+^{(2)}]$ and $\tilde \goF_{c^{(2)}}$ be the
pull back under $p_2$ of the normalization of
$\goF_{c^{(12)}}\mid_{{\{0\}}\times D^*}$ at $[p_-^{(1)}] = [
p_+^{(1)}]$. Let $\goX = (\e
C,\{p^{(1)}_-,p^{(1)}_+,p^{(2)}_-,p^{(2)}_+\} \sqcup P,
\{x^{(1)}_-,x^{(1)}_+,x^{(2)}_-,x^{(2)}_+\} \sqcup \vec \eta)$.

 The glueing construction in the non-abelian case applied to
$\goF_{c^{(1)}}$ (respectively to $\goF_{c^{(2)}}$) and then to
$\goF_{c^{(12)}}$ results in two two-variable versions of
\eqref{formalsolution1}
 for any element of ${\mathcal V}_{\mu_1,\mu_1^\dagger,
\mu_2,\mu_2^\dagger, \nu}^{\dagger}(\goX)$. Explicitly for an
element $\langle \Psi | \in {\mathcal V}_{\mu_1,\mu_1^\dagger,
\mu_2,\mu_2^\dagger, \nu}^{\dagger}(\goX)$ we get a section
$\langle \widehat{\Psi}^{(1)}|$ of $\mathcal{V}_{\mu_2,
\mu_2^\dagger,\nu}^\dagger (\goF_{c^{(1)}})$ given by $$ \langle
\widehat{\Psi}^{(1)}| \Phi^{(1)}\rangle = \sum_{d=0}^\infty
\left\{ \sum_{i=1}^{m_d} \langle \Psi|v^{(1)}_i(d)\otimes
v^i_{(1)}(d) \otimes \Phi^{(1)}\rangle
\right\}(\tau^{(1)})^{\Delta_{\mu_1} +d}, $$ where $\{
v^{(1)}_1(d), \ldots, v^{(1)}_{m_d}(d)\}$ is a basis of
$\mathcal{H}_{\mu_1}(d)$ and $\{ v^1_{(1)}(d), \ldots,
v^{m_d}_{(1)}(d)\}$ is the dual basis of
$\mathcal{H}_{\mu_1^\dagger}(d)$ and $| \Phi^{(1)}\rangle$ is any
section of $\mathcal H_{\mu_2,\mu_2^\dagger,\nu}$. The holomorphic
section $\langle \widehat{\Psi}^{(1)}|$ is covariant constant
along the fibers of $p_1$. Next we construct a holomophic section
$\langle \widehat{\Psi}^{(12)}|$ of $\mathcal{V}_\nu^\dagger
(\goF_{c^{(12)}})$ determined by $$ \langle \widehat{\Psi}^{(12)}|
\Phi\rangle = \sum_{e=0}^\infty \left\{ \sum_{j=1}^{m_e} \langle
\widehat{\Psi}^{(1)}|v_j^{(2)}(e)\otimes v^j_{(2)}(e) \otimes
\Phi\rangle \right\}(\tau^{(2)})^{\Delta_{\mu_2} +e}, $$ where $\{
v^{(2)}_1(e), \ldots, v^{(2)}_{m_e}(e)\}$ is a basis of
$\mathcal{H}_{\mu_2}(e)$ and $\{ v^1_{(2)}(e), \ldots,
v^{m_e}_{(2)}(e)\}$ is the dual basis of
$\mathcal{H}_{\mu_2^\dagger}(e)$ and $| \Phi\rangle$ is any element in
$\mathcal H_{\nu}$. This section is covariant constant along the fibers of
$p^{(1)}$.

Similarly, starting form $\langle \Psi | \in {\mathcal
V}_{\mu_1,\mu_1^\dagger, \mu_2,\mu_2^\dagger,
\nu}^{\dagger}(\goX)$ we can construct a
holomorphic section of
$\mathcal{V}_{\mu_1,\mu_1^\dagger,\nu}(\goF_{c^{(2)}})$, covariant constant
along the fibers of $p_2$, which is
determined by $$ \langle \widehat{\Psi}^{(2)}| \Phi\rangle =
\sum_{e=0}^\infty \left\{ \sum_{j=1}^{m_e} \langle
\Psi^{(2)}|v_j^{(2)}(e)\otimes v^j_{(2)}(e) \otimes \Phi\rangle
\right\}(\tau^{(2)})^{\Delta_{\mu_2} +e}. $$

Then, starting form $\langle \widehat{\Psi}^{(2)}|$ we can
construct a holomorphic section $\langle \widehat{\Psi}^{(21)}|$
of $\mathcal{V}_\nu^\dagger (\goF_{c^{(12)}})$ which is covariant
constant along the fibers of $p^{(2)}$, and given by $$ \langle
\widehat{\Psi}^{(21)}| \Phi\rangle = \sum_{d=0}^\infty \left\{
\sum_{i=1}^{m_d} \langle \widehat{\Psi}^{(2)}|v^{(1)}_i(d)\otimes
v^i_{(1)}(d) \otimes \Phi \rangle
\right\}(\tau^{(1)})^{\Delta_{\mu_1} +d}. $$ Now $\langle
\widehat{\Psi}^{(12)}|$ and $\langle \widehat{\Psi}^{(21)}|$ is
given by the same power series $$
 \sum_{d=0}^\infty
\sum_{i=1}^{m_d} \sum_{j=1}^{m_e} \langle \Psi|v^{(1)}_i(d)
\otimes v^i_{(1)}(d) \otimes v^{(2)}_j(e) \otimes v^j_{(2)}(e)
\otimes \Phi \rangle (\tau^{(1)})^{\Delta_{\mu_1} +d}
(\tau^{(2)})^{\Delta_{\mu_2} +e}. $$ Hence, the two sections
$\langle \widehat{\Psi}^{(12)}|$ and $\langle
\widehat{\Psi}^{(21)}|$ of $\mathcal{V}_\nu^\dagger
(\goF_{c^{(12)}})$ coincide, and in fact they are covariant
constant.

From this we in particular see that the connection in ${\mathcal
V}_{ \nu}^{\dagger}(\tilde \goF_{c^{(12)}})$ is flat. By applying
the same argument to the abelian theory, we get that the glueing
construction is also independent of the order of the glueing in
the abelian case and the connection in $\Vdagab(\goF_{c^{(12)}})$
is flat over $D\times D$. The theorem now follows.

 \eproof

\section{Verification of the axioms}\label{verification}

It is now straight forward to check the axioms of a modular
functor given the results obtained in the previous sections.

\begin{theorem}\label{Main}
The functor $V^{\fg}_{\ell}$ from the category of labeled marked
surfaces to the category finite dimensional vector spaces is a
modular functor.
\end{theorem}

\proof In order to check axiom {\em MF1}, we only need to check
that the disjoint union isomorphisms
satisfies associativity, but
this follows from associativity of the isomorphisms between the
corresponding ${\mathcal H}^\dagger$'s and ${\mathcal F}$'s.

We have that the glueing isomorphism $I_c$ from Theorem
\ref{Glueingiso} is compatible with
\begin{itemize}
  \item The disjoint union isomorphisms: Proposition \ref{duiglue}.
  \item The glueing isomorphisms them self, i.e. the glueing
  isomorphisms should commute: Theorem \ref{glueglue}
\end{itemize}
Hence axiom {\em MF2} is checked.

Axiom {\em MF3} is trivial, since we define $V^{\fg}_{\ell}(\emptyset)
= \C{}$. Axiom {\em MF4} and {\em MF5} follows from Corollary 3.5.2 (1) and (2) in \cite{Ue2}.
\eproof

\section{Appendix. Families of stable curves and glueing}\label{App}

In order to define the functor $V^{\fg}_{\ell}$, we only need to
refer to Riemann Surface and families of such which form smooth
complex manifolds. However, in defining the glueing morphisms we
also need to consider the so called {\em stable nodal curves}.
These are one dimensional algebraic sub-varieties (algebraic curves) of complex
projective space, such that the singularities
are locally analytically isomorphic to a neighbourhood of the origin
of $xy=0$. A neighbourhood of a node is obtained by patching
together at the origins of two small disks $D_1=\{\; x; |\; |x|<
\varepsilon_1\;\}$ and $D_2=\{\; y; |\; |y|< \varepsilon_2\;\}$.
By reversing the process, from a neighbourhood of a node we obtain
two disconnected small disks. This process is called
normalization or desingularization of the node. Thus, by a
normalization of a nodal curve $C$ we obtain a compact Riemann
surface $\widetilde{C}$ and holomorphic map $\nu : \widetilde{C}
\rightarrow C$ such that for a node $P$, the inverse image
$\nu^{-1}(P)$ consists of two distinct points $P_+$ and $P_-$. A
curve which has only nodes as singularities is called a nodal
curve.

\subsection{Stable curves}
We begin by introducing the notion of stable curves.

\begin{definition}\label{stablecurve}
The data ${\mathfrak X} = (C;\, Q_1, Q_2, \ldots , Q_N)$
consisting of
 an algebraic curve $C$ and points $Q_1, \ldots, Q_N$ on $C$
are called an {\it ($N$-)pointed stable curve,} if the following conditions are satisfied.
\begin{enumerate}
\item[(1)] \quad  The curve $C$ is a nodal curve.
\item[(2)] \quad
 $Q_1,Q_2, \ldots, Q_N$ are non-singular points of the curve $C$.
\item[(3)] \quad
  If an irreducible component $C_i$ is a Riemann
sphere ${\mathbf P}^1$ (resp. a rational curve with one double point,
resp. an elliptic curve), the sum of the number of
intersection points of $C_i$ and other components and the number of $Q_j$'s
on $C_i$ is at least three (resp. one, resp. one).
\item[(4)] \quad
 $\dim_\bolc H^1(C,\mathcal O_C) = g$.
 \end{enumerate}
\end{definition}

Note that  the condition (3) is equivalent to saying that
the Euler characteristic of each component of the complement of
the points $Q_j$ on it and the nodes is negative.

A pointed stable curve with formal neighbourhoods is defined
in an analogous way as a pointed Riemann surface with formal
neighbourhoods (see Definition \ref{v1order}). Also we
can define a family of pointed stable curves with formal
neighbourhoods.
\begin{definition}\label{stablefamily}
The data $\goF = (\pi : \mathcal{C} \rightarrow \mathcal{B};
s_1,s_2,\ldots, s_N; {\eta}_1 , {\eta}_2 , \ldots, {\eta}_N )$ is
called a family of ($N$-)pointed stable curves of genus $g$ with
formal neighbourhoods, if the following conditions are satisfied.
\begin{enumerate}
\item[(1)] Both  $\mathcal{C}$ and $\mathcal{B}$ are connected complex manifolds,
$\pi : \mathcal{C} \rightarrow
\mathcal{B}$ is a proper flat holomorphic map and $s_1,s_2, \ldots,s_N$ are
holomorphic sections of $\pi$.
\item[(2)] For each point $b \in \mathcal{B}$ the data $(\mathcal{C}_b :=\pi^{-1}(b); s_1(b),
 s_2(b), \ldots,
s_N(b))$ is an $N$-pointed stable curve of genus $g$.
\item[(3)] For each $j$, ${\eta}_j $ is an $\mathcal O_{\mathcal{B}}$-algebra isomorphism
$$
    {\eta}_j  : \widehat{\mathcal O}_{_{\mathcal{C}}/s_j} = \varprojlim_{n \to \infty}
\mathcal O_{\mathcal{C}}/I_{j}^{n} \simeq \mathcal O_{\mathcal{B}}[[\xi]],
$$
where $I_{j}$ is the defining ideal of $s_j(\mathcal{B})$ in $Y$.
\end{enumerate}
\end{definition}

The only families of pointed stable curves we need in this paper
are all constructed explicitly
from families of pointed Riemann surfaces via the glueing
process discussed in the begining of section \ref{shofvacandglue}.

\subsection{Sheaf of vacua for a family of pointed stable curves}

For a family of stable curves with formal neighbourhoods $\goF = (
\pi : \mathcal{C} \rightarrow \mathcal{B} ; \vec{s}; \vec{\eta})$
we can define the sheaf of vacua $\Vdaglam(\goF)$ and the sheaf of
covacua $\Vlam(\goF)$ just as in Definition \ref{ShofVdef}. If the
family contains smooth curves (Riemann surfaces), then we have
 the connection given by Definition \ref{conomega} on the
complement of the locus of nodal curves and the connection has a
regular singularity along this locus (see \S5.3 of \cite{Ue2}).

Using this connection with regular singularities it is shown in
\cite{Ue2}, that
 Theorem \ref{localfreeness} is also valid for a family of stable curves.

\begin{theorem}\label{stabllocalfreeness}
For a stable family of stable curves with formal neighbourhoods
$\goF = ( \pi : \mathcal{C} \rightarrow \mathcal{B} ; \vec{s};
\vec{\eta})$, the sheaves $ \Vdaglam(\goF)$ and $\Vlam(\goF)$ are
locally free sheaves of $\Ob$-modules of finite rank over
$\mathcal B$. They are duals to of each other.
\end{theorem}

We include this theorem here for completeness. We do not need this
result for the constructions in this paper.

\subsection{The sheaf of abelian vacua associated to families of pointed stable
 curves}

For a family of stable curves with formal neighbourhoods $\goF = (
\pi : \mathcal{C} \rightarrow \mathcal{B} ; \vec{s}; \vec{\eta})$
we can define the sheaf of abelian vacua $\Vdagab(\goF)$
just as in Definition \ref{dfn4.1}.
We have following result.

\begin{theorem}[{[2, Theorem 5.2]}]
\label{thm5.2ab}
The abelian vacua construction applied to $\goF$ gives a holomorphic
line bundle $\Vdagab(\goF)$
over $\mathcal B$.
\end{theorem}

Next let us consider a nodal curve $C$ with node $P$. Let
$\widetilde{C}$ be the Riemann surface obtained by resolving the
singularity at $P$ and let $\pi : \widetilde{C} \rightarrow C$ be
the natural holomorphic mapping. Then $\pi^{-1}(P)$ consists of two
points $P_+$ and $P_-$. Let $$
\gX=(C;q_1,\ldots,q_N;\xi_1,\ldots,\xi_N) $$ be a pointed nodal
curve with formal neighbourhoods and we let $$ \widetilde{\gX} =
(\widetilde{C};P_+,P_-,q_1,\ldots,q_N;z,w, \xi_1,\ldots,\xi_N) $$ be
the associated pointed Riemann surface with formal neighbourhoods.
Define an element $|0_{+,-}\rangle \in \cF\otimes \cF$ by
\begin{equation}
\label{0+-} |0_{+,-}\rangle = |0\rangle \otimes |-1\rangle
 - |-1\rangle \otimes |0 \rangle .
\end{equation}
The natural  inclusion
\begin{eqnarray*}
 \cF_N & \hookrightarrow & \cF_{N+2} \\
 |u \rangle & \mapsto & |0_{+,-}\rangle  \otimes |u \rangle
\end{eqnarray*}
defines a natural linear mapping
$$
\iota_{+,-}^* : \cFd_{N+2}  \rightarrow   \cFd_N .
$$
\begin{theorem}[{[2, Theorem 3.5]}]
\label{thm3.5} The natural mapping $\iota^*_{+,-}$  induces  a
natural isomorphism
$$
\Vdagab(\widetilde{\gX}) \cong \Vdagab(\gX).
$$
\end{theorem}

Let us now consider the case where we have one marked point on $C$
$$
\gX=(C;Q;\xi)
$$
and we let
$$
\widetilde{\gX} = (\widetilde{C};P_+,P_-,Q;z,w, \xi)
$$
be the associated 3-pointed curve with formal neighbourhoods.
We can define $\langle \omega(\gX,\{\alpha, \beta\})|$ similar to
the non-singular case, by choosing a basis
$$(\vec{\alpha}, \vec{\beta})=(\alpha_1, \ldots, \alpha_{g-1},\alpha_g, \beta_1,
\ldots, \beta_{g-1})$$ of $H_1(C,\bZ)$, in such a way that
$\alpha_1$, $\alpha_2, \ldots,\alpha_{g-1}$ and $\beta_1$, $\beta_2,
\ldots, \beta_{g-1}$ is the image of a symplectic basis of
$H_1(\widetilde{C},\bZ)$ under natural map to $H_1(C,\bZ)$ and
$\alpha_g$ corresponds to the invariant cycle of a flat deformation
of the curve $C$. Then we can choose a basis $\{\omega_1, \ldots,
\omega_{g-1}, \omega_g, \omega_{g+1}, \omega_{g+2}, \ldots \}$ of
$H^0(C, \omega(*Q))$ such that $\{\pi^*\omega_1, \ldots,
\pi^*\omega_{g-1}, \pi^*\omega_{g+1}, \pi^*\omega_{g+2}, \ldots \}$
is a normalized basis of $H^0(\widetilde{C},
\omega_{\widetilde{C}}(*Q))$ as in (\ref{betaone2}), (\ref{omegaQ1})
and (\ref{omegaQ}) where we put $$ \pi^*\omega_{g+n} =
\omega_Q^{(n)}, \quad n=1,2, \ldots, $$ and $\pi^*\omega_g$ is a
meromorphic one-form on $\widetilde{C}$ which has poles of order one
at $P_+$ and $P_-$ with residue $-1$ and 1, respectively, is
holomorphic outside $P_\pm$ and $$ \int_{P_+}^{P_-}\pi^*\omega_g =
1. $$
Then put $$ \langle \omega(\gX,\{\alpha, \beta\})| = \langle
\cdots \wedge e(\omega_{m} )\wedge \cdots \wedge e(\omega_2)\wedge
e(\omega_1) |. $$ The proof of Lemma 3.1 of \cite{AU1} applies also
in this case and shows that $\langle \omega(\gX,\{\alpha, \beta\})|$
is an element of $\Vdagab(\gX)$. Let
$$ \widehat{\gX} = (\widetilde{C}, Q;
\xi) . $$ Then, by applying Theorem \ref{thm3.4} at the points
$P_\pm$ we have a canonical isomorphism
$$ \iota^* :
\Vdagab(\widetilde{\gX}) \cong \Vdagab(\widehat{\gX}) .
$$
\begin{theorem}[{[2, Theorem 6.5]}]
\label{thm6.5} Under the above assumptions and notation we have that
$$
 \iota^* \circ  (\iota_{+,-}^*)^{-1}  (\langle \omega(\gX,\{\alpha, \beta\}) |) =
 (-1)^g  \langle \omega(\widehat{\gX},\{\widehat{\alpha}, \widehat{\beta}\})|
$$
where $\{\widehat{\alpha}, \widehat{\beta}\}= \{\alpha_1, \ldots,
\alpha_{g-1}, \beta_1, \ldots, \beta_{g-1}\}$.
\end{theorem}

Let $\gF=(\pi : \cC\rightarrow \cB, s_1, \ldots, s_N, \xi_1, \ldots,
\xi_N)$ be a family of $N$-pointed Riemann surfaces. For any point
$b \in \cB$ there exists an open neighbourhood $U_b$ such that
$\pi^{-1}(U_b)$ is topologically  trivial so that we can choose
smoothly varying symplectic bases
$$
(\alpha_1(t), \ldots, \alpha_g(t), \beta_1(t), \ldots, \beta_g(t)),
\quad t \in \cB.
$$
Then we can define, using propagation of vacua,
$$\langle \omega(\gX_t, \{\alpha(t), \beta(t)\})| \in \Vdagab(\gF)_t$$
where $\gX_t=(\pi^{-1}(t), s_1(t), \xi_1)$.

\begin{theorem}[{[2, Theorem 6.6]}]
\label{thm6.6} The section $\langle \omega(\gX_t, \{\alpha(t),
\beta(t)\})|$ is a holomorphic section of $\Vdagab(\gF)$ over $U_b$.
\end{theorem}

\subsection{The glueing of vacua construction}
Let us briefly recall the glueing of vacua. First let us consider
the non-abelian case. We use freely the notation in
\S\ref{shofvacandglue}. Let $\goF = (\pi : \mC \ra \mathcal B;
s_-,s_+,{\vec s}; \eta_-,\eta_+,\vec{\eta})$ be a family of
pointed Riemann surfaces with formal neighbourhoods on $\Sib$ over
a simply-connected base $\mathcal B$ and we let $\goF_{c} = (\pi_c
: \mathcal C_c \ra \mathcal B_c,\vec s, \vec \eta)$ be a \stable
and \fm family of $N$-pointed curves with formal neighbourhoods
obtained by applying the glueing construction pointwise over
$\mathcal B$ to $\goF$ where $\mathcal{B}_c = \mathcal{B} \times
D$ with the unit disk $D$.

By the isomorphism \eqref{factfam} holomorphic sections of
$\oplus_{\mu}\Vdagmulam(\goF)$ over $\mathcal {B}$ may be regarded
as holomorphic sections of $\Vdaglam(\goF_c|_{{\mathcal B}\times
\{0\}})$ over $\mathcal{B} \times \{0\}$. The glueing of vacua is
an isomorphism from holomorphic section of
$\oplus_{\mu}\Vdagmulam(\goF)$ over $\mathcal{B}$ to a holomorphic
sections of $\Vdaglam({\tilde \goF}_c)$ over ${\tilde {\mathcal
B}}_c$, which is covariant constant along the direction of $\tD$.
The construction is as follows.

By Lemma \ref{L2.2.12} we can choose a basis $\{ v_1(d), \ldots,
v_{m_d}(d)\}$ of $\mathcal{H}_\mu(d)$ and the dual basis $\{
v^1(d), \ldots, v^{m_d}(d)\}$ of $\mathcal{H}_{\mu^\dagger}(d)$
such that
\begin{equation} \label{nonabelianpair}
(v^j(d)|v_k(d)) = \delta^j_k .
\end{equation}
For a holomorphic section $ \langle \Psi| \in \Vdagmulam(\goF)$
over $\mathcal{B}$ we define a formal series $\langle
\widehat{\Psi}|$ by
\begin{equation} \label{formalsolution1}
\langle \widehat{\Psi}| \Phi\rangle = \sum_{d=0}^\infty \left\{
\sum_{i=1}^{m_d} \langle \Psi|v_i(d)\otimes v^i(d) \otimes
\Phi\rangle \right\}\tau^{\Delta_\mu +d},
\end{equation}
for all $| \Phi\rangle \in \mathcal H_{\vec \lambda}.$ Here the
fractional power $\tau^{\Delta_\mu +d}$ is clearly well defined on
$\tD$. This formal power series converges and defines in fact a
holomorphic section of $\Vdaglam(\goF_c)$. Namely, we have the
following theorem.
\begin{theorem}[{[19, Theorem 5.3.4]}] \label{thm5.3.4}
The formal power series $\langle \widehat{\Psi}|$ is a formal solution of
the differential equation
\begin{equation}\label{fuchsian}
\left(\tau \frac{d}{d\tau} - T[\vec{l}] + a(\vec{l})\right)
\langle \widehat{\Psi}|=0,
\end{equation}
of Fuchsian type where $$ \vec{l}=\left(
l_1(\xi_1)\frac{d}{d\xi_1}, \ldots, l_N(\xi_N) \frac{d}{d\xi_N}
\right) $$ is an $N$-tuple of formal vector fields such that
$\theta(l) =\tau \dfrac{d}{d\tau}$. Moreover, $\langle
\widehat{\Psi}|$ converges and defines a holomorphic section of
$\Vdaglam(\tilde \goF_c)$ over ${\tilde {\mathcal B}}_c$.
\end{theorem}

The differential equation \eqref{fuchsian} gives the holomorphic
connection along the direction of $\tD$. Hence, the above formal
power series is a holomorphic section of $\Vdaglam(\goF_c)$ which
is covariant constant along the direction of $\tD$.

Let us now discuss the glueing in the abelain case.

Recall we have a perfect pairing
$$
\{\phantom{X}|\phantom{X}\}_+ : \mathcal{F}_d(p) \times \mathcal{F}_d(-p-1) .
$$
Let $\{v_i(d,p)\}_{i=1,\ldots,m_d}$ be a basis of $\cF_d(p)$
for any $p \in \bZ$ and
$\{v^i(d,p)\}_{i=1,\ldots,m_d}$ be the dual basis of $\cF_d(-p-1)$
with respect to the pairing $\{\phantom{X}|\phantom{X}\}_+$.

For a holomorphic section $\langle \psi|$ of
$\Vdagab(\goF)$
define $\langle \widetilde{\psi} |$ by
\begin{equation}
\label{formalsolution2} \langle \widetilde{\psi} | u \rangle =
\sum_{p \in \bZ}\Bigl\{ \sum_{d=0}^\infty
\sum_{i=1}^{m_d}(-1)^{p+d} \langle \psi |v_i(d,p)\otimes
v^i(d,-p-1)\otimes u\rangle \Bigr\}\tau^{d+p(p+1)/2}.
\end{equation}
Then the formal power series converges and defines a holomorphic section of
$\langle \widetilde{\psi}|  \in \Vdagab(\goF_c)$. Moreover this section is
covariant constant along the directions of $D$.

\begin{theorem}
\label{glueabcova} The glueing construction gives an
isomorphism between sections of $\Vdagab(\goF)$ and sections of $\Vdagab(\goF_c)$
which are covariant constant along the directions of $D$.
\end{theorem}

This theorem follows directly form Theorem 5.3 in \cite{AU1}.
Finally, we analyze the preferred section for the families.
Suppose
we have a continuous basis $(\alpha_i(t),\beta_i(t))$ of
$H_1(\pi^{-1}(t),\bZ)$, $t\in (0,1)\subset D$, such that we get a
well defined limit as $t$ goes to zero, which gives a symplectic
basis, say $(\alpha_1(0),\ldots,
\alpha_{g-1}(0),\alpha_{g}(0),\beta_{1}(0),\ldots,\beta_{g-1}(0))$
of $H_1(C_0,\bZ)$ as described above for nodal curves and
$\beta_g(0) = 0$. Let $\gX_t=(\pi^{-1}(t), s_1(t), \ldots s_N(t),
 \xi_1, \ldots \xi_N)$.

\begin{theorem}[{[2, Theorem 6.7]}]
\label{thm6.7} We have that
$$\langle \omega(\gX_0, \{\alpha(0), \beta(0)\})| =
\lim_{t\to 0} \langle \omega(\gX_t, \{\alpha(t), \beta(t)\})|.$$
\end{theorem}

\end{document}